\documentclass{lmcs}


\input{macros-nuc-lmcs}
\graphicspath{{PIC/}}
\begin{document}

\title[Nucleus I]{\raisebox{-3.5ex}{Nucleus I: Adjunction spectra}\\[1.5ex]
in recommender systems and descent}

\author[D.~Pavlovic]{Dusko Pavlovic}[a]
\author{Dominic J.\,D.\,Hughes}[b,c]

\address{University of Hawaii, Honolulu,  HI}
\thanks{$^{a}$Supported by NSF and AFOSR.}
\address{Apple Inc., Cupertino, CA}
\address{UC Berkeley, Berkeley, CA}
\thanks{$^{c}$Visiting scholar at Logic Group. Thanks to Wes Holliday and Dana Scott}

\begin{abstract}
Recommender systems build user profiles using concept analysis of usage matrices. The concepts are mined as spectra and form Galois connections. Descent is a general method for spectral decomposition in algebraic geometry and topology which also leads to generalized Galois connections. Both recommender systems and descent theory are vast research areas, separated by a technical gap so large that trying to establish a link would seem foolish. Yet a formal link emerged, all on its own, bottom-up, against authors' intentions and better judgment. Familiar problems of data analysis led to a novel solution in category theory. The present paper arose from a series of earlier efforts to provide a top-down account of these developments.
\end{abstract}

\maketitle

\tableofcontents

\section{Introduction}\label{Sec:intro}

\subsection{Idea and result}

\begin{definition} We say that an adjunction $F = (\adj F:\BBb\to \AAa)$ is \emph{nuclear}\/ when the right adjoint $\radj F$ is monadic and the left adjoint $\ladj F$ is comonadic. 
\end{definition}

\begin{figure}[!ht]
\begin{center} 
\hspace{-2em}  
\raisebox{2cm}{$\begin{tikzar}[row sep=2.5cm,column sep=3cm]
\AAa \arrow[phantom]{d}[description]{\dashv}  
\arrow[loop, out = 135, in = 45, looseness = 4]{}[swap]{\lft F} 
\arrow[bend right = 13]{d}[swap]{\ladj F} 
\arrow{r}[description]{\mbox{\Huge$\simeq$}} 
\& 
\Emc\BBb F 
\arrow[phantom]{d}[description]{{\dashv}}   \arrow[bend right = 13]{d}[swap]{\lnadj{F}}  
\\
\BBb  \arrow[loop, out = -45, in=-135, looseness = 6]{}[swap]{\rgt F} 
\arrow{r}[description]{\mbox{\Huge$\simeq$}}  
\arrow[bend right = 13]{u}[swap]{\radj F} 
\& 
\Emm \AAa F \arrow[bend right = 13]{u}[swap]{\rnadj{F}}
\end{tikzar}$}
\hspace{4em} \includegraphics[height=4cm,angle=90,origin=c]{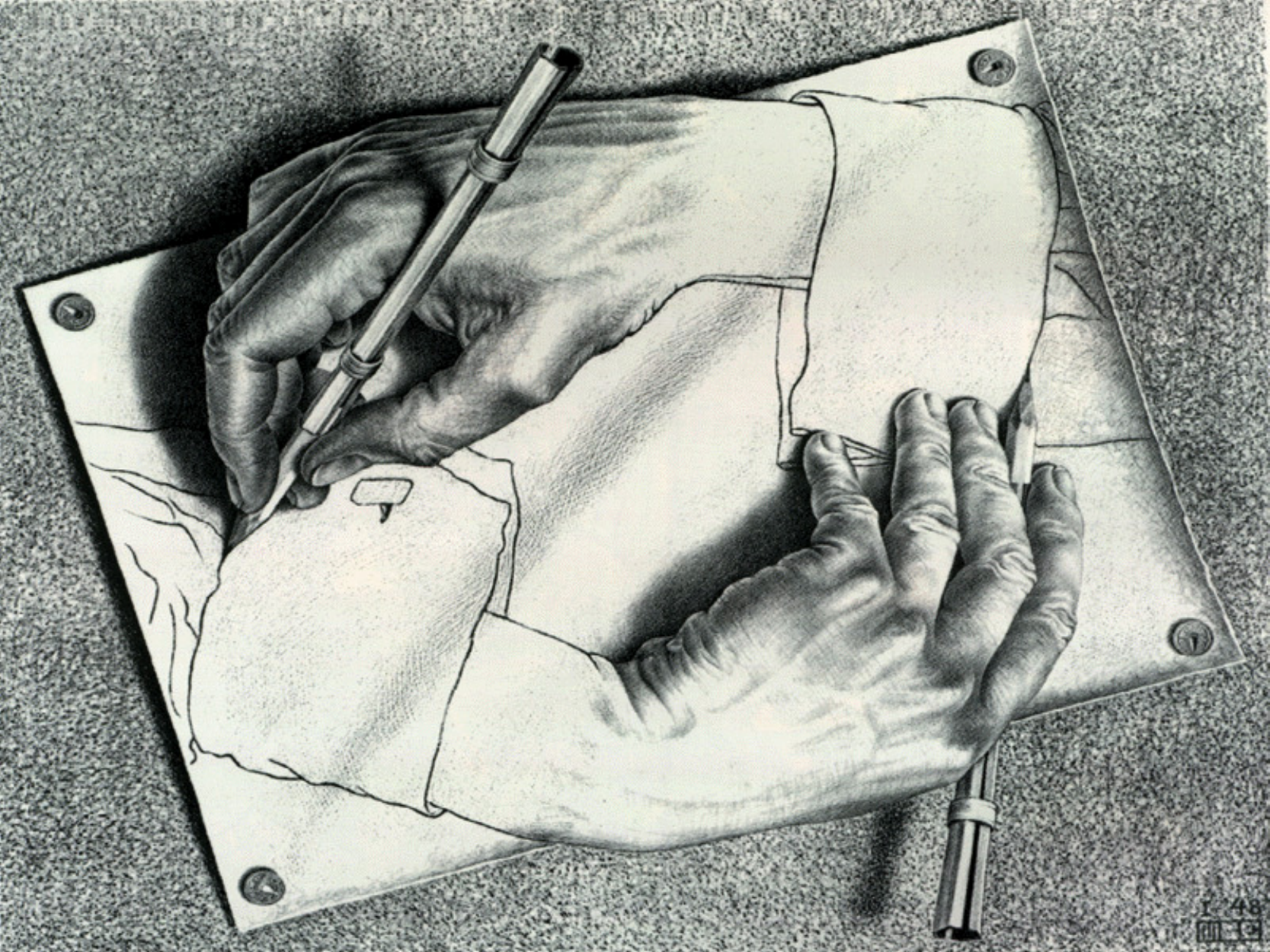}
\caption{An adjunction $\left(\adj F\right)$ is nuclear when $\AAa \simeq \Emc \BBb F$ and $\BBb \simeq \Emm \AAa F$.
}
\label{Fig:escher}
\end{center}
\end{figure}

\para{Explanation.} For an adjunction $F = (\adj F:\BBb\to \AAa)$, we say that 
\begin{itemize}
\item $\radj F$ is monadic when $\BBb$ is equivalent to the category $\Emm \AAa F$ of algebras for the monad $\lft F = \radj F \ladj F:\AAa\to \AAa$, and that
\item $\ladj F$ is comonadic when $\AAa$ is equivalent to the category $\Emc \BBb F$ of coalgebras for the comonad $\rgt F = \ladj F \radj F:\BBb\to \BBb$.
\end{itemize}
The adjunction $F$ is thus nuclear when it displays how the categories $\AAa$ and $\BBb$ determine each other along it: $\AAa$ as $\Emc\BBb F$ and $\BBb$ as $\Emm \AAa F$.  The situation is reminiscent of Maurits Escher's \emph{``Drawing hands''} in Fig.\ref{Fig:escher}. In a formal sense, nuclear adjunctions are the lax version of equivalences.

\para{Claim.} Any adjunction $F$ induces a nuclear adjunction $\NucL F$
\beq\label{eq:NucL}
\prooftree 
F = (\adj F:\BBb\to \AAa)
\justifies
\NucL F = \left(\nadj F\colon \Emm \AAa F \to \Emc \BBb F \right)
\endprooftree
\eeq
where $\rnadj F$ is formed by composing the forgetful functor $\Emm \AAa F \to \AAa$ with the comparison functor $\AAa \to \Emc \BBb F$, and $\lnadj F$ by composing $\Emc \BBb F \to \BBb$ with $\BBb\to \Emm \AAa F$. Hence the left-hand square in Fig.~\ref{Fig:nucadj}. 
\begin{figure}[!ht]
\[\begin{tikzar}[row sep=2.5cm,column sep=3cm]
\AAa \arrow[phantom]{d}[description]{\dashv}  
\arrow[loop, out = 135, in = 45, looseness = 4]{}[swap]{\lft F} 
\arrow[bend right = 13]{d}[swap]{\ladj F} 
\arrow{r}
\& 
\Emc\BBb F \arrow[loop, out = 135, in = 45, looseness = 2.5]{}[swap]{\Lft F}
\arrow[phantom]{d}[description]{{\dashv}}   \arrow[bend right = 13]{d}[swap]{\lnadj{F}}  \arrow{r}[description]{\mbox{\Huge$\simeq$}} \& \left(\Emm\AAa F\right)^{\Rgt F} \arrow[bend right = 13]{d}[swap]{\lnnadj F} \arrow[phantom]{d}[description]{\dashv}
\\
\BBb  \arrow[loop, out = -45, in=-135, looseness = 6]{}[swap]{\rgt F} 
\arrow{r}
\arrow[bend right = 13]{u}[swap]{\radj F} 
\& 
\Emm \AAa F \arrow[loop, out = -45, in=-135, looseness = 6]{}[swap]{\Rgt F} \arrow[bend right = 13]{u}[swap]{\rnadj{F}}  \arrow{r}[description]{\mbox{\Huge$\simeq$}} \& \left(\Emc\BBb F\right)^{\Lft F}\arrow[bend right = 13]{u}[swap]{\rnnadj F} 
\end{tikzar}\]
\caption{The nucleus construction induces an idempotent monad on adjunctions.}
\label{Fig:nucadj}
\end{figure}
We prove that the right-hand square is an equivalence of adjunctions.

\para{Upshot.}  The equivalences in Fig.~\ref{Fig:nucadj} present $\lft F$-algebras as algebras over $\rgt F$-coalgebras and $\rgt F$-coalgebras as coalgebras over $\lft F$-algebras. Simplifying these presentations furnishes the \emph{simple}\/ nucleus construction in Sec.~\ref{Sec:simple}, which opens up a new view of $\lft F$-algebras and $\rgt F$-coalgebras, complementing the familiar Eilenberg-Moore construction \cite{Eilenberg-Moore}. This new view was used as a  programming tool in \cite{PavlovicD:LICS17} and  as a mathematical tool in \cite{PavlovicD:dede}. The resulting reconstruction of monadicity and comonadicity in terms of idempotent splittings echoes Par\'e's explanations in terms of absolute colimits \cite{PareR:absolute-coeq,PareR:absolute}, and contrasts with Beck's fascinating but less transparent treatment in terms of split coequalizers \cite{BarrM:ttt,BeckJ:thesis}. Applications branch in many directions, spanning a wide gamut from descent theory in algebraic geometry to concept analysis on the web. Curiously, the presented results did not evolve from pure mathematics to its applications, but the other way around. Conceptual advances were driven by the ideas evolving in response to practical problems.

\subsection{Background} Nuclear adjunctions have been studied since the early days of category theory, albeit without a name. The problem of characterizing situations when the left adjoint of a monadic functor is comonadic is the topic of Michael Barr's paper in the proceedings of the legendary Battelle conference \cite{BarrM:algCoalg}. 
Ever since, the coalgebras over algebras and the algebras over coalgebras have been emerging in the wide range of applications, from semantics of computation \cite[and the references therein]{KurzA:algcoalg,jacobs1994coalgebras,JacobsB:bases} on one end, to algebraic geometry \cite{Caenepeel:galois,mesablishvili2006monads} and homotopy theory on the other \cite{BalmerP:annalen,HessK:general}. The monad-comonad couplings on two sides of an adjunction could be gleaned already on the primordial examples of adjunctions \cite{KanD:adj} and the efforts to grasp them propelled the early research relating them with monads and comonads \cite{Applegate-Tierney:iterated,BeckJ:thesis,HuberP:homotopy}. There are, however, different ways in which monad-comonad couplings arise from adjunctions, and from other monads and comonads. In \cite{Applegate-Tierney:models}, Applegate and Tierney  considered monads and comonads formed on the two sides of the comparison functors from adjunctions to their final resolutions in the induced categories of coalgebras. They found that such couplings in general induce further such couplings and form transfinite towers. We describe this process in Sec.~\ref{Sec:HT}. Confusingly, Applegate-Tierney's towers of adjunctions of comparison functors left a false impression that the adjunctions between algebras over coalgebras and coalgebras over algebras also lead to transfinite towers. This impression blended into folklore and spread through literature.\footnote{There is an interesting exception outside the literature. In a fax message sent to Paul Taylor on 9/9/99  \cite{Lack-Taylor}, a copy of which was kindly provided after the present paper appeared on arxiv, Steve Lack set out to determine the conditions under which the tower of coalgebras over algebras, which "a priori continues indefinitely", settles to equivalence at a finite stage. Within 7 pages of diagrams, the question was reduced to splitting a certain idempotent. While the argument is succinct, it does seem to prove a claim which, together with its dual, seems equivalent to our Prop.~\ref{prop:nuc}. The claim was not pursued any further. This episode from the early life of the nucleus underscores its message: that a concept is within reach whenever there is an adjunction  but it does need to be spelled out to be recognized.} We show that the latter towers settle in one step, displayed in Fig.~\ref{Fig:nucadj}.

\subsubsection*{Name} Although they appeared in many avatars, nuclear adjunctions don't seem to have been named. We call them nuclear in reference to the nuclear operators on Banach spaces, which generalize spectral decompositions of hermitians and the singular value decompositions of linear operators. The terminology was introduced in Grothendieck's thesis \cite{GrothendieckA:memAMS}. 

\subsection{Overview} 
\begin{figure}[!ht]
\begin{center}
\begin{tikzar}[row sep = 4em,column sep = 1em]
\mbox{\textit{\textbf{matrices}}} \&\&  \mbox{\textit{\textbf{extensions}}} \& \mbox{\textit{\textbf{localizations}}} \& \mbox{\textit{\textbf{nuclei}}}
\\[-7.5ex] 
\&\&\& \Mnd \ar[bend right=15,tail]{dl}[swap]{\EM} 
\ar[bend left=15,two heads]{dr}{\MN} 
\ar[phantom]{dr}[rotate = -45]{\top} 
\\
\Mat \ar
{rr}{\MA} 
\&\& 
\Adj  
\ar[bend right=15,two heads]{ur}[swap]{\AM} 
\ar[phantom]{ur}[rotate = 45]{\top}
\ar[bend left=15,two heads]{dr}{\AC} 
\ar[phantom]{dr}[rotate = -45]{\bot} 
\&\&
\Nucl
\ar[bend left=15,tail]{ul}{\NM} 
\ar[bend right=15,tail]{dl}[swap]{\NC} 
\\
\&\&\& \Cmn 
\ar[bend left=15,tail]{ul}{\EC} 
\ar[bend right=15,two heads]{ur}[swap]{\CN} 
\ar[phantom]{ur}[rotate = 45]{\bot} 
\end{tikzar}
\caption{The paths to the nucleus}
\label{Fig:Nuc}
\end{center}
\end{figure}
The paths to the nucleus are mapped in Fig.~\ref{Fig:Nuc}. They lead from observations tabulated in data {\sf Mat}{\small rices} to concepts collected in {\sf Nuc}{\small lei}. They lead through {\sf Adj}{\small unctions} and then branch through {\sf M}{\small o}{\sf n}{\small a}{\sf d}{\small s} and {\sf C}{\small o}{\sf m}{\small o}{\sf n}{\small ads}. The {\sf Nuc}{\small lei} arise at the intersection of the two branches. Sec.~\ref{Sec:FCA} presents the posetal case, which is intuitive and well-known\footnote{In the theory of locales, the word ``nucleus'' is sometimes used to refer to the meet-preserving closure operators and to the sublocales that they induce \cite[\S II.2.2--4]{SimmonsH:framework,JohnstoneP:stone}. This usage is unrelated to the present concept or its background.}. The categorical setting is in Sec.~\ref{Sec:cat}. The nucleus functor is defined in Sec.~\ref{Sec:prop}.  
The fact that it is an idempotent monad is proved in Sec.~\ref{Sec:Theorem}. The \emph{simple nucleus}\/ construction, offering an alternative to the familiar Eilenberg-Moore presentation of algebras and coalgebras, is in Sec.~\ref{Sec:simple}. Sec.~\ref{Sec:HT} touches upon the ways in which the nucleus concept feeds back into the ideas of descent and of homotopy-invariant space decompositions, from which it indirectly originates. The closing section returns to the big picture and comments about the remaining work.

\section{Posetal nuclei: from contexts to concepts}
\label{Sec:FCA}

\subsection{Idea of Formal Concept Analysis} 
\begin{figure}[!ht]
\begin{center}
\includegraphics[height=12cm]{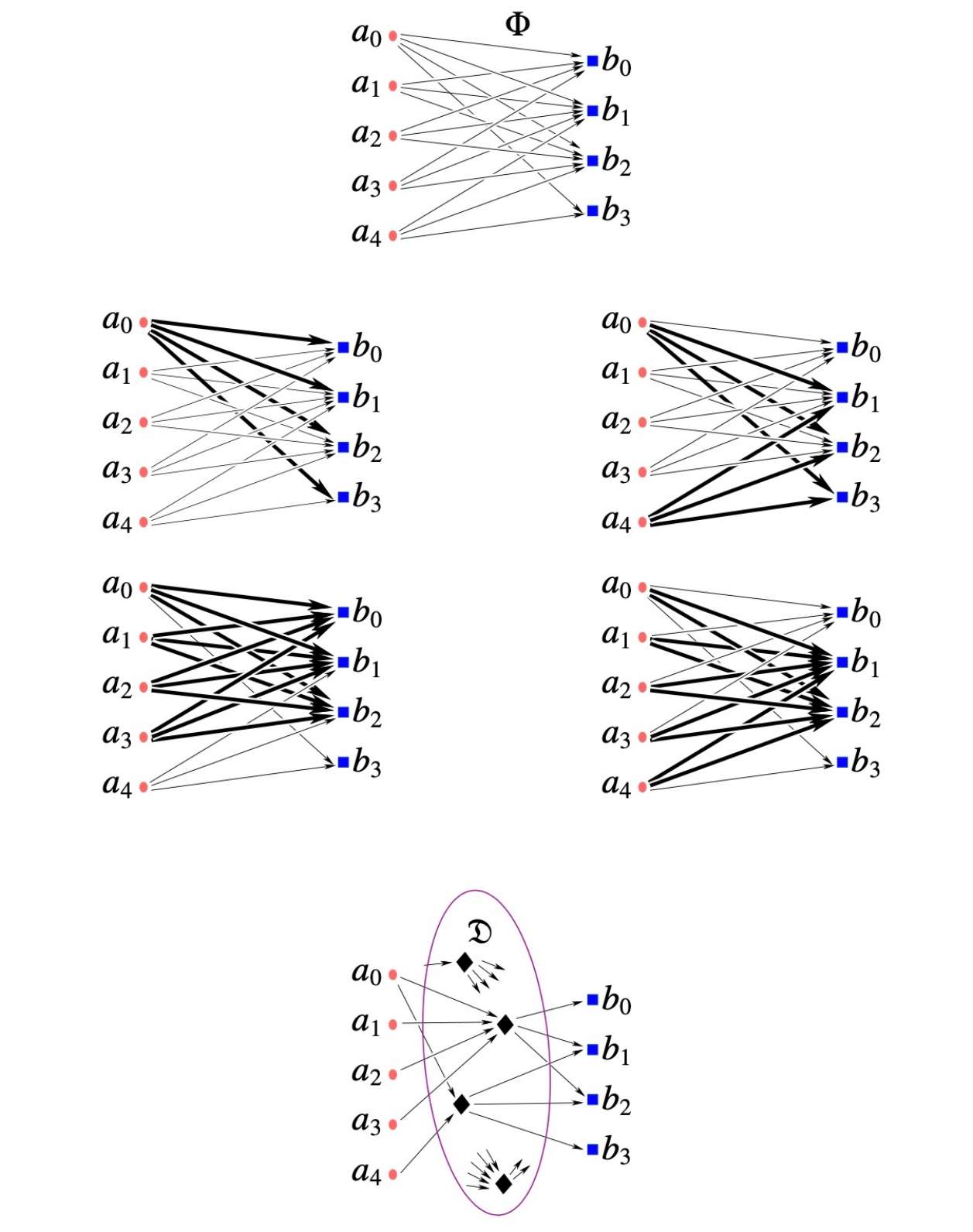}
\caption{A context $\Phi$, its four concepts, and their concept lattice}
\label{Fig:FCA}
\end{center}
\end{figure}

Consider a market with $A$ sellers and $B$  buyers. Their interactions are recorded in an adjacency matrix $A\times B \tto \Phi 2$, where $2$ is the set $\{0,1\}$, and the entry $\Phi_{ab}$ is 1 if the seller $a\in A$ at some point sold goods to the buyer $b\in B$; otherwise it is 0. It is often convenient to move between the adjacency matrices and the induced binary relations along the correspondence
 \bea\label{eq:comprehrel}
 \left(A\times B \tto \Phi 2\right) & \leftrightsquigarrow & \eh \Phi = \{<a,b>\in A\times B\ |\ \Phi_{ab} = 1\}
 \eea
which allows writing $a\eh \Phi b$ instead of $\Phi_{ab} = 1$. In the literature on  Formal Concept Analysis (FCA) \cite{Carpineto-Romano:book,FCA-book,FCA-foundations}, the adjacency matrices $\Phi$ and the corresponding relations $\eh \Phi$ are called \emph{contexts}. The latent \emph{concepts}\/ that they carry are extracted as pairs $<U,V>\in \WP A\times \WP B$ where every seller in $U$ is $\eh \Phi$-related to every buyer in $V$, and vice versa. Fig.~\ref{Fig:FCA} displays an example.  
The binary relation $\eh \Phi \subseteq A\times B$ is presented as a bipartite graph, and the pairs $<U,V>$ are complete subgraphs. E.g., if buyers $a_0$ and $a_4$ have farms, and sellers $b_1$, $b_2$ and $b_3$ sell farming equipment, but seller $b_0$ does not, then the sets $U = \{a_0, a_4\}$ and $V = \{b_1, b_2, b_3\}$ span a complete subgraph of the bipartite graph $\Phi$, which corresponds to the concept  \emph{"farming"}. If the buyers from the set $U'=\{a_0, a_1,a_2,a_3\}$ own cars but the buyer $a_4$ does not, and the sellers $V' = \{b_0, b_1, b_2\}$ sell car accessories but the seller $b_3$ does not, then the sets $<U',V'>$ span another complete subgraph, this time corresponding to the concept \emph{"car"}. The latent concepts are extracted from a context, viewed as a bipartite graph, as its complete subgraphs.

\subsection{Formalizing concept analysis} The subgraph of a bipartite graph $\eh \Phi \subseteq A\times B$ spanned by a pair of sets of nodes $<U,V>\in \WP A\times \WP B$ is complete if
\[
U \ = \ \bigcap_{v\in V} \{x\in A\ |\ x\eh \Phi v\}  \qquad\qquad \qquad V = \bigcap_{u\in U} \{y\in B\ |\  u\eh \Phi y\}\]
Such pairs are ordered by the relation
\bea\label{eq:cutorder}
<U,V>\leq <U',V'> & \iff & U\subseteq U' \ \wedge \ V\supseteq V'
\eea
They are obviously a subposet of the lattice $\Wp A\times \Wpp B$, where $\Wp A$ is the set of subsets of $A$ ordered by the inclusion $\subseteq$ and $\Wpp B$ is the set of subsets of $B$ ordered by reverse inclusion $\supseteq$. We will see in a moment that they are also a quotient of this lattice, and thus form a lattice $\Cut\Phi$. This is the \emph{concept lattice}\/ induced by the \emph{context matrix} $\Phi$. 

\subsection{Posetal matrices} The posetal case of Fig.~\ref{Fig:Nuc} starts from the category
\bea\label{eq:matrp}
|\Matrp | & = & \coprod_{A,B\in \Pos} \Pos(A^\op \times B, \TTwo)
\\[1.5ex]
\Matrp(\Phi, \Psi) & = & \{<h,k>\in \Pos(A,C)\times \Pos(B,D) \ |\ \Phi(a, b) =  \Psi (ha,kb)\}\notag 
\eea
$A,B,C,D$ are posets and $\Phi \in \Pos(A^\op \times B, \TTwo)$ and $\Psi \in \Pos(C^\op \times D, \TTwo)$ are matrices with entries from the poset $\TTwo= \{0\lt1\}$. Note that the sets $A$ of sellers and $B$ of buyers, discussed above, generally come with partial orders, e.g. from previous concept analyses. 

\para{Comprehensive view.} Generalizing \eqref{eq:comprehrel}, context matrices can also be viewed as subposets of product posets using the isomorphisms
\bea\label{eq:compreh-pos}
\Pos (A^\op \times B , \TTwo) &
\begin{tikzar}[row sep = 4em]\hspace{.1ex} \ar[bend left]{r}{\eh{(-)}} \ar[phantom]{r}[description]{\cong} 
\& \hspace{.1ex}
\ar[bend left]{l}{\chi}
\end{tikzar} & \Sub\diagup A\times B^\op
\\[1ex]
\Phi & \mapsto & \eh{\Phi} = \{<x,y> \in A\times B^\op \ |\ \Phi(x,y) = 1\}
\notag
\\[2ex]
 \notag
\chi_S(x,y) = 
{\scriptstyle \left.\begin{cases} 1 & \mbox{ if  } <x,y>\in S\\
0  & \mbox{ otherwise}\end{cases} \right\}}
&\mapsfrom & \Big(S\subseteq A\times B^\op\Big)  
\eea
The poset $\eh \Phi$ can be thought of as the \emph{``comprehension''} of $\Phi$ viewed as a \emph{``predicate''}  \cite[Sec.~III.3]{LawvereFW:compreh,PavlovicD:thesis}. The poset $\eh \Phi$ refines the set-theoretic comprehension schema by imposing the order-closure requirement:
\bea\label{eq:monotone}
a \leq a' \ \wedge\ a' \eh \Phi b'\ \wedge\ b' \leq b & \implies & a\eh \Phi b 
\eea

\subsection{Completions make matrices representable as adjunctions} The concepts are extracted from a context $\eh \Phi \subseteq A\times B$ by finding the pairs of a lower-closed subset of $A$ and an upper-closed subset of $B$ coupled by the context.  The lower-closed and the upper-closed subsets of a poset form the complete semilattices $(\Do A, \bigvee )$ and $(\Up B, \bigwedge)$, where
\beq
\Do A \ = \ \{X\subseteq A\ |\ a \leq a' \in X \implies a \in X\} \qquad\qquad 
\Up B  \ = \  \{Y\subseteq B\ |\ Y\ni b' \leq b  \implies Y \ni b\}  
\label{eq:DoA-Upb-pos}
\eeq
and $\bigvee$ in $\Do A$ and $\bigwedge$ in $\Up B$ are the unions. The embeddings 
\begin{align*}
A &\tto \mnd \Do A & B&\tto\cmn \Up B\\
a &\mapsto \{x\leq a\} & b &\mapsto \{y\geq b\}
\end{align*}
are respectively the $\bigvee$-completion of $A$ and the $\bigwedge$-completion of $B$. These semilattice completions allow extending the context matrix $\eh \Phi\subseteq  A \times B$ from \eqref{eq:comprehrel} along \eqref{eq:compreh-pos} to $\overline \Phi\subseteq  \Do A \times \Up B$, setting
\bea
X\overline \Phi Y & \iff & 
\forall a\in X\ \forall b \in Y.\ a\eh\Phi  b
\eea
Since its domain and codomain are complete semilattices $\Do A$ and $\Up B$, the matrix $\overline \Phi$ is representable in the form 
\beq\label{eq:adjunctions-pos} 
\ladj \Phi X \subseteq Y\  \ \iff X\overline \Phi Y \ \ \iff \ \  X \supseteq \radj \Phi Y\eeq
The adjoints $\adj \Phi$ extract the complete bipartite subgraphs, illustrated in Fig.~\ref{Fig:FCA}, as intersections:
\beq\label{eq:galois-pos}
\begin{tikzar}{}
X \ar[mapsto]{dd} \& \Do A \ar[bend right=15]{dd}[swap]{\ladj \Phi}\ar[phantom]{dd}[description]\dashv \& 
\displaystyle \bigcap_{y\in Y}\  _\bullet  \Phi y
\\
\\
\displaystyle \bigcap_{x\in X} \ x  \Phi_\bullet
\& \Up B \ar[bend right=15]{uu}[swap]{\radj \Phi}
\&   Y \ar[mapsto]{uu}
\end{tikzar}
\eeq
where $_\bullet   \Phi y =\{x\in A\ |\ x   \Phi y\}$ and $x  \Phi_\bullet  = \{y\in B\ |\  x\eh \Phi y\}$ are the images along the transpositions
\beq
\prooftree
\Phi:A^\op \times B \to \TTwo
\justifies 
_\bullet   \Phi :B\to \TTwo^{A^{\op}}\cong\Do A \qquad \qquad 
 \Phi_\bullet :  A\to\left(\TTwo^{B}\right)^{\op}\cong  \Up B
\endprooftree
\eeq
The poset adjunctions like \eqref{eq:galois-pos} are often called \emph{Galois connections}. They form the category 
\bea\label{eq:adjp}
|\Adjp | & = & \coprod_{A,B\in \Pos} \left\{<\ladj \Phi, \radj \Phi> \in \Pos(A, B) \times \Pos(B,A)\ |\ \ladj \Phi x \leq y \iff x\leq \radj \Phi y
\right\}
\\[1ex]
\Adjp\big(\Phi, \Psi\big) & = & \big\{<H,K>\in \Pos(A,C)\times \Pos(B,D) \ |\ K\ladj \Phi  = \ladj \Psi H\  \wedge\  H \radj \Phi =  \radj \Psi H\big\}\notag 
\eea
The matrix extension at the first step of the posetal case of Fig.~\ref{Fig:Nuc} is thus  
\bea
\MAp : \Matrp & \longrightarrow &  \Adjp\\
\Phi & \longmapsto & \left(\adj \Phi :\Up B\to \Do A\right) \mbox{ defined in  \eqref{eq:galois-pos}}\notag
\eea 

\subsection{Posetal adjunctions induce closure and interior operators} The further steps through the posetal Fig.~\ref{Fig:Nuc} go to the categories of closure and interior operators, respectively:
\bea
|\Mndp | & = & \coprod_{A\in \Pos} \left\{\lft \Phi \in \Pos(A, A)\ |\  x \leq \lft \Phi x \geq \lft\Phi\lft \Phi(x) \right\}\label{eq:mndp}
\\
\Mndp\big(\lft\Phi_{A}, \lft\Psi_{C}\big) & = & \left\{H\in \Pos(A,C) \ |\ H \lft \Phi = \lft \Psi H\right\}\notag
\\[2ex]
 |\Cmnp | & = & \coprod_{B\in \Pos} \left\{\rgt\Phi \in \Pos(B, B)\ |\ x \geq \rgt\Phi x \leq 
\rgt\Phi\rgt\Phi x
\right\}\label{eq:cmnp}
\\
\Cmnp\big(\rgt\Phi_{B}, \rgt\Psi_{D}\big) & = & \left\{ K\in \Pos(B,D) \ |\ K\rgt \Phi  = \rgt \Psi K\right\}\notag
\eea
along the functors
\beq\label{eq:huber-pos}
\begin{tikzar}[row sep=1ex,column sep=4em]
\Cmnp 
\& \Adjp \ar[thick]{l}[swap]{\ACp} \ar[thick]{r}{\AMp} 
\&
\Mndp
\\
\ladj \Phi\radj\Phi 
\& \left( \adj \Phi\right) \ar[mapsto,shorten=3ex]{l} \ar[mapsto,shorten=3ex]{r} 
\& \radj\Phi \ladj \Phi
\end{tikzar}
\eeq
%
\subsection{Closure and interior operators have adjunction resolutions} The functors from adjunctions to closure and interior operators have full and faithful left adjoints:  
\bea\label{eq:resolutions-pos}
\Cmnp & \tto{\hspace{2em}\ \ECp\hspace{2em}}  \Adjp\  \oot{\hspace{2em} \EMp\hspace{2em}} & \Mndp\\[2ex]
&\hspace{2.9em}\Big( \begin{tikzar}{}
A \ar[bend right=20,two heads]{r}\ar[phantom]{r}[description]{\top} \& \Emm{A}\Phi \ar[bend right=20,hook]{l}
\end{tikzar}\Big)\hspace{1.2em}  \longmapsfrom & \hspace{.8em} \lft \Phi_{A} \notag\\
\rgt\Phi_{B}\hspace{.8em} & \longmapsto \hspace{1.1em}  \Big(\begin{tikzar}{}
\Emc{B}\Phi \ar[bend right=20,hook]{r} \ar[phantom]{r}[description]{\top} \& B \ar[bend right=20,two heads]{l}
\end{tikzar}\Big)\hspace{2.9em} &\notag 
\eea
where the posets
\beq\label{eq:EM-pos}
\Emm{A}\Phi = \{x\in A\  |\ \lft  \Phi x = x\} \qquad \qquad \qquad 
\Emc{B}\Phi = \{y\in B\  |\ \rgt  \Phi y = y\} 
\eeq
inherit the order from $A$ and $B$. The claimed adjunctions $\AMp\dashv\EMp$ and $\ACp\dashv\ECp$ are easily checked, as are the natural isomorphisms $\AMp\EMp \cong \id$ and $\ACp\ECp \cong \id$, making $\Mndp$ and $\Cmnp$ into reflective subcategories of $\Adjp$. 

\subsection{Concept lattices as nuclei} The posetal version of Fig.~\ref{Fig:Nuc}
is comprised of the following steps
\beq\label{eq:steps}
\prooftree
\prooftree
\prooftree
\Phi : A^\op \times B \to \TTwo
\justifies
\MAp\Phi = \AdjL\Phi = \Bigg( \begin{tikzar}{}
\Do A \ar[bend right=20]{r}[swap]{\ladj \Phi} \ar[phantom]{r}[description]{\top} \& \Up B \ar[bend right=20]{l}[swap]{\radj \Phi}
\end{tikzar}\Bigg) 
\endprooftree
\justifies
\EMnd_0\AdjL\Phi = \Bigg( \begin{tikzar}{}
\Do A \ar[bend right=20,two heads]{r}
\ar[phantom]{r}[description]{\top} \& \Emm{\Do A}\Phi \ar[bend right=20,hook]{l}
\end{tikzar}\Bigg)
\qquad \qquad \qquad \qquad \ECmn_0\AdjL\Phi = \Bigg( \begin{tikzar}{}
\Emc{\Up B}\Phi \ar[bend right=20,hook]{r} \ar[phantom]{r}[description]{\top} \& \Up B \ar[bend right=20,two heads]{l}
\end{tikzar}\Bigg)
\endprooftree
\justifies
\NucL_0 \Phi= \Bigg(  \begin{tikzar}{}
\Up B^{\rgt \Phi} \ar[bend right=20,tail,two heads]{r}[swap]{\lnadj \Phi} \ar[phantom]{r}[description]{\cong} \& \Do A^{\lft \Phi} \ar[bend right=20,tail,two heads]{l}[swap]{\rnadj \Phi}
\end{tikzar}\Bigg)
\endprooftree
\eeq
where $\EMnd_0 = \EMp\circ \AMp$, and $\ECmn_0 = \ECp \circ \ACp$ are the idempotent monads displaying $\Mndp$ and $\Cmnp$ as reflective subcategories of $\Adjp$. A peculiarity of the posetal case is that they are also idempotent comonads, arising from the adjunctions $\EMp\dashv\AMp$ and $\ECp\dashv\ACp$, from which it follows not only that $\Mndp$ and $\Cmnp$ are also coreflective in $\Adjp$, but also that they are equivalent. A concept lattice can thus be obtained using $\EMnd_0$ and presented in terms of a closure operator, or using $\ECmn_0$ and presented in terms of an interior operator, or using $\NucL_0 \cong \EMnd_0\ECmn_0\cong\ECmn_0\EMnd_0$ and presented in terms of an adjunction, reducing the context $\Phi\colon A^{\op}\times B \tto \Phi \TTwo$ to the concept lattice
\bea\label{eq:nucone}
\Cut \Phi & = &  \big\{<X,Y> \in \Do A\times \Up B\ |\ X = \radj \Phi Y \ \wedge\ \ladj  \Phi X = Y\big\}
\eea
The elements of this lattice generalize Dedekind cuts of rational numbers, and the nucleus construction generalizes the Dedekind-MacNeille completion to matrices \cite{Banaschewski-Bruns,DedekindR:zahlen,MacNeille}. All of this is well-known and valid not only for categories of posets but also for posetally enriched categories. Our present concern is the general case, where the categories $\Mnd$ and $\Cmn$ are not equivalent.

\section{From categorical contexts to adjunctions, monads, and comonads}
\label{Sec:cat}

\para{The need for categorical concept analysis.} The posetal nuclei from the preceding section lift to poset-enriched categories with no major surprises, though with significant new applications \cite{BradleyT:enriched,PavlovicD:ICFCA12,PavlovicD:Samson13,WillertonS:nucleus}. Large-scale concept analyses in recommender systems \cite{Recommender-handbook}, however,  require going beyond the posetal and numeric contexts, to proper categories. The reason is that the posetal contexts do not capture out-of-band dependencies or hidden variables. The vector space models preclude hidden variables by imposing a tacit but very consequential assumption that all analyzed sources are statistically independent. This assumption is implied by the linearity of the operations on context matrices which aggregate correlations as sums of \emph{products}\/ of the observed frequencies. While the source interference is normally measured as the deviation from the product distribution,  in the vector space model, the product distribution is built-in. This conceptual shortcoming, briefly described at the outset of our work \cite{PavlovicD:CALCO15}, in the meantime grew into a practical problem of tremendous impact, leading to wholesale correlation amplifications, information cascades, and echo chambers. On the theoretical side, the task could be stated in terms of Grothendieck's Galois descent \cite{GrothendieckA:fibrations59,GrothendieckA:SGA1,Janelidze-Tholen:facets,VistoliA:descent} but required substantial reinterpretation and narrowing. The basic ideas, sketched in  \cite{PavlovicD:CALCO15,PavlovicD:LICS17}, evolved bottom-up, from implementations to constructions. Here we attempt to provide a top-down account.

\subsection{Lifting concept analysis to categories --- and not further}
Spanned over categories $\AAa$ and $\BBb$, data matrices become functors $\AAa^{\op}\times \BBb\tto\Phi \Set$, assigning to each pair of objects $<a,b>$ the set $\Phi_{ab}$ of observations linking them. Whether the categories $\AAa$ and $\BBb$ represent sellers and buyers like in Sec.~\ref{Sec:FCA}, or structural components and functional modules of a machine like in \cite{PavlovicD:CALCO15}, or some other pair of types, the matrix entries record the interactions between the entities of type $\AAa$ and the entities of type $\BBb$ and the arrows in each of the categories capture the dependencies within each of the types. The arrow part of the functor $\Phi$ thus tracks how the dependencies propagate through the interactions. In category theory, the functors generalizing matrices in this way are variably called distributors, profunctors, or bimodules \cite[Vol.~1, Sec.~7.8]{BorceuxF:handbook}. We continue to call them matrices not so much to expand the terminological choices but to remember the application at hand. 

\para{The 2-dimensional structure.} Given the matrices $\AAa^{\op}\times \BBb\tto\Phi \Set$ and $\BBb^{\op}\times \CCc\tto\Psi \Set$, the usual matrix composition is generalized using the \emph{coend}\/ operation
\bea
\left(\Phi;\Psi\right)_{ac} & = & \supp_{y\in \BBb} \left(\Phi_{ay}\times \Psi_{yc}\right)
\eea
See \cite{LoregianF:coend} for details. The compositions accumulate and propagate data dependencies between $\AAa$ and $\CCc$ across the intermediary $\BBb$. The main point is that, for each $\AAa$ and $\BBb$, the $\AAa^{\op}\times \BBb$-matrices form a category recording their correlations, and that the correlations can be composed, functorially up to natural isomorphism. The categorical matrices (i.e. distributors, profunctors, bimodules) thus form a \emph{bicategory}\/ \cite[Vol.~1,Sec.~7.7]{BenabouJ:bicats,BorceuxF:handbook}. The categories of adjunctions, monads, and comonads form 2-categories, meaning that the morphisms between any pair of adjunctions (resp. monads, comonads) also form a category, and the representable data correlations captured there also compose functorially, this time not up to isomorphism but strictly  \cite[Vol.~1,Sec.~7.1]{BorceuxF:handbook}. The 2-categorical structure of adjunctions, monads, and comonads goes back to \cite{Auderset,Street-Lack:monads,Street-Shanuel,StreetR:monads}, and could be gleaned already in  \cite{MarandaJM:adj}. 

\para{The nucleus construction is 1-dimensional.} The nucleus construction combines monadic and comonadic adjunctions. To lift a monad morphism to a morphism  between the induced monadic adjunctions, we use the 2-cell component of the monad morphism, which stipulates how the 1-cell distributes over the two monads. Ditto for lifting a comonad morphism to a morphism between the induced comonadic adjunctions. The problem that arises for the nucleus construction is that the two distributivities, one needed to map monad morphisms to adjunction morphisms between  final resolutions, the other one to map the comonad morphisms similarly, are required to run in opposite directions. To compose the two 2-monads on the 2-category of  adjunctions, one making them monadic, the other making them comonadic, the 2-dimensional structure must be restricted to invertible distributivity 2-cells. But then the invertible distributivities can be absorbed into 1-cells and the 2-cells can be left out. This will turn out to be the feature making the nucleus monad strongly idempotent. In 2-dimensional category theory, a monad is usually weakly idempotent, in the lax sense that its algebras are adjoint to the unit \cite{KockA:zoeberlein,StreetR:fib-bicat,Zoeberlein}. The nucleus monad is strongly idempotent, in the sense that its algebras are equivalences. It is a categorical construction, not 2-categorical. The following definitions and the results that they enable provide the evidence.

\subsection{Matrices (distributors, profunctors, bimodules)} The category of matrices over categories with the morphisms up to invertible 2-cells is defined in Fig.~\ref{Fig:Mat}.
\begin{figure}[!ht]
\begin{center}
\bear
|\Mat| & = & \coprod_{\AAa, \BBb\in \CAT} \CAT(\AAa^{o}\times \BBb, \Set)\\ 
\Mat (\Phi, \Psi) & = & \coprod_{\substack{H\in \CAT(\AAa,\CCc)\\
K\in \CAT(\BBb,\DDd)}}\left\{\begin{tikzar}[row sep = .25ex,column sep = 1ex]\AAa^{\op}\times \BBb \ar{ddrr}[description]{\Phi} \ar{rrrr}{H^{\op}\times K} \&\&\&\& \CCc^{\op}\times \DDd\ar{ddll}[description]{\Psi}\\
\&\&\Leftarrow\gamma\Rightarrow\\
\&\&\Set
\end{tikzar} \right\}\notag
\eear
%
\caption{$\gamma\in \Niso\Big(\Phi,\Psi(H^{\op}\times K)\Big)$ are natural isomorphisms}
\label{Fig:Mat}
\end{center}
\end{figure}

\para{Comprehensive view.} The posetal comprehension \eqref{eq:compreh-pos} lifts in categories to the equivalence of presheaves and discrete fibrations  
\bea\label{eq:compreh-cat}
\Cat (\AAa^\op \times \BBb , \Set) &
\begin{tikzar}[row sep = 4em]\hspace{.1ex} \ar[bend left]{r}{\eh{(-)}} \ar[phantom]{r}[description]{\simeq} 
\& \hspace{.1ex}
\ar[bend left]{l}{\Xi}
\end{tikzar} & \Dfib\diagup \AAa\times \BBb^\op 
\eea
The details are similar to the posetal case  and also widely known, as the categorical fibrations, generalizing the posetal lower sets, go back to Grothendieck \cite{GrothendieckA:fibrations59,GrothendieckA:SGA1} and have been reviewed many times \cite{JacobsB:book,PavlovicD:thesis,Streicher:fibered}. For completeness, we reproduce the details of \eqref{eq:compreh-cat} in Appendix \ref{Appendix:compreh}. Since switching between the $\Set$-functors and their comprehensions is as routinely in category theory as the switching between the predicates as $\TTwo$-functions and their comprehensions is in logic, we often use $\Do \AAa$ and $\Up \BBb$ to refer to both sides of the equivalences
\beq
\Set^{\AAa^{\op}} \simeq \Do \AAa = \Dfib\diagup \AAa \qquad \qquad \qquad \qquad \left(\Set^{\BBb}\right)^{\op}\simeq \Up\BBb = \left(\Dfib\diagup \BBb^{\op}\right)^{\op}
\eeq
where $\Set^{\XXx}$ abbreviates the category  $\Cat (\XXx, \Set)$. The embeddings 
\begin{align*}
\AAa &\tto \mnd \Do \AAa & \BBb&\tto\cmn \Up \BBb\\
a &\mapsto \left(\AAa/a\tto{\Dom} \AAa\right) & b &\mapsto \left(b/\BBb\tto{\Cod} \BBb\right)
\end{align*}
are this time, respectively, the $\supp$-completion of $\AAa$ and the $\inff$-completion of $\BBb$. As usually, $\AAa/a$ denotes the category of $\AAa$-arrows into $a$, whereas $b/\BBb$ denotes the category of $\BBb$ arrows from $b$. The fact that $\Do \AAa$ is a $\supp$-completion and that $\Up \BBb$ is a $\inff$-completion means that every $\left(\XXx\tto X\AAa\right)\in \Do\AAa$ and $\left(\YYy\tto Y\BBb\right)\in\Up\BBb$ satisfy
\beq
X = \supp\left(\XXx\tto X \AAa\tto \mnd \Do \AAa\right) \qquad \qquad\qquad\qquad Y = \inff\left(\YYy\tto Y\BBb\tto \cmn \Up\BBb\right)
\eeq

\subsection{Completions make matrices representable as adjunctions again} \label{Sec:repres}
A matrix $\Phi\colon\AAa^\op \times \BBb \to \Set$ can be represented as the following adjunction between the completions
\beq\label{eq:galois-cat}
\begin{tikzar}{}
\XXx\tto X\AAa \ar[mapsto]{dd} \& \Do \AAa \ar[bend right=15]{dd}[swap]{\ladj \Phi}\ar[phantom]{dd}[description]\dashv \& 
\displaystyle \inff\ \Big(\YYy\tto Y\BBb \tto{_\bullet\Phi} \Do\AAa\Big) 
\\
\\
\displaystyle \supp\ \Big(\XXx\tto{X}\AAa \tto{\Phi_\bullet} \Up\BBb\Big)
\& \Up \BBb \ar[bend right=15]{uu}[swap]{\radj \Phi}
\& \YYy\tto Y \BBb \ar[mapsto]{uu}
\end{tikzar}
\eeq
where $\Phi_\bullet$ and $_\bullet\Phi$ are the two transposes of $\Phi$
\beq\label{eq:deriv}
 \prooftree
 \prooftree
 \Phi\colon\AAa^\op \times \BBb \to \Set
 \justifies
 \Phi_\bullet \colon \AAa \to \Bigl(\Set^\BBb\Bigr)^\op \simeq \Up\BBb \qquad \qquad
_\bullet \Phi \colon \BBb \to \Set^{\AAa^\op} \simeq \Do\AAa  
 \endprooftree
 \justifies
\ladj \Phi \colon \Do \AAa \tto{\hspace{3.85em}} \Up \BBb
 \qquad \qquad\  \ \radj\Phi \colon \Up\BBb \tto{\hspace{3.2em}} \Do\AAa  \endprooftree
\eeq
and the definitions of $\ladj \Phi$ and $\radj\Phi$ are based on the fact that any left adjoint must preserve $\supp$, and that any right adjoint must preserve $\inff$. In other words, $\ladj \Phi$ is the left Kan extension  of $\Phi_\bullet$ along $\mnd\colon \AAa\to\Do \AAa$, whereas $\radj \Phi$ is the right Kan extension of $_\bullet\Phi$ along $\cmn\colon \BBb\to \Up\BBb$. It is useful to notice that both are computed pointwise, just like both extensions in \eqref{eq:galois-pos} are computed as set intersections. The reason is that dualizing along the equivalence $\Up\BBb\simeq\left(\Set^{\BBb}\right)^{\op}$ gives
\bear
\supp\ \Big(\XXx\tto{X}\AAa \tto{\Phi_\bullet} \Up\BBb\Big) & = & \inff \Big(\XXx^{\op}\tto{X^{\op}}\AAa^{\op} \tto{\Phi^{\op}_\bullet} \Set^{\BBb}\Big)
\eear
But the limits in $\Set^\BBb$ are pointwise, since the Yoneda lemma gives for any diagram $D\colon \DDd\to \Set^{\BBb}$
\[
\left(\inff D\right)b \ \  = \ \  \Set^\BBb\left(\cmn b, \inff D\right)\ \  =\ \ \Con(b, D)
\]
In words, the limit of $D$ is the functor $\BBb\to \Set$ whose value at $b$ is the set of commutative cones from $b$ to $D$. The matrix extension
\bea
\MA : \Mat & \longrightarrow &  \Adj\\
\Phi & \longmapsto & \left(\adj \Phi :\Up \BBb\to \Do \AAa\right) \mbox{ defined in  \eqref{eq:galois-cat}}\notag
\eea 
takes us to the category of adjunctions in Fig.~\ref{Fig:Adj}.
\begin{figure}[!ht]
\begin{center}
\bear
|\Adj| & = & \coprod_{\AAa, \BBb\in \CAT}\  \ \coprod_{\substack{\ladj F\in \CAT(\AAa,\BBb)\\
\radj F\in \CAT(\BBb,\AAa)}}\ \ 
\left\{\left<
\begin{matrix}\eta\\[-1ex] \\[-.5ex] \varepsilon
\end{matrix}\right>\ \in\ \  \begin{matrix}\Nat(\Id, \radj F\ladj F)\\ \times\\ \Nat(\ladj F\radj F, \Id)\end{matrix}\ \ \ \Bigg| \ \ 
\begin{tikzar}[row sep=1.2em,column sep=2em]
\& \radj F \ladj F \radj F \ar{dr}[description]{\radj F\varepsilon}\\
\radj F\ar{ur}[description]{\eta\radj F} \ar[equal]{rr} \&\& \radj F
\\[-3ex]
\ladj F\ar{dr}[description]{\ladj F\eta} \ar[equal]{rr}\&\& \ladj F\\
\& \ladj F \radj F\ladj F  \ar{ur}[description]{\varepsilon \ladj F}  \end{tikzar}
 \ \ \right\}
\\[2ex]
\Adj (F, G) & = & \coprod_{\substack{H\in \CAT(\AAa,\CCc)\\
K\in \CAT(\BBb,\DDd)}}\left\{\left<\begin{matrix} \ladj \upsilon \\[-.5ex] \\[-.5ex] \radj\upsilon \end{matrix}\right>\  \in\ \ \begin{matrix}\Niso(K\ladj F ,\ladj G H)\\ \times\\\Niso(H\radj F,\radj G K) \end{matrix}\ \  \  \ \Bigg|\ \  
\begin{tikzar}[row sep=1.75em,column sep=.35em]
H\radj F\ladj F   \ar[leftrightarrow]{dr}[description]{\radj\upsilon \ladj F} 
\&
\ladj G\radj GK \ar[leftrightarrow]{dr}[description]{\ladj G\radj \upsilon} \ar{rr}[description]{\varepsilon^GK}\&\& K \ar[leftarrow]{dd}[description]{K\varepsilon^F} 
\\ 
\& \radj G K\ladj F  \ar[leftrightarrow]{dr}[description]{\radj G\ladj \upsilon} \& \ladj G H \radj F \ar[leftrightarrow]{dr}[description]{\ladj \upsilon\radj F} 
\\
H \ar{uu}[description]{H\eta^F} \ar{rr}[description]{\eta^G H} \&\& \radj G\ladj G H \& K\ladj F\radj F\ar{uu}[description]{K\varepsilon^F}
 \end{tikzar}\ \ \right\}\notag
 \eear
\caption{$\Nat(X,Y)$ are the natural transformations, $\Niso(X,Y)$ the natural isomorphisms from $X$ to $Y$}
\label{Fig:Adj}
\end{center}
\end{figure}

\subsection{From adjunctions to monads and comonads}\label{Sec:AM-AC} Lifting definitions (\ref{eq:mndp}--\ref{eq:cmnp}) from posets to categories requires imposing the commutativity conditions in Fig.~\ref{Fig:Mnd}, which were automatic in posets.
\begin{figure}[!ht]
\begin{center}
\bear
|\Mnd| & = & \coprod_{\raisebox{-0.75ex}{$\scriptstyle\AAa\in \CAT$}} \coprod_{\lft T\in \CAT(\AAa,\AAa)}\left\{\left<
\begin{matrix}\eta\\[-.25ex] \\[-.5ex] \mu
\end{matrix}\right>\ \in\ \  \begin{matrix}\Nat(\Id, \lft T)\\[-.25ex] \times\\[-.25ex] \Nat(\lft T\lft T, \lft T)\end{matrix}\ \ \ \Bigg| \ \ 
\begin{tikzar}[row sep=2em,column sep=1.5em]
\& \lft T\lft T\lft T  \ar{dl}[description]{\lft T\mu} \ar{dr}[description]{\mu\lft T} \\
\lft T\lft T \ar{r}[description]{\mu}\&\lft T \ar[equal]{d}\& \lft T\lft T \ar{l}[description]{\mu}
\\
\& \lft T \ar{ul}[description]{\eta\lft T} \ar{ur}[description]{\lft T\eta} 
 \end{tikzar}
 \ \ \right\}
\\[2.5ex] 
\Mnd \left(\lft T_{\AAa},\lft S_{\CCc}\right) & = & \coprod_{H\in \CAT(\AAa,\CCc)}\  \left\{\ \ \chi \in \Niso(H\lft T,  \lft SH)\ \  \Bigg|\ \ \begin{tikzar}[row sep=2em,column sep=3.2em]
\& H \lft T \ar[leftrightarrow]{dd}[description]{\chi} \& H\lft T\lft T \ar{l}[description]{H \mu^T} \ar[leftrightarrow]{d}[description]{\chi \lft T} 
\\
H \ar{ur}[description]{H \eta^T} \ar{dr}[description]{\eta^S H} 
\&\& \lft S H\lft T \ar[leftrightarrow]{d}[description]{\lft S \chi}
\\
\&\lft S H \& \lft S\lft S H \ar{l}[description]{\mu^S H}
\end{tikzar}\ \ \ \right\}\notag
\eear
\caption{Monads generalize closure operators \eqref{eq:mndp}.}
\label{Fig:Mnd}
\end{center}
\end{figure}
\begin{figure}[!ht]
\begin{center}
\bear
|\Cmn| & = & \coprod_{\raisebox{-0.75ex}{$\scriptstyle\BBb\in \CAT$}} \coprod_{\rgt T\in \CAT(\BBb,\BBb)} \left\{\left<
\begin{matrix}\varepsilon\\[-.25ex] \\[-.5ex] \nu
\end{matrix}\right>\ \in\ \  \begin{matrix}\Nat(\rgt T,\Id)\\[-.25ex] \times\\[-.25ex] \Nat(\rgt T,\rgt T\rgt T)\end{matrix}\ \ \ \Bigg| \ \
\begin{tikzar}[row sep=2em,column sep=1.5em]
\& \rgt T\rgt T\rgt T  \ar[leftarrow]{dl}[description]{\rgt T\nu} \ar[leftarrow]{dr}[description]{\nu\rgt T} \\
\rgt T\rgt T \ar[leftarrow]{r}[description]{\nu}\&\rgt T \ar[equal]{d}\& \rgt T\rgt T \ar[leftarrow]{l}[description]{\nu}
\\
\& \rgt T \ar[leftarrow]{ul}[description]{\varepsilon\rgt T} \ar[leftarrow]{ur}[description]{\rgt T\varepsilon} 
 \end{tikzar}
  \ \ \right\}
\\[2.5ex] 
 \Cmn \left(\rgt T_{\BBb},\rgt S_{\DDd}\right) & = & \coprod_{K\in \CAT(\BBb,\DDd)}\ \ 
\left\{ \kappa \in \Niso(K\rgt T, \rgt S K)\ \ \Bigg|\ \ \begin{tikzar}[row sep=2.5em,column sep=3.2em]
\& K \rgt T \ar{dl}[description]{K \varepsilon^T}  \ar{r}[description]{K \nu^T } \& K \rgt T\rgt T
\\
K  \&\&  
\rgt S K\rgt T \ar[leftrightarrow]{u}[description]{\kappa \rgt T} 
\\
\& \rgt S K \ar[leftrightarrow]{uu}[description]{\kappa} \ar{r}[description]{ K \nu^S} \ar{ul}{\varepsilon^S K} 
\& \rgt S\rgt S K \ar[leftrightarrow]{u}[description]{\rgt S\kappa}
\end{tikzar}\ \ \ \right\}
\notag
\eear
\caption{Comonads generalize interior operators \eqref{eq:cmnp}.}
\label{Fig:Cmn}
\end{center}
\end{figure}
The object parts of the functors generalizing \eqref{eq:huber-pos} are now
\beq\label{eq:huber-cat}
\begin{tikzar}[row sep=1ex,column sep=4em]
\Cmn 
\& \Adj \ar[thick]{l}[swap]{\AC} \ar[thick]{r}{\AM} 
\&
\Mnd
\\
\left<\ladj F\radj F,\varepsilon,\nu\right> 
\& \left<\ladj F, \radj F, \eta,\varepsilon\right> \ar[mapsto,shorten=3ex]{l} \ar[mapsto,shorten=3ex]{r} 
\& \left<\radj F\ladj F,\eta,\mu\right>
\end{tikzar}
\eeq
where the the comonad (chain) evaluation $\nu$ and the monad (cochain) evaluation $\mu$ are defined
\beq\label{eq:mnd-cmn-ev}
\mu = \left( \radj F\ladj F\radj F\ladj F \tto{\radj F \varepsilon{\ladj F }}\radj F\ladj F\right)   
\qquad \qquad\qquad 
\nu = \left(\ladj F\radj F \tto{\ladj F \eta{\radj F}} \ladj F \radj F\ladj F\radj F \right) 
\eeq
The arrow part  of \eqref{eq:huber-cat}, whereby the adjunction morphisms induce the monad and comonad morphisms, is displayed in Fig.~\ref{Fig:huber}.
\begin{figure}[!ht]
\begin{center}
\[\begin{tikzar}[row sep=3em,column sep=2em]
\&\Cmn \&\&\&\& \ar[thick]{llll}[swap]{\AC}  \Adj  \ar[thick]{rrrr}{\AM} \&\&\&\&\Mnd
\\[-5ex]
\&\&\&\&\AAa \ar{rr}[description]{H} \ar{d}[description]{\ladj F} 
\&\& 
\CCc\ar{d}[description]{\ladj G} 
\ar[Leftrightarrow,bend right=15,shorten=1.5mm,shift left=.666ex]{dll}[description]{\ladj \upsilon}
\&\& 
\AAa \ar{rr}[description]{H} \ar{dd}[description]{\lft F} 
\&\& 
\CCc\ar{dd}[description]{\lft G}
\ar[Leftrightarrow,bend right=15,shorten=1.5mm,shift left=.666ex]{ddll}[description]{\chi}
\\
\BBb \ar{rr}[description]{K} \ar{dd}[description]{\rgt F} 
\&\& 
\DDd\ar{dd}[description]{\rgt G}
\ar[Leftrightarrow,bend right=15,shorten=1.5mm,shift left=.666ex]{ddll}[description]{\kappa}
\&\&
\BBb \ar{rr}[description]{K} \ar{d}[description]{\radj F} 
\&\& 
\DDd\ar{d}[description]{\radj G}
\ar[Leftrightarrow,bend right=15,shorten=1.5mm,shift left=.666ex]{dll}[description]{\radj \upsilon}
\ar[thick,mapsto,shorten=2mm]{rr}
\&\&\hspace{1em} 
\\
\&\&\hspace{1em}\&\&\AAa \ar{rr}[description]{H} \ar{d}[description]{\ladj F} 
\ar[thick,mapsto,shorten=2mm]{ll}
\&\& 
\CCc\ar{d}[description]{\ladj G} 
\ar[Leftrightarrow,bend right=15,shorten=1.5mm,shift left=.666ex]{dll}[description]{\ladj \upsilon}
\&\& 
\AAa \ar{rr}[description]{H} \&\& \CCc
\\
\BBb \ar{rr}[description]{K} \&\& \DDd \&\& \BBb \ar{rr}[description]{K} \&\& \DDd
\end{tikzar}\]
\caption{$\kappa = \left(K\ladj F\radj F\stackrel{\ladj\upsilon\radj F}\longleftrightarrow \ladj G H \radj F\stackrel{\ladj G\radj \upsilon}\longleftrightarrow \ladj G \radj G K\right)$ and $\chi = \left(H\radj F\ladj F\stackrel{\radj\upsilon\ladj F}\longleftrightarrow \radj G K \ladj F\stackrel{\radj G\ladj \upsilon}\longleftrightarrow \radj G \ladj G H\right)$}
\label{Fig:huber}
\end{center}
\end{figure}

\subsection{From monads and comonads to adjunctions} To lift \eqref{eq:resolutions-pos} from posets to categories, we first replace the fixpoint posets in \eqref{eq:EM-pos} with the categories of algebras $\Emm \AAa T$ and coalgebras $\Emc\BBb T$, as specified in Fig.~\ref{Fig:EM-EC}. 
{\small 
\begin{figure}[!ht]
\begin{center}
\begin{align*}
|\Emm \AAa T| & =  \coprod_{a\in |\AAa|}\ \left\{\ 
\alpha\in \AAa(\lft Ta, a)\ \ \Big|\ \  
\begin{tikzar}[row sep=1em,column sep=.75em]
\lft T\lft T a \ar{rr}[description]{\mu} \ar{dd}[description]{\lft T \alpha} \&\& \lft T a\ar{dd}[description]{\alpha}\\
\& a \ar[bend right,shift left=.5ex]{ur}[description]{\eta} \ar[bend right,shift left=.6ex,equals]{dr}
\\
\lft T a \ar{rr}[description]{\alpha} \&\& a 
\end{tikzar}\ \ \ 
\right\}
& 
|\Emc \BBb T| & = \coprod_{b\in |\BBb|}\ \left\{\ 
\beta\in \BBb(b,\rgt T b)\ \ \Big|\ \  
\begin{tikzar}[row sep=1em,column sep=0.75em] 
\rgt T \rgt T  B \ar[leftarrow]{rr}[description]{\nu} \ar[leftarrow]{dd}[description]{\rgt T  \beta} \&\& \rgt T  b\ar[leftarrow]{dd}[description]{\beta}\\
\& b \ar[leftarrow,bend right,shift left=.5ex]{ur}[description]{\varepsilon} \ar[bend right,shift left=.6ex,equals]{dr}
\\
\rgt T  b \ar[leftarrow]{rr}[description]{\beta} \&\& b 
\end{tikzar}\ \ \ 
\right\}\\[3ex]
\Emm \AAa T(\alpha, \gamma) & =\ \ \ \ \ \  \ \  \left\{h\in \AAa(a,c)\ \Big|\   \begin{tikzar}[row sep=2.5em,column sep=2.5em]
\lft Ta \ar{r}[description]{\lft T h} \ar{d}[description]{\alpha} \& \lft T c \ar{d}[description]{\gamma} \\
a \ar{r}[description]{h} \& c \end{tikzar}\ \ \ 
\right\}
&
\Emc \BBb T(\beta, \delta) & =\ \ \ \  \  \ \ \ \left\{k\in \BBb(b,d)\ \Big|\   \begin{tikzar}[row sep=2.5em,column sep=2.5em]
\rgt T b \ar[leftarrow]{r}[description]{\rgt T  k} \ar[leftarrow]{d}[description]{\beta} \& \rgt T  c \ar[leftarrow]{d}[description]{\delta} \\
b \ar[leftarrow]{r}[description]{k} \& d \end{tikzar}\ \ \ 
\right\}
\end{align*} 
\caption{The final resolutions of the monad $\lft T$ and the comonad $\rgt T$}
\label{Fig:EM-EC}
\end{center}
\end{figure}}
The object parts of $\EC$ and $\EM$ are displayed in Figures ~\ref{Fig:huber-obj} 
\begin{figure}[!ht]
\begin{center}
\[\begin{tikzar}[row sep=2em,column sep=1em]
\Cmn \&\&\&\& \ar[leftarrow,thick]{llll}[swap]{\EC}  \Adj  \ar[leftarrow,thick]{rrrr}{\EM} \&\&\&\& \Mnd
\\[-1ex]
\&\&\&\Emc \BBb T \arrow[phantom]{dd}[description]{\dashv} \arrow[bend right = 13]{dd}[swap]{\raisebox{2mm}{$\Emc U T$}}
\&\& 
\AAa \arrow[phantom]{dd}[description]{\dashv}\arrow[loop, out = 135, in = 45, looseness = 4]{}[swap]{\lft T} \arrow[bend right = 13]{dd}[swap]{\raisebox{2mm}{$\Emm U T$}}
\\
\rgt T \ar[thick,mapsto,shorten=10mm]{rrr}  
\&\hspace{3em}\&\&\hspace{.2em}\& \&\hspace{.2em}\&\&\hspace{3em}\& \lft T\ar[thick,mapsto,shorten=10mm]{lll}
\\
\&\&\&\BBb  \arrow[loop, out = -45, in=-135, looseness = 6]{}[swap]{\rgt T}   
\arrow[bend right = 13]{uu}[swap]{\raisebox{-3mm}{$\Klc U T$}} 
\&\& \Emm \AAa T \arrow[bend right = 13]{uu}[swap]{\raisebox{-3mm}{$\Klm U T$}}
\end{tikzar}\]
\caption{The object part of $\EM$ and $\EC$. The $U$-functors are specified in Fig.~\ref{Fig:huber-units}.}
\label{Fig:huber-obj}
\end{center}
\end{figure}
and \ref{Fig:huber-units}. The arrow parts are in Fig.~\ref{Fig:huber-arrow}, where 
\begin{figure}[!ht]
\begin{center}

$\begin{tikzar}[row sep=3em,column sep=4em]
\BBb \ar{r}[description]{K} \ar{dd}[description]{\rgt S} 
\& 
\DDd\ar{dd}[description]{\rgt T}
\ar[Leftrightarrow,bend right=15,shorten=1.5mm,shift left=.666ex]{ddl}[description]{\kappa}
\&
\BBb \ar{r}[description]{K} \ar{d}[description]{U_{S}} 
\& 
\DDd\ar{d}[description]{U_{T}} 
\ar[Leftrightarrow,bend right=15,shorten=1.5mm,shift left=.666ex]{dl}[description]{\kappa}
\\
\&\hspace{.1em} 
\ar[thick,mapsto,shorten=2mm]{r}{\EC}  
\&
\Emc \BBb S 
\ar{r}[description]{K^\chi} 
\ar{d}[description]{U^{S}} 
\& 
\Emc \DDd T \ar{d}[description]{U^{T}}
\ar[Leftrightarrow,bend right=15,shorten=1.5mm,shift left=.666ex]{dl}[description]{\id}
\\
\BBb \ar{r}[description]{K} 
\& 
\DDd 
\&
\BBb \ar{r}[description]{K} 
\& 
\DDd 
\end{tikzar}$
\qquad
\qquad
$\begin{tikzar}[row sep=3em,column sep=4em]
\AAa \ar{r}[description]{H} \ar{d}[description]{U^S}  
\& 
\CCc\ar{d}[description]{U^T} 
\ar[Leftrightarrow,bend right=15,shorten=1.5mm,shift left=.666ex]{dl}[description]{\chi}
\&
\AAa \ar{r}[description]{H} \ar{dd}[description]{\lft S}
\& 
\CCc\ar{dd}[description]{\lft T}
\ar[Leftrightarrow,bend right=15,shorten=1.5mm,shift left=.666ex]{ddl}[description]{\chi}
\\
\Emm \AAa S \ar{r}[description]{H^\chi} \ar{d}[description]{U_S} 
\& 
\Emm \CCc T \ar{d}[description]{U_T}
\ar[Leftrightarrow,bend right=15,shorten=1.5mm,shift left=.666ex]{dl}[description]{\id}
\& \hspace{.1ex} \ar[thick,mapsto,shorten=2mm]{l}{\EM}  
\\
\CCc \ar{r}[description]{H} 
\& 
\AAa 
\&\CCc \ar{r}[description]{H} 
\& 
\AAa
\end{tikzar}$

\caption{The arrow part of the functors in Fig.~\ref{Fig:huber-obj}}
\label{Fig:huber-arrow}
\end{center}
\end{figure}
where $K^{\kappa}$ and $H^\chi$ are 
\beq\label{eq:KkHh}
\prooftree
b\tto\beta \rgt S b
\justifies
K^{\kappa} \beta  = \left(Kb \tto{K\beta} K\rgt S b \stackrel\kappa \longleftrightarrow \rgt T Kb \right)
\endprooftree
\qquad\qquad\qquad
\prooftree
\lft S a\tto \alpha a
\justifies
H^\chi \alpha  = \left(\lft T Ha \stackrel\chi \longleftrightarrow H\lft S a\tto{H\alpha} Ha \right)
\endprooftree\eeq

\para{The functors $\EC$ and $\EM$ require opposite distributivities.} 
Note that the functor components $K^{\kappa}$ and $H^{\chi}$ of $\EC(K,\kappa)$ and $\EM(H,\chi)$, displayed in Fig.~\ref{Fig:huber-arrow} and defined in \eqref{eq:KkHh}, require different directions of $\kappa$ and $\chi$, because
\begin{itemize}
\item $<K,\kappa>\in\Cmn(\rgt S, \rgt T)$ lifts $K\colon\BBb\to\DDd$ to $K^\kappa\colon\Emc\BBb S\to\Emc\DDd T$ using $\kappa:K\rgt S\to \rgt TK$, whereas 
\item $<H,\chi>\in\Mnd(\lft S, \lft T)$ lifts $H\colon\AAa\to\CCc$ to $H^\chi \colon\Emm\AAa S\to\Emm\BBb T$ using $\chi:\lft TH\to H\lft S$.
\end{itemize}
Although $\kappa$ and $\chi$ are 2-cell components of morphisms in different categories, $\Cmn$ and $\Mnd$ respectively, their images $\EC(K,\kappa)$ and $\EM(H,\chi)$ land in the same category $\Adj$. To accommodate both $\EC$ and $\EM$, the 2-cell components of the morphisms in $\Adj$, and consequently in $\Mnd$, and $\Cmn$, are required to be invertible. The functors $\KC$ and $\KM$, mentioned in the following section and described in detail later, also require opposite distributivities, as do $\KC$ and $\EC$ on one hand, and $\KM$ and $\EM$ on the other.

\begin{figure}[!ht]
\begin{center}
\[\begin{tikzar}[row sep=1.5mm,column sep=.1mm]
\&\Emc\BBb F \arrow[phantom]{dddddddd}[description]{\dashv} \arrow[bend right = 13]{dddddddd}[swap]{\raisebox{2mm}{$\Emc U F$}}
\& \&\&\& \& \&\&\& \AAa \arrow{llllllll}[swap]{\compson^{0}}
\arrow[phantom]{dddddddd}[description]{\dashv}  
\arrow[loop, out = 135, in = 45, looseness = 4,thin]{}[swap]{\lft F} 
\arrow[bend right = 13]{dddddddd}[swap]{\ladj F} 
\arrow[equal]{rrrrrrrr}
\&\hspace{3em} \&\hspace{3em}\&\& \& \&\&\&
\AAa  \arrow[loop, out = 135, in = 45, looseness = 4,thin]{}[swap]{\lft F} 
\arrow[phantom]{dddddddd}[description]{{\dashv}}   
\arrow[bend right = 13]{dddddddd}[swap]{\raisebox{2mm}{$\Emm U F$}}  
\\ 
{\scriptstyle \left(b\stackrel\beta\rightarrow \rgt Fb\right)} \ar[mapsto,bend right = 13]{dddddd}
\& \&
{\scriptstyle \left(\rgt Fy \stackrel\nu \rightarrow \rgt F \rgt F y\right)}\&
\&\& \& \&\&\&\&\hspace{3em}\&\&\&\&\& \& {\scriptstyle x} \ar[mapsto,bend right = 13]{dddddd} 
\&\& {\scriptstyle a}
\\ \hspace{3em}\\ \hspace{3em} \\ \hspace{3em} \\ \\ \\ 
{\scriptstyle b} \& \&
{\scriptstyle y} \ar[mapsto,bend right = 13]{uuuuuu}
\&\&\& \& \&\&\&\&\&\&\&\&\&\& {\scriptstyle\left(\lft F\lft Fx\stackrel \mu\rightarrow \lft Fx\right)}\&\& {\scriptstyle\left(\lft Fa \stackrel \alpha\rightarrow a\right)}\ar[mapsto,bend right = 13]{uuuuuu}
\\ 
\& \BBb \arrow[loop, out = -45, in=-135, looseness = 6,thin]{}[swap]{\rgt F} \arrow[bend right = 13]{uuuuuuuu}[swap]{\raisebox{2mm}{$\Klc U F$}}  \arrow[equal]{rrrrrrrr}
\&\hspace{3em} \&\hspace{3em}\&\hspace{3em}\& \& \&\&\& \BBb  
\arrow[loop, out = -45, in=-135, looseness = 6,thin]{}[swap]{\rgt F} 
\arrow{rrrrrrrr}[swap]{\compson_{1}}  
\arrow[bend right = 13]{uuuuuuuu}[swap]{\radj F} 
\& \hspace{3em}\&\&\&\& \&\&\&
\Emm \AAa F 
\arrow[bend right = 13]{uuuuuuuu}[swap]{\raisebox{2mm}{$\Klm U F$}}  
\end{tikzar}\]
\caption{The adjunction units $\EC\circ\AC (F)\ot F\to \EM\circ \AM (F)$ of $\AC\dashv \EC$ and $\AM\dashv \EM$}
\label{Fig:huber-units}
\end{center}
\end{figure}

\subsection{Monads and comonads are reflective in adjunctions} 
The adjunction equipment of 
\beq \AC\dashv \EC\colon \Cmn\to \Adj\qquad\qquad \mbox{and}\qquad\qquad \AM\dashv \EM\colon\Mnd\to \Adj\eeq is displayed in Fig.~\ref{Fig:huber-units}, with the \emph{comparison}\/ functors defined
\beq \label{eq:comparison}
\compson^{0}a = \left(\ladj Fa\tto{\ladj F\eta} \ladj F\radj F\ladj Fa\right)\qquad\qquad\mbox{and}\qquad\qquad \compson_{1}b = \left(\radj F\ladj F\radj F b\tto{\radj F\varepsilon} \radj F b\right) 
\eeq
The constructions underlying the functors $\AC$ and $\AM$ go back to Huber \cite{HuberP:monad}. Soon after Huber, Kleisli in  \cite{KleisliH}, and Eilenberg and Moore in \cite{Eilenberg-Moore} provided two different inverses of Huber's constructions, decomposing arbitrary (co)monads into adjunctions. The construction by Eilenberg and Moore induces the functors $\EC$ and $\EM$ above, whereas Kleisli's construction induces functors $\KC$ and $\KM$, which turned out to be adjoint to $\AC$ and $\AM$ on the other side:  
\beq\label{eq:KL-EM} \KM\dashv \AM\dashv \EM\qquad\qquad\mbox{and}\qquad \qquad  \KC\dashv \AC\dashv \EC\eeq
The (co)monad morphisms and the adjunction morphisms with respect to which Kleisli's and Eilenberg-Moore's constructions are functorial were introduced in \cite{MarandaJM:adj} and the induced 2-categories were studied in \cite{Auderset,StreetR:monads}. The details are still being studied \cite[e.g.]{VerkruysseJ:EM-Morita,Vidal:kleisli}. Lifting ideas from homological algebra, the multiple adjunctions that induce the same (co)monad have been called its \emph{resolutions}\/ \cite[Sec.~0.6]{Lambek-Scott:book}. It was established already in \cite{MarandaJM:adj} that Kleisli's resolution was initial, whereas Eilenberg-Moore's was final, meaning that the two constructions provide the adjoint reflections displayed in Fig.~\ref{Fig:AdjMndCmn}.
\begin{figure}[!ht]
\begin{center}
\begin{tikzar}[column sep = 8em]
\Cmn 
\ar[bend left = 12,phantom]{r}[description]{\scriptstyle \top}
\ar[bend right = 12,phantom]{r}[description]{\scriptstyle \top} 
\ar[bend right=20,tail]{r}[swap]{\KC} 
\ar[bend left=20,tail]{r}{\EC} 
\ar[twoheadleftarrow]{r}[description]{\AC}
\&
\Adj
\& 
\Mnd 
\ar[bend left = 12,phantom]{l}[description]{\scriptstyle \top}
\ar[bend right = 12,phantom]{l}[description]{\scriptstyle \top} 
\ar[bend right=20,tail]{l}[swap]{\EM} 
\ar[bend left=20,tail]{l}{\KM} 
\ar[twoheadleftarrow]{l}[description]{\AM}
\end{tikzar}
\caption{Both comonads and monads are reflective and coreflective in adjunctions}
\label{Fig:AdjMndCmn}
\end{center}
\end{figure}
Note that Fig.~\ref{Fig:Nuc} only displays half of the monad and comonad resolution picture\footnote{This is mostly to avoid clutter, and because the initial resolutions only really  concern us towards the end.}. 
Since both Kleisli's and Eilenberg-Moore's  constructions invert Huber's 
\beq \AC\circ \EC = \AC\circ \KC = \Id_{\Cmn}\qquad\qquad\mbox{and}\qquad\qquad \AM\circ \EM = \AM\circ \KM = \Id_{\Mnd}   \eeq
they induce idempotent monads 
\begin{align*}
\ECmn & = \EC\circ\AC & \EMnd &= \EM\circ\AM\\
\KCmn & = \KC\circ\AC & \KMnd &= \KM\circ\AM
\end{align*}
The first two are studied next, in Thm.~\ref{thm:Theorem}. For the record, note that both monads and comonads are localized as subcategories of adjunctions in \emph{two}\/ extremal ways and provide a typical example of Lawvere's \emph{unity-and-identity of the opposites} \cite{Kelly-Lawvere:taco,LawvereFW:taco,LawvereFW:taco-unity}.

\section{The nucleus functor}\label{Sec:prop}

Fig.~\ref{Fig:construction} defines the object part of the nucleus construction $\NucL$.
\begin{figure}[!ht]
\begin{center}
$\prooftree
F\  =\  \left(\adj F \ \colon \ \ \BBb\to \AAa\right)\hspace{1em}
\justifies
\NucL(F) = \left(\nadj F\colon \Emm\AAa F \to \Emc\BBb F\right)
\endprooftree$
\hspace{5em}
$\begin{tikzar}[row sep=2.5cm,column sep=3cm]
\AAa \arrow[phantom]{d}[description]{\dashv}  
\arrow[loop, out = 135, in = 45, looseness = 4,thin]{}[swap]{\lft F} 
\arrow[bend right = 13]{d}[swap]{\ladj F} 
\arrow[thin,dashed]{r}{\compson^{0}} 
\& 
\Emc\BBb F \arrow[bend right = 13,pos=0.8,thin]{dl}{\Emc V F} 
\arrow[phantom]{d}[description]{{\dashv}}   
\arrow[bend right = 13]{d}[swap]{\lnadj{F}}  
\\
\BBb  
\arrow[loop, out = -45, in=-135, looseness = 6,thin]{}[swap]{\rgt F} 
\arrow[thin,dashed]{r}[swap]{\compson_{1}}  
\arrow[bend right = 13]{u}[swap]{\radj F} 
\& 
\Emm \AAa F 
\arrow[bend right = 13]{u}[swap]{\rnadj{F}}  \arrow[bend right = 13,crossing over,thin]{ul}[swap,pos=.75]{\Klm U F} 
\end{tikzar}$
\caption{{\small The object part of $\NucL$ is comprised of $\lnadj F = \compson_1\circ \Emc V F$ and $\rnadj F = \compson^0\circ\Emm U F$ as defined in Fig.~\ref{Fig:huber-obj}}}
\label{Fig:construction}
\end{center}
\end{figure}
Fig.~\ref{Fig:construction-arrow} defines the arrow part, which is
$$\prooftree
\hspace{1em}\Upsilon\ \ \ =\ \ \ \left<H,K,\ladj\upsilon,
\radj \upsilon\right>\ \ \colon\ \  \ \ F\to G\hspace{1.5em}
\justifies
\NucL(\Upsilon) =  \left< K^{\upsilon},H^{\upsilon}, \lnadj\upsilon, \rnadj \upsilon \right>\colon \NucL(F)\to \NucL(G)
\endprooftree$$
\begin{figure}[!ht]
\begin{center}
\begin{gather*}
\prooftree
b\tto\beta \ladj F\radj F b
\justifies
K^{\upsilon} \beta  = \left(Kb \tto{K\beta} K\ladj F\radj F b \stackrel{\ladj\upsilon\radj F} \longleftrightarrow \ladj G H\radj F b \stackrel{\ladj G\radj \upsilon} \longleftrightarrow \ladj G\radj G Kb\right)
\endprooftree 
\hspace{4em} \lnadj \upsilon\  =\  \ladj G\left(\radj G\ladj \upsilon\circ\radj \upsilon \ladj F\right)
\\[2ex]
\prooftree
\radj F\ladj F a\tto \alpha a
\justifies
H^\upsilon \alpha  = \left(\radj G\ladj G Ha  \stackrel{\radj G\ladj \upsilon} \longleftrightarrow 
\radj G K\ladj F a\stackrel{\radj \upsilon \ladj F} \longleftrightarrow H \radj F\ladj F a\tto{H\alpha}Ha \right)
\endprooftree
\hspace{4em} \rnadj \upsilon\  =\  
\radj G\left(\ladj G\radj \upsilon\circ  \upsilon \radj F\right)
\end{gather*}
\\[3ex]
\begin{tikzar}[row sep=4em,column sep=8em]
\AAa \ar{r}[description]{H} \ar{d}[description]{\ladj F} 
\& 
\CCc\ar{d}[description]{\ladj G}
\ar[Leftrightarrow,bend right=15,shorten=1.5mm,shift left=.666ex]{dl}[description]{\ladj \upsilon}
\&
\Emc \BBb F \ar{r}[description]{K^{\upsilon}} \ar{d}[description]{\lnadj F} 
\& 
\Emc \DDd G\ar{d}[description]{\lnadj G} 
\ar[Leftrightarrow,bend right=15,shorten=1.5mm,shift left=.666ex]{dl}[description]{\lnadj \upsilon}
\\
\BBb \ar{r}[description]{K} \ar{d}[description]{\radj F}  
\& \DDd \ar{d}[description]{\radj G} 
\ar[Leftrightarrow,bend right=15,shorten=1.5mm,shift left=.666ex]{dl}[description]{\radj \upsilon}\ar[thick,mapsto,shorten=4mm]{r}{\NucL} 
\&
\Emm \AAa F \ar{r}[description]{H^{\upsilon}} 
\ar{d}[description]{\rnadj F} 
\& 
\Emm \CCc G \ar{d}[description]{\rnadj G}
\ar[Leftrightarrow,bend right=15,shorten=1.5mm,shift left=.666ex]{dl}[description]{\rnadj \upsilon}
\\
\AAa \ar{r}[description]{H} 
\& 
\CCc 
\&
\Emc \BBb F \ar{r}[description]{K^{\upsilon}} 
\& 
\Emc \DDd G
\end{tikzar}
\caption{{\small The arrow part of $\NucL$ maps $<H,K,\ladj \upsilon,\radj \upsilon>$ to $<K^{\upsilon},H^{\upsilon},\lnadj \upsilon,\rnadj \upsilon>$ as defined in \eqref{eq:KkHh}}}
\label{Fig:construction-arrow}
\end{center}
\end{figure}

\noindent Prop.~\ref{prop:nuc}  below claims that Figures \ref{Fig:construction} and \ref{Fig:construction-arrow} together define a functor $\NucL\colon \Adj\to \Adj$.  Before proving that, we discharge a useful lemma that follows from the definition of $\lnadj F$ and $\rnadj F$ alone.

\subsection{Nucleus decomposes descent monads and comonads}
Descent monads and comonads emerge as a common denominator of a variety of structures echoing the monad-comonad couplings in descent theory  \cite{BalmerP:annalen,Caenepeel:galois,HessK:general,MesablishviliB:comonDesc}. They will come in handy in Sections \ref{Sec:NucTheorem}
 and \ref{Sec:desc}.

\begin{definition} For an adjunction $F=\left(\adj F\colon \BBb\to \AAa\right)$, the monad $\lft F = \radj F\ladj F$ on $\AAa$ and $\rgt F = \ladj F\radj F$ on $\BBb$ with the final (Eilenberg-Moore) resolutions
\[U(F) =\left(\Emm U F\dashv \Klm U F\colon \Emm \AAa F\to \AAa\right)\qquad\qquad\qquad\qquad V(F) =\left(\Emc V F\dashv \Klc V F\colon \BBb\to \Emc \BBb F\right)\]
the induced comonad and monad \emph{of $F$-descent} are defined
\[\Rgt F = \Emm U F \Klm U F\colon \Emm \AAa F\to \Emm\AAa F\qquad\qquad\qquad \qquad \Lft F = \Klc V F \Emc V F \colon \Emc\BBb F \to \Emc \BBb F\]
\end{definition} 

\begin{lemma}\label{lemma:nuc-decom}
The descent monad $\Rgt F$ and the descent comonad $\Lft F$ also arise as composites of the nucleus components $\lnadj F = \compson_1\circ \Emc V F$ and $\rnadj F = \compson^0\circ\Emm U F$.
\beq
\Emm U F \Klm U F \ =\ \Rgt F\ =\ \lnadj F \rnadj F\qquad\qquad\qquad\qquad 
\Klc V F \Emc V F \ =\ \Lft F\ =\ \rnadj F \lnadj F
\eeq
\end{lemma}

\begin{proof} Fig.~\ref{Fig:decomp} and the definitions of $\lnadj F$ and $\rnadj F$ give
\[
\lnadj F \rnadj F \ =\ \compson_{1}\Emc V F \compson^{0}\Emm U F\ =\ \Emm U F\Klm U F\qquad\qquad\qquad\qquad
\rnadj F \lnadj F \ =\ \compson^{0}\Emm U F\compson_{1}\Emc V F \ =\ \Klc V F\Emc V F
\]
where $\compson_{1}\Emc V F \compson^{0} =\Emm U F$ and $\compson^{0}\Emm U F\compson_{1} = \Klc V F$ follow from \eqref{eq:comparison}.
\begin{figure}[!ht]
\begin{center}
\[\begin{tikzar}[row sep=2.5cm,column sep=3cm]
\AAa 
\arrow[thin]{r}{\compson^{0}} 
\& 
\Emc\BBb F \arrow[bend right = 13,pos=0.8,thin]{dl}[swap]{\Emc V F}
\arrow[leftarrow,bend left = 13,pos=0.8,thin]{dl}{\Klc V F} 
\ar[phantom]{dl}[pos=.8,description,rotate = -45]{\dashv}
\arrow[loop, out = 135, in = 45, looseness = 2.5]{}[swap]{\Lft F}
\arrow[bend right = 13]{d}[swap]{\lnadj{F}}  
\\
\BBb  
\arrow[thin]{r}[swap]{\compson_{1}}  
\& 
\Emm \AAa F 
\arrow[loop, out = -45, in=-135, looseness = 6]{}[swap]{\Rgt F} 
\arrow[bend right = 13]{u}[swap]{\rnadj{F}}  \arrow[bend right = 13,crossing over,thin]{ul}[swap,pos=.75]{\Klm U F} 
\arrow[leftarrow,bend left = 13,crossing over,thin]{ul}[pos=.75]{\Emm U F} 
\ar[phantom]{ul}[pos=.75,description,rotate = 45]{\dashv}
\end{tikzar}\]
\caption{The monad $\Lft F = \Klc V F \Emc V F$ and the comonad $\Emm U F \Klm U F \ =\ \Rgt F$ also decompose by $\rnadj F$ and $\lnadj F$ }
\label{Fig:decomp}
\end{center}
\end{figure}
\end{proof}

\subsection{Nucleus is an adjunction}
\begin{prop}\label{prop:nuc} 
The nucleus functor $\NucL\colon\Adj\to\Adj$ is well-defined:
\begin{itemize}
\item $\NucL(F)$, specified in Fig.~\ref{Fig:construction}, is an adjunction; and
\item $\NucL(H,K,\ladj \upsilon,\radj \upsilon)$ in Fig.~\ref{Fig:construction-arrow} is an adjunction morphism.
\end{itemize}
\end{prop}

\begin{proof}
Checking that $\left< K^{\upsilon},H^{\upsilon}, \lnadj\upsilon, \rnadj \upsilon \right>$ satisfies the two commutativity requirements of adjunction morphisms in Fig.~\ref{Fig:Adj}, as soon as $\left<H, K, \ladj\upsilon, \radj \upsilon \right>$ satisfies them, is straightforward. Checking the adjunction equipment of $\nadj F$ requires work. The object parts of the functors $\rnadj F$ and $\lnadj F$ from Fig.~\ref{Fig:construction} are unfolded in Fig.~\ref{Fig:nucdef} again. 
\begin{figure}[!ht]
\begin{center}
$\begin{tikzar}[row sep=.07cm,column sep=1cm]
{\scriptstyle x}\arrow[mapsto,dashed,thin]{r} 
 \& \left(\begin{minipage}[c]{3em}\scriptsize\centering$\ladj F x\hspace{1.5em}$\\ $\downarrow \ladj F \eta$\\ 
$\ladj F \radj F \ladj F x$ \end{minipage}\right)  
\& \left(\begin{minipage}[c]{2em}\scriptsize\centering$\radj F \ladj F x$\\ $\downarrow \alpha$\\ 
$x\hspace{1em}$ \end{minipage}\right) \arrow[mapsto]{l} 
 \\ \vspace{5ex} 
\arrow[dashed,thin]{r}{\compson^0} \AAa \&  \Emc \BBb F \&  \Emm \AAa F \arrow[thick]{l}[swap]{\rnadj{F}} \ar[bend left = 30,thin]{ll}[pos=0.8]{\radj U}
 \\[8ex] 
\BBb \arrow[dashed,thin]{r}{\compson_1} \& \Emm \AAa F \&  \arrow[thick]{l}[swap]{\lnadj{F}} \Emc \BBb F\ar[bend right = 30,thin]{ll}[swap,pos=.8]{\ladj V} \\
{\scriptstyle y} \arrow[mapsto,dashed,thin]{r} 
\&
\left(\begin{minipage}[c]{3em}\scriptsize\centering$\radj F \ladj F\radj F y$\\ $\downarrow \radj F \varepsilon$\\ 
$\radj F y\hspace{2em}$ \end{minipage}\right)  \& 
\left(\begin{minipage}[c]{2em}\scriptsize\centering$y\hspace{1em}$\\ $\downarrow \beta$\\ 
$\ladj F \radj F y$
 \end{minipage}\right) 
\arrow[mapsto]{l} 
\end{tikzar}$
\caption{The definitions of $\rnadj F$ and $\lnadj F$}
\label{Fig:nucdef}
\end{center}
\end{figure}
The arrow part of $\rnadj F$ is the same as $\ladj F$ and the arrow part of $\lnadj F$ is the same as $\radj F$. For the functors $\lnadj F$ and $\rnadj F$, we now prove that the correspondence
\begin{align}\label{eq:bigadjj}
\Emm \AAa F \big( \lnadj{F} \beta, \alpha\big) \ & \ \cong\ \   \Emc \BBb F \big( \beta,  \rnadj{F}\alpha\big) \\
f & \mapsto  \overline f = \ladj F f\circ \beta \nonumber
\end{align}
is a natural bijection. More precisely, the claim is that
\begin{enumerate}[a)]
\item $f$ is an algebra homomorphism if and only if $\overline f$ is a coalgebra homomorphism: each of the following squares commutes if and only if the other one commutes
\bea\label{eq:bigadj}
\begin{tikzar}[row sep=2em,column sep=1em]
\radj F\ladj F \radj F y \ar{rr} {\radj F\ladj F f} \ar{dd} {=\lnadj{F}\beta}[swap]{\radj F \varepsilon} \& \& \radj F \ladj F x \ar{dd} {\alpha}\\ \& \hspace{1em}  \\
\radj F y \ar{rr}[swap]f \&\& x
\end{tikzar}
\hspace{1.5em}& \iff &\hspace{1.5em}
\begin{tikzar}[row sep=2em,column sep=1em]
\ladj F \radj F y \ar{rr} {\ladj F \radj F \overline  f}  \&\&  \ladj F \radj F \ladj F x \\ \& \hspace{1em} \\
y \ar{rr}[swap]{\overline  f} \ar{uu} {\beta} \&\& \ladj F x \ar{uu} {\rnadj{F}\alpha = }[swap]{\ladj F \eta}
\end{tikzar}
\eea
\item the map $f\mapsto \overline f$ is a bijection, natural along the coalgebra homomorphisms on the left and along the algebra homomorphisms on the right.
\end{enumerate}
Claim (a) is proved as Lemma~\ref{lemma:a}.  The bijection part of claim (b) is proved as Lemma~\ref{lemma:b}. The naturality of the correspondence is straightforward. \end{proof}

\begin{lemma}
\label{lemma:b}
For an arbitrary adjunction $F = \left(\adj F: \BBb\to \AAa\right)$, any algebra $\radj F \ladj F x\tto \alpha x$, and any coalgebra $y \tto \beta \ladj F \radj F y$ in $\BBb$, the mappings
\[
\begin{tikzar}[row sep=3.5em,column sep=4.5em]
\AAa(\radj F y, x) \arrow[bend left = 20]{r}{\overline{(-)}}  
\&  \BBb(y, \ladj F x) \arrow[bend left = 20]{l}{\underline{(-)}}
\end{tikzar}
\]
defined by
\[ \overline{f} = \ladj F f \circ \beta \qquad \qquad \qquad \underline g = \alpha \circ \radj F g
\]
induce a bijection between the subsets
\bear
\big\{f\in \AAa(\radj F y, x)\ |\ f = \alpha\circ \radj F \ladj F f \circ \radj F \beta \big\} & \cong & \big\{g\in \BBb(y,\ladj F x)\ |\ g = \ladj F \alpha\circ \ladj F \radj F g \circ \beta \big\}
\eear
illustrated in the following diagram.
\bear
\begin{tikzar}[row sep=2em,column sep=1.6em]
\radj F\ladj F \radj F y \ar{rr} {\radj F\ladj F f} 
\&\& \radj F \ladj F x \ar{dd} {\alpha}\\ 
\\
\radj F y \ar{rr}[swap]f \ar{uu}{\radj F \beta}
\ar[thin]{rruu}[swap]{\radj F g}
\&\& x
\end{tikzar}
& 
\leftrightsquigarrow  & 
\begin{tikzar}[row sep=2em,column sep=1.6em]
\ladj F \radj F y \ar{rr} {\ladj F \radj F g}  \ar[thin]{rrdd}[swap]{\ladj F f}
\&\&  \ladj F \radj F \ladj F x \ar{dd}{\radj F {\alpha}}\\ 
\\
y \ar{rr}[swap]{g} \ar{uu}{\beta} \&\& \ladj F x 
\end{tikzar}
\eear
\end{lemma}

\begin{proof}
Following each of the mappings "there and back" gives
\begin{gather*}
f\hspace{2em} \longmapsto\hspace{2em} \overline f =\ladj F f \circ \beta\hspace{2em} \longmapsto\hspace{2em} \underline{\overline f} = \alpha \circ \radj F\ladj F f \circ \radj F \beta = f\\
g\hspace{2em} \longmapsto\hspace{2em} \underline g = \alpha \circ \radj F g\hspace{2em} \longmapsto\hspace{2em} \overline{\underline g} = \ladj F \alpha \circ \ladj F\radj F g  \circ \beta = g
\end{gather*}
\end{proof}

\begin{lemma}
\label{lemma:a}
For any adjunction $F = \left(\adj F: \BBb\to \AAa\right)$, algebra $\radj F \ladj F x\tto \alpha x$ in $\AAa$, coalgebra $y \tto \beta \ladj F \radj F y$ in $\BBb$, arrow $f\in \AAa(\radj F y, x)$ and $\overline f = \ladj F f\circ \beta \in \BBb(x,\ladj F y)$, if any of the squares (1-4) in Fig.~\ref{Fig:squares} commutes, then they all commute.
\begin{figure}[!ht]
\begin{center}
\bear
\begin{tikzar}[row sep=1.5em,column sep=1em]
\radj F\ladj F \radj F y \ar{rr}{\radj F\ladj F f} \ar{dd} [swap]{\radj F \varepsilon} \&\& \radj F \ladj F x \ar{dd} {\alpha}\\
\& (1)\\
\radj F y \ar{rr}[swap]f \&\& x
\end{tikzar}
\hspace{2em}&  &\hspace{2em}
\begin{tikzar}[row sep=1.5em,column sep=1em]
\ladj F \radj F y \ar{rr} {\ladj F \radj F \overline  f}  \&\&  \ladj F \radj F \ladj F x \\ \& (4)\\
y \ar{rr}[swap]{\overline  f} \ar{uu} {\beta} \&\& \ladj F x \ar{uu} [swap]{\ladj F \eta}
\end{tikzar}\\[2ex]
(a)\Updownarrow\hspace{6em} & & \hspace{6em}\Updownarrow(c)\\[2ex]
\begin{tikzar}[row sep=1.5em,column sep=1em]
\radj F\ladj F \radj F y \ar{rr} {\radj F\ladj F f} \&\& \radj F \ladj F x \ar{dd} {\alpha}\\ \& (2)\\
\radj F y \ar{rr}[swap]f \ar{uu}{\radj F \beta}\&\& x
\end{tikzar}
\hspace{2em}& \stackrel{\displaystyle(b)}\Leftrightarrow  & \hspace{2em}
\begin{tikzar}[row sep=1.5em,column sep=1em]
\ladj F \radj F y \ar{rr} {\ladj F \radj F \overline  f}  \&\&  \ladj F \radj F \ladj F x \ar{dd}{\radj F {\alpha}}\\ \& (3) \\
y \ar{rr}[swap]{\overline  f} \ar{uu}{\beta} \&\& \ladj F x 
\end{tikzar}
\eear
\caption{Proof schema for \eqref{eq:bigadj}}
\label{Fig:squares}
\end{center}
\end{figure}
In particular, a square on one side of any of the equivalences  (a--c) commutes if and only if the square on the other side of the equivalence commutes. 
\end{lemma}

\begin{proof}
The claims are established as follows.

\paragraph{$(1)\stackrel{(a)}\Rightarrow (2)$:} Using the commutativity of (1) and $(\ast)$ the counit equation $\varepsilon \circ \beta = \id$ for the coalgebra $\beta$, we derive (2) as 
\[\alpha\circ \radj F\ladj F f \circ \radj F \beta\ \  \stackrel{(1)} =\ \  f\circ \radj F \varepsilon\circ \radj F \beta\ \  \stackrel{(\ast)}=\ \  f\]

\paragraph{$(2)\stackrel{(a)}\Rightarrow (1)$} is proved by chasing the following diagram:
\[
\begin{tikzar}[row sep=2em,column sep=1.3em]
\radj F\ladj F \radj F \ladj F \radj F y \ar{rrrrrr}{\radj F \ladj F\radj F\ladj F f} \ar{dddddd}[swap]{\radj F\varepsilon} \& \& \& \& \& \& \radj F\ladj F \radj F \ladj F x \ar{ddll}{\radj F \ladj F \alpha} \ar{dddddd}{\radj F\varepsilon}\\
\& \& \& (2) \\
\&\& \radj F \ladj F \radj F y \ar{dd}{\radj F \varepsilon} \ar{uull}[swap]{\radj F \ladj F \radj F \beta} \ar{rr}{\radj F \ladj F f} \& \& \radj F \ladj F x \ar{dd}[swap]{\alpha}\\
\&(\dag) \& \& (1) \& \& (\ddag) \\
\&\& \radj F y\ar{rr}[swap]f \ar{ddll}{\radj F \beta} \&\& x \\
\& \& \& (2) \\
\radj F \ladj F \radj F y \ar{rrrrrr}[swap]{\radj F \ladj F f} \& \& \& \& \& \& \radj F \ladj F x \ar{uull}[swap]\alpha
\end{tikzar}
\]
The top and the bottom trapezoids commute by assumption (2), whereas the left-hand trapezoid (denoted (\dag)) and the outer square (denoted $(\Box)$) commute by the naturality of $\varepsilon$. The right-hand trapezoid (denoted (\ddag)) commutes by the cochain condition for the algebra $\alpha$. It follows that the inner square (denoted (1)) must also commute:
\bear
f\circ \radj F \varepsilon & \stackrel{(2)} = & \alpha \circ \radj F \ladj F f \circ \radj F \beta\circ \radj F \varepsilon \\
& \stackrel{(\dag)} = & \alpha \circ \radj F \ladj F f \circ  \radj F \varepsilon \circ \radj F \ladj F \radj F \beta \\
& \stackrel{(\Box)} = & \alpha \circ  \radj F \varepsilon  \circ \radj F \ladj F \radj F \ladj F f \circ \radj F \ladj F \radj F \beta\\
& \stackrel{(\ddag)} = & \alpha \circ  \radj F\ladj \alpha  \circ \radj F \ladj F \radj F \ladj F f \circ \radj F \ladj F \radj F \beta\\
& \stackrel{(2)} = & \alpha \circ   \radj F \ladj F f 
\eear

\paragraph{$(4)\stackrel{(c)}\Leftrightarrow (3)$} is proven dually to $(1)\stackrel{(a)}\Leftrightarrow (2)$ above. The duality consists of reversing the arrows, switching $\radj F$ and $\ladj F$, and also $\alpha$ and $\beta$, and replacing $\varepsilon$ with $\eta$.

\paragraph{$(2) \stackrel{(b)} \Leftrightarrow (3)$} follows from Lemma~\ref{lemma:b}.
\end{proof}

\section{The nucleus monad is idempotent}\label{Sec:Theorem}

The claim is that the nucleus of a nucleus is equivalent to it. An easy way to prove this is to show that the construction of the nucleus $\NucL(F)$ can be decomposed into the final (Eilenberg-Moore) resolutions of the monad and the comonad induced by the adjunction $F=\left(\adj F\right)$. More precisely, first resolving the monad $\lft F$ induced by $F$ into $\EMnd(F)$ and then resolving the comonad $\Rgt F$ induced by $\EMnd(F)$ to $\ECmn\left(\EMnd(F)\right)$ is the same as first resolving $F$ to $\ECmn(F)$ and then to $\EMnd\left(\ECmn(F)\right)$. The two pairs of resolutions are displayed in Fig.~\ref{Fig:NucL-units}.
\begin{figure}[!ht]
\begin{center}
\[\begin{tikzar}[row sep=1.2mm,column sep=6mm]
\Emc\BBb F \arrow[loop, out = 135, in = 45, looseness = 4,thin]{}[swap]{\Lft F}  \arrow[bend right = 13,pos=0.51]{dddddddd}
\& \&\&\& \Emc\BBb F \arrow[loop, out = 135, in = 45, looseness = 4,thin]{}[swap]{\Lft F}  \arrow[equal]{llll}
 \arrow[bend right = 13]{dddddddd}
\& \&\&\& \AAa \arrow{llll}[swap]{\compson^{0}}
\arrow[loop, out = 135, in = 45, looseness = 4,thin]{}[swap]{\lft F} 
\arrow[bend right = 13]{dddddddd}\arrow[equal]{rrrr}
\&\&\&\& \AAa 
\arrow[bend right = 13]{dddddddd}
\arrow[loop, out = 135, in = 45, looseness = 4,thin]{}[swap]{\lft F} \arrow{rrrr}[swap]{\compstwo_{0}} 
 \& \&\&\&
\EMC \AAa F  
\arrow[bend right = 13]{dddddddd}
\\ \\ \\ \\ \\
\hspace{0.1ex}
\&\&\&\&\hspace{0.1ex}\ar[Rightarrow,shorten=12mm]{llll}[swap]{\eta} \&\&\&\&\hspace{0.1ex}\ar[Rightarrow,shorten=10mm]{llll}[swap]{\eta} \ar[Rightarrow,shorten=10mm]{rrrr}{\eta} \&\&\&\&\hspace{0.1ex} \ar[Rightarrow,shorten=12mm]{rrrr}{\eta}  \&\&\&\&\hspace{0.1ex}
 \\ \\ \\ 
\EMM\BBb F 
\arrow[bend right = 13,pos=0.45]{uuuuuuuu} \arrow[phantom]{uuuuuuuu}[description,pos=0.4]{\EMnd\circ \ECmn(F)\ \ } 
\&\&\&\& 
\BBb \arrow[loop, out = -45, in=-135, looseness = 6,thin]{}[swap]{\rgt F} \arrow[bend right = 13]{uuuuuuuu}
\arrow[equal]{rrrr}
\arrow{llll}{\compstwo^{1}}
\arrow[phantom]{uuuuuuuu}[description]{\ECmn(F)}
\& \&\&\& \BBb  
\arrow[loop, out = -45, in=-135, looseness = 6,thin]{}[swap]{\rgt F} 
\arrow{rrrr}[swap]{\compson_{1}}  
\arrow[bend right = 13]{uuuuuuuu}
\arrow[phantom]{uuuuuuuu}{F}
\&\&\&\& \Emm \AAa F \arrow[equal]{rrrr} \arrow[loop, out = -45, in=-135, looseness = 6,thin]{}[swap]{\Rgt F} \arrow[bend right = 13]{uuuuuuuu}
\arrow[phantom]{uuuuuuuu}[description]{\EMnd(F)}   
\& \&\&\&
\Emm \AAa F \arrow[loop, out = -45, in=-135, looseness = 6,thin]{}[swap]{\Rgt F}
\arrow[bend right = 13]{uuuuuuuu}  \arrow[phantom]{uuuuuuuu}[description]{\ \ \ECmn\circ \EMnd(F)}     
\end{tikzar}\]
\caption{Alternating the monad and comonad localizations of an adjunction $F$}
\label{Fig:NucL-units}
\end{center}
\end{figure}
Now we claim that the adjunctions $\EMnd\left(\ECmn(F)\right)$ and $\ECmn\left(\EMnd(F)\right)$ on the two ends are equivalent, and that they are equivalent to $\NucL(F)$. Since $\ECmn$ and $\EMnd$ are idempotent monads, the claim that they commute implies that their composite is also an idempotent monad, which is also their intersection in the lattice of idempotent monads. This composite is the nucleus monad $\NucL$.

\subsection{From weak to strong equivalences through absolute completions}
The claim that the resolutions $\ECmn\circ\EMnd(F)$ and $\EMnd\circ\ECmn(F)$ are equivalent is true either 
\begin{enumerate}[a)]
\item for arbitrary categories and weak (Morita-)equivalences, or 
\item for absolutely (Cauchy-)complete categories and strong equivalences.
\end{enumerate}
See \cite[Vol.~1, Sec.~7.9]{BorceuxF:handbook} for definitions of weak equivalences and absolute completions. In the setting from (a), nuclei could be studied entirely in terms of matrices, without getting into adjunctions, monads, or comonads\footnote{Grothendieck's generalization of Galois theory \cite[Expos\'{e} VI]{GrothendieckA:fibrations59,GrothendieckA:SGA1} was originally stated in that way.  The theory of monadicity was developed as a more salient version  \cite{BeckJ:LNM80,BeckJ:thesis,Benabou-Roubaud}. Turning the tables, descent could be viewed as monadicity without monads \cite{PavlovicD:Como}.}. But the technical overhead of working with weak equivalences is significant and many simple ideas get obscured by irrelevant implementation details. Approach (b) is easier, since the absolute completions of categories are simple, succinct, always available, and make weak equivalences strong.\footnote{The tradeoffs like (a) vs (b) often arise in mathematics. Scientists could work with the sequences of rational numbers that arise from their measurements and use weak equality to compare them. Or they can construct the Cauchy-completion of the rationals, the real numbers, and work with the strong equality of those. The latter approach is usually preferred.} Appendix~\ref{appendix:idempot} provides an overview.  To get a simple strong equivalence between a nucleus and its nucleus, we pursue approach (b) and construct them over absolutely complete categories. This, however, does not mean that the result needs to be restricted to absolutely complete categories. Since the absolute completion functor $\Kar{(-)}\colon \Cat\to \Cat$ is an idempotent monad itself (as explained in Appendix~\ref{appendix:idempot}), which readily lifts to adjunctions (and also to monads and comonads), we can precompose all functors in sight with the absolute completion monad and work with absolutely complete categories with no loss of generality. In particular, we precompose the monads $\ECmn, \EMnd,\NucL\colon \Adj \to \Adj$ with the absolute completion monad $\Kar{(-)}\colon \Adj\to\Adj$, lifted from categories to adjunctions, and define functors $\EECmn, \EEMnd,\NNucL\colon \Adj \to \Adj$ by
\beq
\EECmn(F) = \ECmn(\Kar F)\qquad \qquad\NNucL(F) = \NucL(\Kar F) \qquad\qquad \EEMnd(F) = \EMnd(\Kar F)
\eeq
where $\Kar\CCc$ denotes the absolute completion of a category $\CCc$ and $\Kar F\colon \Kar\CCc\to\Kar\DDd$ is the lifting of the functor $F\colon \CCc\to\DDd$. Checking that $\EEMnd$ and $\EECmn$ are idempotent monads is straightforward, and checking that $\NNucL$ is strongly equivalent to both of their composites and thus the idempotent monad at their intersection, will be the task of the theorem below. Since the liftings are unique, the underlining is omitted whenever the confusion is unlikely. We reiterate that this is inessential for the main constructions and could be avoided at the cost of more involved arguments.

\subsection{The Nucleus Theorem}\label{Sec:NucTheorem}
\begin{thm}\label{thm:Theorem} 
The idempotent monads
\beq\label{eq:ECmn-EMnd-defn} 
\EECmn = \left(\Adj\eepi{\Kar\AC}\Cmn\mmono{\EC} \Adj\right) \qquad\mbox{and}\qquad \EEMnd = \left(\Adj\eepi{\Kar\AM}\Mnd\mmono{\EM} \Adj\right)   
\eeq
commute and their composites are isomorphic to the nucleus 
\beq\label{eq:ECmn-EMnd} \EECmn \circ \EEMnd(F) \cong \NNucL(F) \cong  \EEMnd \circ \EECmn(F)   
\eeq
naturally in $F$. With the leftmost and the rightmost adjunctions in Fig.~\ref{Fig:NucL-units} isomorphic, the monad units $\eta$ displayed there can be folded into the commutative square in $\Adj$, displayed in  Fig.~\ref{Fig:Nuc-fac}, presenting the unit of the monad $\NNucL$.
\begin{figure}[!ht]
\begin{center}
\begin{tikzar}[row sep = 2.5em,column sep = 5em]
\&\Kar \AAa \ar[bend right = 20,thin]{d} 
\ar[leftarrow,bend left = 20,thin]{d}\ar[phantom]{d}[description,pos=0.55]{\mbox{\scriptsize\ \ $\EEMnd(F)$}} 
\ar{ddrr} {\compson^0} 
\\
\& \Emma \AAa F 
\ar[equal]{ddrr}
\\
\Kar\AAa \ar[leftarrow,bend left = 20,thin]{d} \ar[equal]{uur} \ar[bend right = 20,thin]{d}
\ar[phantom]{d}[description]{\mbox{\ \footnotesize $\Kar F$\ }} 
\&\&\& 
\Emca \BBb F \ar[bend right = 20,thin]{d}
\ar[leftarrow,bend left = 20,thin]{d} \ar[phantom]{d}[description,pos=0.55]{\mbox{\ \scriptsize $\NNucL(F)$\ }}
\\
\Kar\BBb 
\ar[equal]{ddrr} \ar{uur}[swap]{\compson_1}
\&\&\& 
\Emma \AAa F 
\\
\&\& \Emca\BBb F \ar[equal,crossing over]{uur} \ar[leftarrow,crossing over]{uull}[swap]{\compson^0}
\ar[bend right = 20,thin]{d}
\ar[leftarrow,bend left = 20,thin]{d}
\ar[phantom]{d}[description,pos=0.45]{\mbox{\scriptsize $\EECmn(F)$\,}}
\\
\&\& \Kar \BBb 
 \ar{uur}[swap]{\compson_1}
\end{tikzar}
\caption{$\NNucL(F)$ factorized into $\EEMnd\circ\EECmn(F)$ and $\EECmn\circ\EEMnd(F)$}
\label{Fig:Nuc-fac}
\end{center}
\end{figure}
\end{thm}

\begin{proof} The claim is that the leftmost resolution $\EEMnd\circ\EECmn(F)$ and the rightmost resolution $\EECmn\circ\EEMnd(F)$ in Fig.~\ref{Fig:NucL-units} are isomorphic, and that they are both isomorphic to $\NNucL(F)$, as displayed in Fig.~\ref{Fig:NucL-commut}.
\begin{figure}[!ht]
\begin{center}
\[\begin{tikzar}[row sep=1.4mm,column sep=6mm]
\Emca\BBb F \arrow[loop, out = 135, in = 45, looseness = 4,thin]{}[swap]{\Lft F} \arrow[phantom]{dddddddd}[description]{\dashv} \arrow[bend right = 13]{dddddddd}[swap,pos=0.495]{\raisebox{-.75ex}{$\scriptstyle\EmM U F$}}
\& \&\&\& \& \&\&\& \Emca\BBb F \arrow[equal]{llllllll} \arrow[loop, out = 135, in = 45, looseness = 4,thin]{}[swap]{\Lft F}
\arrow[phantom]{dddddddd}[description]{\dashv}  
\arrow[bend right = 13]{dddddddd}[swap]{\lnadj F} 
\arrow{rrrrrrrr}[description]{\mbox{\Huge$\simeq$}}
\&\&\&\& \& \&\&\&
\EMCa\AAa F   
\arrow[phantom]{dddddddd}[description]{{\dashv}}   
\arrow[bend right = 13]{dddddddd}[swap]{\EmC U F}  
\\ \\ \\  \\  
\\ \\ \\ \\ 
\raisebox{-1ex}{$\EMMa \BBb F$}  \arrow[bend right = 13]{uuuuuuuu}[swap]{\raisebox{-2ex}{$\scriptstyle\KlM U F$}} 
\& \&\&\& \& \&\&\& \Emma \AAa F \arrow{llllllll}[description]{\mbox{\Huge$\simeq$}}  
\arrow[loop, out = -45, in=-135, looseness = 6,thin]{}[swap]{\Rgt F} 
\arrow[equal]{rrrrrrrr}  
\arrow[bend right = 13]{uuuuuuuu}[swap]{\rnadj F} 
\&\&\&\&\& \&\&\&
\Emma \AAa F \arrow[loop, out = -45, in=-135, looseness = 6,thin]{}[swap]{\Rgt F} 
\arrow[bend right = 13]{uuuuuuuu}[swap]{\KlC U F}  
\end{tikzar}\]
\caption{The equivalences claimed in Thm.~\ref{thm:Theorem}}
\label{Fig:NucL-commut}
\end{center}
\end{figure}
Recall from Lemma~\ref{lemma:nuc-decom} that $\NNucL(F) = \left(\nadj F\colon \Emma \AAa F\to \Emca \BBb F\right)$ as a resolution of the descent comonad $\Rgt F$ on $\Emma \AAa F$.  We show that for the induced comparison functor $\compss^0\colon \Emca \BBb F \to \EMCa \AAa F$ is an equivalence. $\compss^0$ is defined in Fig.~\ref{Fig:Eleft} by instantiating the standard comparison functor construction from \eqref{eq:comparison} to $\Rgt F$.
\begin{figure}[!ht]
\begin{center}
\bear\Emca \BBb F \hspace{2em}  & \tto{\hspace{2.5em} \displaystyle \compss^0 \hspace{3.5em}}\hspace{-3.5em}& \hspace{7em} \EMCa \AAa F\\[3ex] 
\begin{tikzar}{}
y \arrow{dd}{\displaystyle\beta} \\ \\ 
\ladj F \radj F y
\end{tikzar}
\hspace{1.5em}& \longmapsto &
\begin{tikzar}[row sep=2.5em,column sep=1.3em]
\radj F \ladj F \radj F y \ar{dd}[swap] {\radj F\ladj F \radj F \beta}
\ar{rrrr}[name = U,below]{\displaystyle\lnadj F\beta}{
\radj F \varepsilon_y\ =}  
\& \& \&\& \radj F y  \ar{dd}{\radj F \beta} \\
 \\
\radj F \ladj F \radj F \ladj F \radj F y  \ar{rrrr}[name = D]{\displaystyle\lnadj F\rnadj F\lnadj F\beta}[swap]{
= \ \radj F \varepsilon_{\rgt F y}} \&\&\&\&  \radj{F}\ladj F \radj F y
\arrow[thick,to path = (U) -- (D)\tikztonodes]{}{\displaystyle =\, \lnadj F \eta_{\beta}}[swap]{\displaystyle \compss^0 \beta}
\end{tikzar}\eear
\caption{The comparison functor from $\NNucL(F)$ to $\EECmn\circ\EEMnd(F)$}
\label{Fig:Eleft}
\end{center}
\end{figure}
The right adjoint $\compss_0$ is defined in  Fig.~\ref{Fig:Eright},
\begin{figure}[!ht]
\begin{center}
\bear\EMC \AAa F \hspace{3em}  & \tto{\hspace{1.5em} \displaystyle \compss_0 \hspace{1.5em}}& \hspace{2.7em} \Emc \BBb F\\[3ex] 
\begin{tikzar}[row sep=3.5em,column sep=1em]
\radj F \ladj F x \ar{dd}[swap] {\radj F\ladj F d}
\ar{rrrr}[name = U,below]{\displaystyle\alpha}  
\& \& \&\& x  \ar{dd}{d} \\
 \\
\radj F \ladj F \radj F \ladj F x  \ar{rrrr}[name = D]{\displaystyle\lnadj F\rnadj F \alpha}[swap]{= \ \radj F \varepsilon_{\ladj F x}} \&\&\&\&  \radj{F}\ladj F x
\arrow[thick,to path = (U) -- (D)\tikztonodes]{}{\displaystyle \delta}
\end{tikzar}
\hspace{1em}& \longmapsto &\hspace{1em}
\begin{tikzar}[row sep=4em,column sep=2em]
y  \arrow[dashed]{dd}[swap]{\displaystyle \compss_0 \delta} \arrow[tail]{r}{e} \& \ladj F x  \arrow{dd}{=\ladj F \eta}[swap]{\displaystyle \rnadj F \alpha}
\arrow[tail,shift right = 1ex]{rr}[swap]{\ladj F \eta} \arrow[shift left = 1ex]{rr}{\ladj F d} \arrow[shift left = 1ex,bend left = 35,two heads]{l}{r} \&\& \ladj F \radj F\ladj F x  \arrow{dd}{\ladj F \eta}[swap]{\displaystyle \rnadj F \lnadj F \rnadj F \alpha = } \arrow[shift right = 1ex,bend right = 35,two heads]{ll}[swap]{\varepsilon}\\ \\ 
\ladj F \radj F y \arrow[tail]{r}{\ladj F \radj F e} \& \ladj F\radj F\ladj F x  
\arrow[shift left = 1ex]{rr}{\ladj F\radj F\ladj F d} \arrow[tail,shift right = 1ex]{rr}[swap]{\ladj F\radj F\ladj  F \eta} \arrow[shift left = 1ex,bend left = 35,two heads]{l}{\ladj F \radj F r} \&\& \ladj F \radj F\ladj F \radj F\ladj F x  \arrow[shift right = 1ex,bend right = 35,two heads]{ll}[swap]{\ladj F \radj F \varepsilon} 
\end{tikzar}
\eear
\caption{The right adjoint of the comparison from $\NNucL(F)$ to $\EECmn\circ\EEMnd(F)$}
\label{Fig:Eright}
\end{center}
\end{figure}
where $d$ is the structure map of the coalgebra $\alpha\tto d \lnadj F\rnadj F \alpha$ in $\Emma \AAa F$ and $y$ is defined by splitting the idempotent $\varepsilon \circ \ladj F d$. \textbf{\emph{This is where the absolute completeness of $\BBb$ is used.}} To show that  $\compss^0$ and $\compss_0$ make  $\EMCa \AAa F$ and $\Emca \BBb F$  equivalent, we construct natural isomorphisms $\compss_0 \compss^0 \cong \Id$ and $\compss^0 \compss_0 \cong \Id$.  Towards $\compss_0 \compss^0 \cong \Id$, note that instantiating $\compss^0\beta:\lnadj F \beta \to \lnadj F\rnadj F\lnadj F \beta$ (the right-hand square in Fig.~\ref{Fig:Eleft}) as $\delta:\alpha\to \lnadj F\rnadj F \alpha$ (the left-hand square in Fig.~\ref{Fig:Eright}) reduces the right-hand equalizer in Fig.~ \ref{Fig:Eright} to Fig.~\ref{Fig:beta}.
\begin{figure}[!ht]
\begin{center}
\[\begin{tikzar}[row sep=3em,column sep=2em]
y  \arrow[dashed]{dd}{\displaystyle = \beta}[swap]{\displaystyle \compss_0 \compss^0 \beta} \arrow[tail]{rr}[swap]{\beta} \&\& \ladj F \radj F y  \arrow{dd}{\ladj F \eta}
\arrow[tail,shift right = 1ex]{rr}[swap]{\ladj F \eta} \arrow[shift left = 1ex]{rr}{\ladj F \radj F \beta} \arrow[shift left = 1ex,bend left = 35,two heads]{ll}{\varepsilon} \&\& \ladj F \radj F\ladj F \radj F y  \arrow{dd}
{\ladj F \eta} \arrow[shift right = 1ex,bend right = 35,two heads]{ll}[swap]{\varepsilon}\\ \\ 
\ladj F \radj F y \arrow[tail]{rr}[swap]{\ladj F \radj F \beta} 
\&\& \ladj F\radj F\ladj F \radj F y  
\arrow[shift left = 1ex]{rr}{\ladj F\radj F\ladj F \radj F \beta} \arrow[tail,shift right = 1ex]{rr}[swap]{\ladj F\radj F\ladj  F \eta} \arrow[shift left = 1ex,bend left = 35,two heads]{ll}{\ladj F \radj F\varepsilon} \&\& \ladj F \radj F\ladj F \radj F y  \arrow[shift right = 1ex,bend right = 35,two heads]{ll}[swap]{\ladj F \radj F \varepsilon} 
\end{tikzar}
\]
\caption{$\compss_0 \compss^0 \cong \id$}
\label{Fig:beta}
\end{center}
\end{figure}
The commutativity of Fig.~\ref{Fig:beta} is equivalent to the fact that $\beta$ is a $\rgt F$-coalgebra \cite[Sec.~3.6]{BarrM:ttt}. To construct the isomorphism $\compss^0 \compss_0 \cong \Id$, take an arbitrary coalgebra $\alpha\tto \delta \lnadj F \rnadj F\alpha$  from $\EMCa \AAa F$ and consider the functor $\compss^0$ in Fig.~\ref{Fig:Eleft} instantiated to $\beta = \compss_0 \delta$. By extending the right-hand side of this instance of Fig.~\ref{Fig:Eleft} by the $\radj F$-image of the right-hand side of Fig.~\ref{Fig:Eright}, we get Fig.~\ref{Fig:delta}.
\begin{figure}[!ht]
\begin{center}
\[
\begin{tikzar}[row sep=3em]
\radj F \ladj F \radj F y \ar{dd}[swap] {\radj F\ladj F \radj F \compss_0 \delta}
\ar{rr}[name = U,below]{}{\radj F \varepsilon}  
\& \& 
\radj F y  \arrow[dashed]{dd}{\radj F \compss_0 \delta} \arrow[tail]{rr}{\radj F e} \&\& \radj F \ladj F x  \arrow{dd}{\radj F \ladj F \eta}
\arrow[tail,shift right = 1ex]{rr}[swap]{\radj F \ladj F \eta} \arrow[shift left = 1ex]{rr}{\radj F \ladj F d} \arrow[shift left = 1ex,bend left = 30,two heads]{ll}{\radj F r} \&\& \radj F \ladj F \radj F\ladj F x  \arrow{dd}{\radj F \ladj F \eta} \arrow[shift right = 1ex,bend right = 35,two heads]{ll}[swap]{\radj F \varepsilon}
\\ \\ 
\radj F \ladj F \radj F \ladj F \radj F y  \ar{rr}[name = D]{}[swap]{\radj F \varepsilon} \&\&
\radj F \ladj F \radj F y \arrow[tail]{rr}{\radj F\ladj F\radj F e} \&\& \radj F \ladj F\radj F\ladj F x  
\arrow[shift left = 1ex]{rr}{\radj F \ladj F\radj F\ladj F d} \arrow[tail,shift right = 1ex]{rr}[swap]{\radj F \ladj F\radj F\ladj  F \eta} \arrow[shift left = 1ex,bend left = 30,two heads]{ll}{\radj F\ladj F\radj F r} \&\& \radj F \ladj F \radj F\ladj F x  \arrow[shift right = 1ex,bend right = 35,two heads]{ll}[swap]{\radj F \ladj F \radj F \varepsilon} 
\arrow[thick,to path = (U) -- (D)\tikztonodes]{}{}[swap]{\displaystyle \compss^0 \compss_0 \delta}
\end{tikzar}
\]
\caption{Construction of $\compss^0 \compss_0 \delta$}
\label{Fig:delta}
\end{center}
\end{figure}
The claim is now that $x\mmono d \radj F \ladj F x$ equalizes the parallel pair $<\radj F \ladj F\eta, \radj F \ladj F d>$ in the first row. Since $y \mmono e \ladj F x$ was defined in Fig.~\ref{Fig:Eright} as a split equalizer of the pair $<\ladj F\eta,  \ladj F d>$, and all functors preserve split equalizers, it follows that $\radj F y \mmono{\radj F e} \radj F \ladj F x$ is also an equalizer of the same pair $<\radj F \ladj F\eta, \radj F \ladj F d>$. Hence the isomorphism $x\cong \radj F y$, which gives $\compss^0 \compss_0 \delta \cong \delta$. To prove the claim that $x\mmono d \radj F \ladj F x$ equalizes the first row, note that, just like the coalgebra $y\tto \beta \ladj F \radj F y$ in $\Emc \BBb F$ was determined up to isomorphism by the split equalizer in $\BBb$, shown in  Fig.~\ref{Fig:beta}, the coalgebra $\alpha \tto\delta \lnadj F\rnadj F \alpha$ in $\EMC \AAa F$ is determined up to isomorphism by the following split equalizer in $\Emm \AAa F$
\bea\label{eq:deltasplit}
\begin{tikzar}{}
\alpha  \arrow[tail]{rr}{\delta} \&\&\lnadj F \rnadj F \alpha  \arrow[tail,shift right = 1ex]{rr}[swap]{\lnadj F \eta} \arrow[shift left = 1ex]{rr}{\lnadj F \rnadj F \delta} \arrow[shift left = .5ex,bend left = 35,two heads]{ll}{\varepsilon} \&\& \lnadj F \rnadj F\lnadj F \rnadj F \alpha \arrow[shift right = .5ex,bend right = 35,two heads]{ll}[swap]{\varepsilon}\end{tikzar}
\eea
In $\AAa$, \eqref{eq:deltasplit}  unfolds to the lower squares of the Fig.~\ref{Fig:deltasplitA}. 
\begin{figure}[!ht]
\begin{center}
\[
\begin{tikzar}[row sep=3em]
x  \arrow[tail]{rr}{d} \arrow[tail,dashed]{dd}{d} \&\& \radj F \ladj F x  
\arrow[tail,shift right = 1ex]{rr}[swap]{\radj F \ladj F \eta} \arrow[shift left = 1ex]{rr}{\radj F \ladj F d} \arrow[shift left = 1ex,bend left = 35,two heads]{ll}[swap]{\alpha} \arrow[tail]{dd}{\radj F \ladj F\eta} \&\& \radj F \ladj F\radj F \ladj F x  
 \arrow[shift right = 1ex,bend right = 35,two heads]{ll}[swap]{\radj F \varepsilon} \arrow[tail]{dd}{\radj F \ladj F\eta}
 \\
\\
\radj F \ladj F x \arrow[two heads,dashed]{dd}{\alpha} \arrow[tail]{rr}{\radj F \ladj F d} 
\&\& \radj F\ladj F\radj F x  \arrow[two heads]{dd}{\radj F \varepsilon}
\arrow[shift left = 1ex]{rr}{\radj F\ladj F\radj F \ladj F d} \arrow[tail,shift right = 1ex]{rr}[swap]{\radj F\ladj F\radj  F \ladj F \eta} \arrow[shift left = 1ex,bend left = 35,two heads]{ll}[swap]{\radj F \ladj F \alpha} \&\& \radj F \ladj F\radj F \ladj F \radj F x  \arrow[shift right = 1ex,bend right = 35,two heads]{ll}[swap]{\ladj F \radj F \varepsilon} \arrow[two heads]{dd}{\radj F \varepsilon}
\\
\\
x  \arrow[tail]{rr}{d} \&\& \radj F \ladj F x  
\arrow[shift left = 1ex]{rr}{\radj F \ladj F d} \arrow[tail,shift right = 1ex]{rr}[swap]{\radj F \ladj F \eta} \arrow[shift left = 1ex,bend left = 35,two heads]{ll}{\alpha} \&\& \radj F \ladj F\radj F \ladj F x  
 \arrow[shift right = 1ex,bend right = 35,two heads]{ll}[swap]{\radj F \varepsilon}
\end{tikzar}
\]
\caption{The unfolding of \eqref{eq:deltasplit} and its splitting in $\AAa$}
\label{Fig:deltasplitA}
\end{center}
\end{figure}
Since the upper right-hand squares also commute (by the naturality of $\eta$), they also induce the factoring of the split equalizers in the upper left-hand square. But the upper right-hand squares in Fig.~\ref{Fig:deltasplitA} are identical to the right-hand squares in Fig.~\ref{Fig:delta}. The fact that both $\radj F y \mmono{\radj F e} \radj F \ladj F x$ and $x \mmono{\radj d} \radj F \ladj F x$ are split equalizers of the same pair yields the isomorphism $\radj F y \tto \iota x$ in $\AAa$, which turns out to be a coalgebra isomorphism $\compss^0 \compss_0 \delta \ttto \iota \sim \delta$ in $\EMC \AAa F$, as shown in Fig.~\ref{Fig:iota}.
\begin{figure}[!ht]
\begin{center}
\[
\begin{tikzar}{}
\radj F \ladj F x \ar{rrrr} {\alpha}  \ar{dddd}[swap]{\radj F\ladj F d} \&\&\&\& x\ar{dddd} {d} \\
\& \radj F \ladj F \radj F y \ar{ul}[swap]{\radj F \ladj F \iota} \ar{dd}{\radj F\ladj F \radj F \compss_0 \delta}
\ar{rr}
{\radj F \varepsilon}  
\& \& 
\radj F y  \ar{ur}{\iota} \arrow{dd}[swap]{\radj F \compss_0 \delta} 
\\ \\ 
\& \radj F \ladj F \radj F \ladj F \radj F y  \ar{dl}[swap]{\radj F\ladj F \radj F \ladj F \iota} \ar{rr}
{\radj F \varepsilon} \&\&
\radj F \ladj F \radj F y \ar{dr}{\radj F \ladj F \iota} \\
\radj F \ladj F \radj F \ladj F x \ar{rrrr}[swap]{\radj F\varepsilon} \& \& \& \& \radj F \ladj F x
\end{tikzar}
\]
\caption{$\compss^0 \compss_0 \cong \Id$}
\label{Fig:iota}
\end{center}
\end{figure}
Here the outer square is $\delta$, as in Fig.~\ref{Fig:Eright} on the left, whereas the inner square is $\compss^0 \compss_0 \delta$, as in Fig.~\ref{Fig:delta} on the left. The right-hand trapezoid commutes because the middle square in Fig.~\ref{Fig:delta} commutes, and can be chased down to Fig.~\ref{Fig:chase-der} using the fact that $\iota$ is defined by $\radj F e = d \circ \iota$. 
\begin{figure}[!ht]
\begin{center}
\[
\begin{tikzar}[row sep=7em,column sep=4em]
\radj F y \ar[bend left = 20]{rr} {\radj F e} \ar{r}[swap]{\iota} \ar{d}[swap]{\radj F \compss_0\delta} \& x \ar[tail]{r}[swap]{\delta} \&\radj F\ladj F x  \ar[tail]{d}{\radj F\ladj F \eta} \\
\radj F\ladj F\radj F\ladj F y \ar[bend right = 20]{rr}[swap]{\radj F\ladj F \radj F e}
\ar{r}{\radj F\ladj F \iota} \& \radj F\ladj F x \ar[tail]{r}{\radj F\ladj F\delta} 
\ar[equals,bend left = 25]{ur} 
\&\radj F\ladj F\radj F\ladj F x
\ar[bend right = 35,two heads]{l}[swap]{\radj F \varepsilon}
\end{tikzar}
\]
\caption{Commutativity of the right-hand trapezoid in Fig.~\ref{Fig:iota}}
\label{Fig:chase-der}
\end{center}
\end{figure}
The commutativity of the left-hand trapezoid in Fig.~\ref{Fig:iota} follows, because it is an $\radj F\ladj F$-image of the right-hand trapezoid. The bottom trapezoid commutes by the naturality of $\varepsilon$. The top trapezoid commutes because everything else commutes, and $d$ is a monic. The commutative diagram in Fig.~\ref{Fig:iota} thus displays the claimed isomorphism $\compss^0 \compss_0 \delta \tto \iota \delta$. This completes the proof that $\compss^0 \compss_0 \cong \Id$. Together with the proof that $\compss_0 \compss^0 \cong \Id$ in Fig.~\ref{Fig:beta}, it completes the proof that $\compss = \left(\compss^0\dashv \compss_0\right)$ provides the first of the equivalences
\[ \Emc \BBb F \simeq \EMC \AAa F \qquad\qquad\mbox{and}\qquad\qquad \EMM \BBb F \simeq\Emm \AAa F\]
claimed in Fig.~\ref{Fig:NucL-commut}. The second one is dual.
\end{proof}

\para{The adjunction $\NucL(F)=\left(\nadj F: \Emc \BBb F \to \Emm \AAa F\right)$ is monadic and comonadic.} By definition, an adjunction $\adj G\colon \DDd\to\CCc$ is monadic if the comparison functor $\compson^{1}\colon \DDd\to \Emm \CCc G$ in \eqref{eq:comparison} on the right is an equivalence. It is comonadic if the comparison functor $\compson_{0}\colon \CCc\to \Emc \DDd G$ in \eqref{eq:comparison} on the right is an equivalence. Thm.~\ref{thm:Theorem} thus implies that the nucleus $\NucL(F)$ of any adjunction $F$ is monadic and comonadic, as anticipated in Fig.~\ref{Fig:escher}.

\subsection{Consequences of the Nucleus Theorem}
The equivalence $\compss = \left(\adj \compss \right)$ proves the claim in the title of this section.

\begin{corollary}
$\NNucL\colon \Adj\to \Adj$ is an idempotent monad, with the natural isomorphism $\NNucL(F)\cong \NNucL\NNucL(F)$ realized by the unit in Fig.~\ref{Fig:two}.
\end{corollary}

\begin{figure}[!ht]
\begin{center}
\[\begin{tikzar}[row sep=0.75cm,column sep=3cm]
\Emca\BBb F 
\arrow[loop, out = 135, in = 45, looseness = 4,thin]{}[swap]{\Lft F} 
\arrow[bend right = 13]{dd}
\arrow{r}[pos=0.75]{\compss^{0}}[description]{\mbox{\LARGE$\simeq$}} 
\& 
\EMCa\AAa F 
\arrow[bend right = 13]{dd}  
\\
\hspace{.1em} \ar[Rightarrow,shorten=3em,shift right=2]{r}{\eta} \&\hspace{.1em}
\\
\Emma \AAa F  \arrow[loop, out = -45, in=-135, looseness = 6,thin]{}[swap]{\Rgt F} 
\arrow{r}[description]{\mbox{\LARGE$\simeq$}}[swap,pos=0.75]{\compss_{1}}  
\arrow[bend right = 13]{uu}
\arrow[phantom]{uu}[description]{\mbox{\small $\NNucL(F)$}}  
\& 
\EMMa\BBb F  
\arrow[bend right = 13]{uu}  
\arrow[phantom]{uu}[description]{\mbox{\small $\NNucL\NNucL(F)$}}   
\end{tikzar}\]
\caption{The unit $\eta\colon \NNucL(F)\to \NNucL\NNucL(F)$}
\label{Fig:two}
\end{center}
\end{figure}

\begin{corollary}
If $F=\left(\adj F\colon \BBb\to \AAa\right)$ and $G=\left(\adj G\colon \DDd\to \AAa\right)$ induce the same monad, then they have the same nucleus: 
\bea
\radj F\ladj F = \lft T = \radj G\ladj G & \implies & 
\NNucL(F)\cong\NNucL(G)
\eea
In particular, $\radj F\ladj F= \radj G\ladj G$ on $\AAa$  implies $\Emca\BBb F\simeq \Emca \DDd G$.
\end{corollary}

\begin{proof}
The assumption $\radj F \ladj F =\lft T = \radj G\ladj G$ means that $F$ and $G$ have the same final resolution $\EMnd(F) = \left(\Klm U T\vdash \Emm U T\colon \Emm \AAa T\to \AAa\right) =\EMnd(G)$. It follows that they also induce the same descent comonad $\Rgt T =  \Emm U T \Klm U T$ on $\Emm \AAa T$. Lemma~\ref{lemma:nuc-decom} now gives
\[ \lnadj F \rnadj F\ =\ \Rgt T \ =\ \lnadj G \rnadj G\]
whereas Thm.~\ref{thm:Theorem} proves
\[\Emca \BBb F \ \simeq \EMCa \AAa T\ \simeq\ \Emca \DDd G\]
Since the comparison functors commute with the nucleus functors, we have proved $\NNucL(F)\cong\NNucL(G)$. 
\end{proof}

\para{Remark.} Notice that the proof and its dual in essence prove
\bea\label{eq:nucleq-remark}
\EEMnd(F)\cong\EEMnd(G)\ \vee\ \EECmn(F)\cong \EECmn(G)& \implies & \NNucL(F) \cong\NNucL(G)
\eea
Stated in this form, the claim is an obvious   consequence of Thm.~\ref{thm:Theorem}. However, unpacking the implications from the two disjuncts provides useful statements. The first disjunct says that any equivalences $\Kar\AAa \simeq \Kar\CCc$ and $\Emma\AAa F \simeq \Emma \CCc G$, coherent with the forgetful and the free functors, imply $\Emca \BBb F \simeq \Emca \DDd G$. Dually, the second disjunct says that coherent equivalences $\Kar\BBb \simeq \Kar\DDd$ and $\Emca\BBb F \simeq \Emca \DDd G$ imply $\Emma \AAa F \simeq \Emma \CCc G$.

\section{Simple nucleus}\label{Sec:simple}

The equivalences $\Emma \AAa F \simeq \EMMa \BBb F$ and $\Emca \BBb F \simeq \EMCa \AAa F$ from Thm.~\ref{thm:Theorem} present $\lft F$-algebras as $\Lft F$-algebras over $\rgt F$-coalgebras, and $\rgt F$-coalgebras as $\Rgt F$-coalgebras over $\lft F$-algebras. But these presentations are redundant, and discharging the redundancies yields interesting and useful alternatives to the familiar Eilenberg-Moore views of algebras and coalgebras --- in terms of idempotents. In a subtly different form, they were used in \cite{PavlovicD:LICS17} for a particular security application. Here we prepare them for mathematical applications.  
\begin{proposition}\label{Prop:three}
Given an adjunction $F=\left(\adj F : \Kar \BBb\to \Kar \AAa \right)$, consider  the categories
\bea\label{eq:Ec}
|\Ec \AAa F| & = & \coprod_{x\in |\Kar \AAa|}\  \left\{\alpha_x \in \Kar \BBb(\ladj F x, \ladj F x) \ \ \big|\ \ \begin{tikzar}[row sep=1em,column sep=.25em]
\ladj F x \ar{dd}[description]{\alpha_x} \ar{ddrr}[description]{\alpha_x} 
\&
\radj F \ladj F x \ar{ddrr}[description]{\radj F \alpha_x} \ar[two heads]{rr} \&\& x\ar[tail]{dd}[description]{\tilde \alpha_x}
\\ \\
\ladj F x \ar{rr}[description]{\alpha_x} \&\& \ladj Fx \& \radj F\ladj F x
 \end{tikzar}\ \ \right\}  \label{eq:monCoalg}
\\[3ex]
\Ec \AAa F (\alpha_x, \gamma_z) & = &\hspace{2em}  \left\{f \in \Kar \AAa(x,z)\ \Big|\ \ 
\begin{tikzar}[row sep=1.8em,column sep=1.8em]
\ladj F x \ar{r}{\ladj F f} \ar{d}[description]{\alpha_x} \& \ladj F z \ar{d}[description]{\gamma_z}  \\ 
\ladj F x \ar{r}{\ladj F f}\& \ladj F z 
\end{tikzar}\ 
\right\} \notag
\eea

\bea\label{eq:Em}
|\Em \BBb F| & = & \coprod_{u\in |\Kar \BBb|}\  \left\{\beta^u \in \Kar \AAa(\radj F u, \radj F u) \ \ \big|\ \ \begin{tikzar}[row sep=1em,column sep=.25em]
\radj F u \ar{dd}[description]{\beta^u} \ar{ddrr}[description]{\beta^u} 
\&
\ladj F \radj F u \ar{ddrr}[description]{\ladj F \beta^u} \ar[two heads]{rr}[description]{\tilde \beta^u} \&\& u\ar[tail]{dd} 
\\ \\
\radj F u \ar{rr}[description]{\beta^u} \&\& \radj Fu \& \ladj F\radj F u
 \end{tikzar}\ \ \right\}  \label{eq:algComon}
\\[3ex]
\Em \BBb F (\beta^u, \delta^w) & = &\hspace{2em}  \left\{g \in \Kar \BBb(u,w)\ \Big|\ \ 
\begin{tikzar}[row sep=1.8em,column sep=1.8em]
\radj F u \ar{r}{\radj F g} \ar{d}[description]{\beta^u} \& \radj F w \ar{d}[description]{\delta^w}  \\ 
\radj F u \ar{r}{\radj F g}\& \radj F w 
\end{tikzar}\ 
\right\} \notag
\eea
where $x\tto{\tilde \alpha_x}\radj F\ladj F x$ is the transpose of $\ladj F x\tto{\alpha_x}\ladj Fx$, and $\ladj F\radj F u \tto{\tilde \beta^u} u$ is the transpose of $\radj F u \tto\beta \radj F u$. The adjunction $\oadj F: \Em \BBb F \to \Ec \AAa F$ defined in Fig.~\ref{Fig:simplenuc}
\begin{figure}[!ht]
\begin{center}
\[\begin{tikzar}[row sep=1.5mm,column sep=1.5mm]
\AAa \arrow[phantom]{dddddddd}[description]{\dashv}  
\arrow[loop, out = 135, in = 45, looseness = 4,thin]{}[swap]{\lft F} 
\arrow[bend right = 13]{dddddddd}[swap]{\ladj F} 
\arrow[thin,dashed]{rrrrrrrr}{\compstwo^{0}} 
\& \&\&\& \& \&\&\&
\Ec\AAa F 
\arrow[phantom]{dddddddd}[description]{{\dashv}}   
\arrow[bend right = 13]{dddddddd}[swap]{\loadj{F}}  
\\ \&\hspace{3em}\&\&\&\&\& \& <x, \alpha_x> \ar[mapsto,bend right = 13]{dddddd} 
\&\& 
<\radj F u, \ladj F \beta^u>
\\ \hspace{3em}\\ \hspace{3em} \\ \hspace{3em} \\ \\ \\ 
\&\&\&\&\&\&\& <\ladj F x, \radj F \alpha_x  >\&\& <u,\beta^u> \ar[mapsto,bend right = 13]{uuuuuu}
\\ 
\BBb  
\arrow[loop, out = -45, in=-135, looseness = 6,thin]{}[swap]{\rgt F} 
\arrow[thin,dashed]{rrrrrrrr}[swap]{\compstwo_{1}}  
\arrow[bend right = 13]{uuuuuuuu}[swap]{\radj F} 
\& \&\&\&\& \&\&\&
\Em \BBb F 
\arrow[bend right = 13]{uuuuuuuu}[swap]{\roadj{F}}  
\end{tikzar}\]
\caption{The simple nucleus $\oadj F$ of $\adj F$}
\label{Fig:simplenuc}
\end{center}
\end{figure}
with the comparison functors 
\begin{align*}
\compstwo^0\ :\ \AAa & \tto{\hspace{3em}} \Ec\AAa F
& \compstwo_1:\BBb  &  \tto{\hspace{3em}} \Em \BBb F\\[1ex]
x & \longmapsto\ \left<\radj F \ladj F x,\ \  \begin{tikzar}[row sep = 5mm,cramped]
\scriptstyle \ladj F \radj F \ladj F x\ar{d}[swap,pos=0.4]{\varepsilon\ladj F }\\
\scriptstyle \ladj F x \ar{d}[swap,pos=0.4]{\ladj F\eta}\\
\scriptstyle\ladj F \radj F \ladj F x
\end{tikzar}\ \ \  
\right> 
 & u & \longmapsto\ \left<\ladj F\radj F u, \ \ 
 \begin{tikzar}[row sep = 5mm,cramped]
\scriptstyle \radj F \ladj F \radj F u\ar{d}[swap,pos=0.4]{\radj F \varepsilon}\\
\scriptstyle \radj F u \ar{d}[swap,pos=0.4]{\eta\radj F}\\
\scriptstyle\radj F \ladj F \radj F u
\end{tikzar}\ \ \  
 \right>
\\[-2ex]
\end{align*}
is equivalent to the nucleus, i.e.
\bear
\NucL\left(\adj F\right) & \simeq & \left(\oadj F \right)
\eear
\end{proposition}
 
%

\begin{lemma}\label{Lemma:splitting-two} If $\left(\radj F\ladj Fx \tto{\radj F \alpha_x} \radj F \ladj Fx\right) = \left(\radj F\ladj Fx \eepi{\alpha^x} x\mmono{\tilde \alpha_x} \radj F \ladj Fx\right)$, where $\tilde \alpha_{x}$ is the adjoint transpose of $\alpha_{x}$, then $\alpha^{x}\circ \tilde \alpha_{x}= \id$ implies $\alpha^x\circ \eta_x = \id$.
\end{lemma}

\bpr The assumption that $\tilde \alpha_{x}$ is the adjoint transpose implies $\tilde \alpha_x\circ \alpha^x = \radj F \alpha_x \circ \eta_x \circ \alpha^x = \tilde \alpha_x \circ \alpha^x \circ \eta_x \circ \alpha^x$. Pre- and postomposing with $\alpha^{x}\circ \tilde \alpha_{x}= \id$ gives $\alpha^x\circ \eta_x = \alpha^{x}\circ(\tilde \alpha_x \circ \alpha^x \circ \eta_x \circ  \alpha^x) \circ \tilde \alpha_{x} =\alpha^{x}\circ \tilde \alpha_{x}\circ \alpha^{x}\circ \tilde \alpha_{x}= \id$.
\epr

\bprf{\ \bf of Prop.~\ref{Prop:three}}
Splitting the idempotent $\radj F\ladj Fx \tto{\radj F\alpha_x} \radj F\ladj F x$ into $\radj F\ladj Fx \eepi{\alpha^x} x\mmono{\tilde \alpha_x} \radj F\ladj F x$ as above gives
\beq\label{eq:algebra} \alpha^x \circ \tilde \alpha_x = \id_x \qquad \qquad \mbox{ and } \qquad \qquad \alpha^x \circ \eta_x = \id_x\eeq
Analogous reasoning proves  
\beq\label{eq:algebra-morph}
\alpha^x\circ \radj F\ladj F \alpha^x  =  \alpha^x\circ \radj F\varepsilon_{\radj Fx}
\qquad \qquad \mbox{ and } \qquad \qquad \tilde \alpha_x \circ \alpha^x  = \radj F \varepsilon_{\radj F x} \circ \radj F\ladj F \tilde \alpha_x
\eeq
Together, (\ref{eq:algebra}--\ref{eq:algebra-morph}) say that $\radj F \ladj F x\tto{\alpha^x} x$ is a algebra in $\Emma \AAa F$ and that $\tilde \alpha_x\in \Emma \AAa F(\alpha^x, \mu_x)$ is a coalgebra over $\alpha^x$ in $\EMCa \AAa F$. Hence a functor $\Ec \AAa F\to \EMCa \AAa F$. Similar construction yields a similar functor $\Em \BBb F\to \EMMa \BBb F$. Hence the equivalences
\beq\label{eq:simple-descmnd}
\Ec \AAa F \simeq \EMCa \AAa F \qquad \qquad \qquad \qquad \Em \BBb F \simeq \EMMa \BBb F
\eeq
The equivalences
\beq
\Ec \AAa F \simeq \Emca \BBb F \qquad \qquad \qquad \qquad \Em \BBb F \simeq \Emma\AAa F 
\eeq
follow from Thm~\ref{thm:Theorem} and have also been verified directly in \cite{PavlovicD:LICS17}. Every object $<x,\ladj Fx\tto{\alpha_x} \ladj F x>$ of  $\Ec\AAa F$ is thus isomorphic to one in the form $<\radj F y, \ladj F\radj F \tto\varepsilon y \tto\beta \ladj F\radj F y>$, where $\beta$ is a coalgebra in $\Emca \BBb F$. It follows that the inner and the outer squares in the following diagram commute
\beq\label{eq:outsquare}
\begin{tikzar}[column sep = scriptsize, row sep = normal]
\ladj F x \ar{rrrrrr}{\ladj F \eta} 
\ar[leftrightarrow]{dr}[description]{\ladj F\iota} 
\ar{dddddd}[description]{\alpha_x}
\&\&\&\&\&\& \ladj F\radj F\ladj F x \ar{dddddd}[description]{\alpha_x} \ar[leftrightarrow]{dl}[description]{\ladj F\radj F\ladj F \iota} 
\\
\& \ladj F\radj F y 
\ar[two heads]{dd}[description]{\varepsilon} 
\ar{rrrr}{\ladj F\eta\radj F} \&\&\&\& \ladj F\radj F\ladj F \radj F y
\ar[two heads]{dd}[description]{\ladj F\radj F\varepsilon}
\\
\\ 
\&  y \ar[tail]{rrrr}[description]{\beta} 
\ar[tail]{dd}[description]{\beta} \&\&\&\& \ladj F\radj F y\ar{dd}[description]{\ladj F\radj F\beta}
\\ 
\\
\& \ladj F\radj F y 
\ar{rrrr}{\ladj F\eta \radj F} \&\&\&\& \ladj F\radj F\ladj F \radj F y 
\\
\ladj F x \ar{rrrrrr}{\ladj F \eta} \ar[leftrightarrow]{ur}[description]{\ladj F\iota} 
\&\&\&\&\&\& \ladj F\radj F\ladj F x 
\ar[leftrightarrow]{ul}[description]{\ladj F\radj F\ladj F\iota} 
\end{tikzar}\eeq
and that they are linked by  an isomorphism $\begin{tikzar}{} x \ar[leftrightarrow]{r}[description]{\iota}\&\radj Fy\end{tikzar}$ in $\AAa$. Transferring the nuclear adjunction \mbox{$\nadj F :\Emma \AAa F \to \Emca \BBb F$} along the equivalences yields the nuclear adjunction $\oadj F :\Em \BBb F \to \Ec \AAa F$, with the natural correspondence
\bear
\Em \BBb F(\loadj F \alpha_x, \beta^u) & \cong & \Ec \AAa F(\alpha_x, \roadj F\, \beta^u)\\
\left(\ladj F x \tto f u\right) & \mapsto & \tilde f = \left(x\tto{\ \ \eta\ \ } \radj F\ladj F x\tto{\radj F f} \radj F u\right)
\eear 
The transpositions along $\adj F:\BBb\to \AAa$ lift to $\oadj F :\Em \BBb F \to \Ec \AAa F$ because each of the following squares commutes if and only if the other one does: 
\bea\label{eq:tinier}
\begin{tikzar}[row sep=2em,column sep=1em]
\radj F\ladj F x \ar{rr} {\radj F f} \ar{dd} [description]{\radj F \alpha_x} \& \& \radj F u \ar{dd}[description]{\beta^u}\\ \& \hspace{1em}  \\
\radj F\ladj F x \ar{rr}[swap] {\radj F f} \&\& \radj Fu
\end{tikzar}
\hspace{1.5em}& \iff &\hspace{1.5em}
\begin{tikzar}[row sep=2em,column sep=1em]
\ladj F x \ar{rr} {\ladj F\left(\tilde f\right)} \ar{dd}[description]{\alpha_x}  \&\&   \ladj F \radj F u \ar{dd}[description]{\ladj F \beta^u}
\\ 
\& \hspace{1em} \\
\ladj F x \ar{rr}[swap]{\ladj F\left(\tilde f\right)}  \&\&   \ladj F \radj F u
\end{tikzar}
\eea
To see this, suppose that the left-hand side square commutes. Take the $\ladj F$-image of the right-hand side square and precompose it with the outer square from \eqref{eq:outsquare}, as in the following diagram. 
\beq\label{eq:tiny}
\begin{tikzar}[row sep=5em,column sep=normal]
\ladj F x \ar[bend left = 25]{rrrr}{\ladj F\left(\tilde f\right)} 
\ar{d}[description]{\alpha_x} \ar[tail]{rr}{\ladj F \eta} 
\&\& \ladj F\radj F \ladj F x \ar{d}[description]{\ladj F  \radj F \alpha_x} 
\ar{rr}[swap]{\ladj F\radj F f} 
\&\& \ladj F \radj F u\ar{d}[description]{\ladj F \beta^u} 
\\
\ladj F x \ar[bend right = 25]{rrrr}[swap]{\ladj F\left(\tilde f\right)} 
 \ar[tail]{rr}[swap]{\ladj F \eta} \&\&  \ladj F\radj F \ladj F x \ar{rr}{\ladj F\radj F f} 
\&\& \ladj F \radj F u
\end{tikzar}
\eeq 
The two outer paths around this diagram are the paths around the right-hand square in \eqref{eq:tinier}. The right-to-left implication is analogous. \epr

\section{New frameworks for old big pictures}\label{Sec:HT}

\subsection{Descent to nucleus}\label{Sec:desc}
In an (over)simplified form, the idea of descent is to descend to a complex object by covering it with simpler objects. Sheaves can be viewed as a special case. A sheaf is a set $\varsigma$, continuously varying over a topological space $X$. The continuous variation means that it is indexed over the open neighborhoods of $X$ and the elements of $\varsigma B$, for an open neighborhood $B\subseteq X$, can be described by descending from the simpler elements of $\varsigma A_i$ down a covering $B = \bigcup_{i\in I} A_i$. The sets varying over $X$ are formalized as functors $\OOO X^\op\to \Set$, the presheaves,  where $\OOO X$ is the lattice of open neighborhoods of $X$. A presheaf $\varsigma\colon \OOO X ^\op\to \Set$ is a sheaf if every covering $B = \bigcup_{i\in I} A_i$ is mapped to a decomposition 
\bea\label{eq:sheaf}
\varsigma B & \cong & \left\{ <si>_{i\in I}\in \prod_{i\in I} \varsigma A_i\ |\ \forall ij\in I.\ si\restr_{A_{i}\cap A_{j}} =  sj\restr_{A_{i}\cap A_{j}}\right\}
\eea
where $sU\restr_{V}$ is the restriction of $sU\in \varsigma U$, mapped by $\varsigma$-image $\varsigma U \tto\restriction \varsigma V$ of the morphism $U\supseteq V$ in $\OOO X ^{\op}$. This \emph{sheaf condition}\/ formalizes the intuition that the set $\varsigma B$ descends from the family of sets $<\varsigma A_i >_{i\in I}$ just like $B$ descends from its cover $\{A_{i}\}_{i\in I}$: that is, just like $B$ is obtained by gluing $A_i$s and $A_j$s along the intersections $A_i\cap A_j$, the set $\varsigma B$ is obtained by gluing $\varsigma A_{i}$s and $\varsigma A_{j}$s along the intersections $\varsigma (A_{i}\cap A_{j})$. 

The (over)simplified idea of Grothendieck descent \cite{GrothendieckA:fibrations59,GrothendieckA:SGA1} is to lift the concept of gluing from sheaves as continuously variable sets $\OOO X^{\op}\to \Set$ to continuously variable categories $\SSs^{\op}\to \Cat$, where the continuity with respect to $\SSs$ is expressed in terms of coverings of objects by \emph{sieves}\/ in the slices \cite{GiraudJ:descente}. An indexed category $\SSs^{\op}\to \Cat$ satisfying the suitable gluing condition, usually  comprehended as a fibration, is called a \emph{stack}. The concept originated in algebraic geometry  \cite{GiraudJ:cohna,VistoliA:descent}, expanded into homotopy theory \cite{GrothendieckA:pursuing,JardineJF:homotopy-descent} and many areas of mathematics in-between and beyond \cite{Janelidze-Tholen:facets}.  The breadth of applications makes the presentations wildly divergent and implementation-dependent, and it is not easy to see the forest for the trees. The first author proposed a logical interpretation of descent as an interpolation condition \cite{PavlovicD:Como}. The nucleus seems to suggest a further simplification. 

For intuition, consider an indexed category $\Pow \colon \SSs^{\op}\to \Cat$ mapping types (or sets) from a universe $\SSs$ to categories $\Pow X$ whose objects, written $\alpha(x), \gamma(u),\ldots$, can be thought of as observable properties or predicates over $X$.\footnote{If $\SSs =\Esp$ is the category of topological spaces and $\Pow X =\OOO X$ is the lattice of open neighborhoods of a space $X$, then a predicate is a continuous map $\alpha\colon X\to \TTwo$ into the Sierpi\'nsky space $2=\{0,1\}$ with $\{1\}$ open. The comprehension of the indexed category $\OOO \colon \Esp^{\op}\to \Cat$ is then the fibration $\Esp/\TTwo \tto{\Dom} \Esp$. Here we restrict to the discrete case since we just need to build logical intuition.} The morphisms $p \in \Pow X(\alpha, \gamma)$ are then thought of as proofs of $\alpha(x)\Rightarrow \gamma(x)$, along the lines of the Brouwer-Heyting-Kolmogorov interpretation of constructivist logic, echoed in the propositions-as-types interpretation of type theory. See \cite{Coquand-Huet:constructions,JacobsB:book,Martin-LoefP:inttt,PavlovicD:thesis} for formal developments of such interpretations. The starting point of the corresponding categorical developments was Lawvere's definition of \emph{hyperdoctrines} \cite{LawvereFW:dialectica,LawvereFW:compreh}. We call $\Pow$ a hyperdoctrine even if the closed structure, central in the original treatments, presently plays no role. A hyperdoctrine $\Pow$ is thus a functor assigning to every function $f$ in the universe $\SSs$ an adjoint pair of functors $\badj f$, where $\rbadj f$ is interpreted as the substitution operation and $\lbadj f$ provides the corresponding existential quantifiers. The notation in Fig.~\ref{Fig:Pow} was formalized and used in \cite{PavlovicD:thesis,PavlovicD:constructions,PavlovicD:mapsII}.
\begin{figure}[!ht]
\begin{center}
$
\prooftree
f\colon A\to B
\justifies
\badj f\colon \Pow B\to \Pow A
\endprooftree
$
\qquad\qquad
$
\rbadj f\Big(\beta(y)\Big)\  = \ \beta\Big(f(x)\Big)$\qquad\qquad
$\lbadj f\Big(\alpha(x)\Big) \ = \ \exists x:A.\ f(x) \stackrel B =y\wedge \alpha(x)$
\caption{The arrow part of the functor $\Pow\colon \CCC^{\op}\to \Cat$}
\label{Fig:Pow}
\end{center}
\end{figure}
In a high-level logical view, the idea of descent is that the predicates over $B$ can be glued from predicates over $A$ along a covering morphism $f\colon A\to B$. The requirement in \eqref{eq:sheaf}, that a sheaf $\varsigma$ must map every covering $\bigcup_{i\in I} A_i = B$ to a corresponding decomposition of the elements of $\varsigma B$ is generalized to the requirement that $\Pow$ must map every cover by morphisms $\{f_{i}\colon A_{i}\to B\}_{i\in I}$, with the intersections annotated as in Fig.~\ref{Fig:intersect}, to a decomposition of the predicates in the form
\begin{multline}\label{eq:stack}
\Pow B \ \  \simeq \ \  \Bigg\{ <d_{ij}>_{i,j\in I}\in \prod_{i,j\in I} 
\coprod_{%
\substack{x_{i}\in \Pow A_{i} \\
x_{j}\in \Pow A_{j}}} 
\Pow E_{ij}\left(\delta_{ij}^{\ast}x_{i}, \rho_{ij}^{\ast} x_{j}\right)\ \Big|\\ 
\forall i\in I.\ \ \ \eta^{\ast}_{i} d_{ii}=\id_{x_{i}}\ \wedge\  \forall ijk\in I.\ \ \ p_{ij}^{\ast}d_{ij}\circ p_{jk}^{\ast}d_{jk}=p_{ik}^{\ast}d_{ik}\ \Bigg\}
\end{multline}
\begin{figure}[!ht]
\begin{center}
$\begin{tikzar}[column sep = 4em,row sep = 6ex]
\& A_{i}\ar[dashed]{d}[description]{\eta_{i}}
\ar[bend right]{ddl}[description]{\id} \ar[bend left]{ddr}[description]{\id} \\
\& E_{ii} \ar{dl}[description]{\delta_{ii}} \ar{dr}[description]{\rho_{ii}} \ar[phantom]{dd}[description,rotate=135,pos=0.1]{\ulcorner}
\\
A_{i}\ar{dr}[description]{f_{i}} \& \& A_{i}\ar{dl}[description]{f_{i}}\\
\& B
\end{tikzar}$
\qquad\qquad\qquad\qquad
$\begin{tikzar}[column sep = 4em,row sep = 6ex]
\& \widehat E_{ijk} \ar{dl}[description]{p_{ij}} \ar{d}[description]{p_{ik}} \ar{dr}[description]{p_{jk}} \\
E_{ij}\ar{d}[swap,pos=0.65]{\delta_{ij}} 
\ar{dr}[pos=0.75]{\rho_{ij}}\&
E_{ik} \&
E_{jk}\ar{dl}[swap,pos=0.75]{\delta_{jk}} \ar{d}[pos=0.65]{\rho_{jk}} \\
A_{i} \ar{dr}[description]{f_{i}} \ar[leftarrow,crossing over]{ur}[pos=0.75]{\delta_{ik}} \&
A_{j}\ar{d}[description]{f_{j}} \&
A_{k}\ar{dl}[description]{f_{k}} \ar[leftarrow,crossing over]{ul}[swap,pos=0.75]{\rho_{ik}}\\
\& B
\end{tikzar}$
\caption{The intersections of the elements of cover $\{f_{i}\colon A_{i}\to B\}_{i\in I}$ are their pullbacks}
\label{Fig:intersect}
\end{center}
\end{figure}
A hyperdoctrine $\Pow$ satisfying \eqref{eq:stack} is called a \emph{stack} \cite[Ch.~II]{GiraudJ:cohna}\footnote{The book is in French, so stacks are called \emph{champs}. The stack condition goes back to a much earlier seminar on Grothendieck's reconstruction of Teichm\"uller spaces. Although Giraud did not work with indexed categories but with their fibered comprehensions, and hyperdoctrines were not yet defined, the idea of stack condition was there.}. Note that \eqref{eq:stack} only specifies the object part of the \emph{stack condition}. To lift the sheaf condition from sets to categories, the stack condition also has a morphism part, which is not hard to reconstruct. Here, though, we are interested in a special case of the stack condition, and in special stack morphisms. When $\SSs$ has stable coproducts $A = \coprod_{i\in I} A_{i}$, or when we freely adjoin them, then $B$ can be covered by a single morphism $f\colon A\to B$, where $f=[f_{i}]_{i\in I}$, and the stack requirement in \eqref{eq:stack} boils down to the \emph{descent condition}   
\bea\label{eq:descent}
\Pow B &  \simeq & \coprod_{x\in \Pow A}\Bigg\{ d \in \Pow E\left(\delta^{\ast}x, \rho^{\ast}x\right)\ \Big|\  
\eta^{\ast} d=\id_{x}\ \wedge\ p_{0}^{\ast}d\circ p_{1}^{\ast}d =p_{2}^{\ast}d\ \Bigg\}
\eea
The notations are obtained from Fig.~\ref{Fig:intersect} again, renaming $p_{ij}\mapsto p_{0}$, $p_{jk}\mapsto p_{0}$, and $p_{ik}\mapsto p_{2}$ and erasing all other subscripts. Note that the right-hand side still only specifies  the objects of the category that the descent condition requires to be equivalent to $\Pow B$. These objects are  the \emph{descent data}. The sheaf condition in  \eqref{eq:sheaf} evolved to the stack condition in \eqref{eq:stack}, and to the descent condition in \eqref{eq:descent}. The category of descent data along $f$ is usually denoted by $\Desc_{f}$. Its objects, described in \eqref{eq:descent} on the right, are conveniently written in the form $x = \left(\delta^{\ast}x\tto{d_{x}} \rho^{\ast}x\right)$. The morphisms between descent data are:
\bea\label{eq:Desc}
\Desc_{f}(x,y) & = & \left\{\varphi \in \Pow A(x,y)\ \ \ \Big|\ \ 
\begin{tikzar}{}
\delta^{\ast}x \ar{d}[description]{d_{x}} \ar{r}[description]{\delta^{\ast}\varphi} \& \delta^{\ast}y \ar{d}[description]{d_{y}}\\
 \rho^{\ast}x \ar{r}[description]{\rho^{\ast}\varphi} \& \rho^{\ast}y 
\end{tikzar}
\right\}
\eea
Note the similarity with the simple nucleus morphisms in (\ref{eq:Ec}--\ref{eq:Em}). This is not an accident. 
\begin{figure}[!ht]
\begin{center}
\prooftree
\prooftree
\delta^{\ast}x\tto{\ \ d\ \ } \rho^{\ast}x
\justifies
\rho_{!}\delta^{\ast}x \tto{\ \ d'\ \ }x
\endprooftree
\justifies
f^{\ast}f_{!} x 
\tto{\ \ d''\ \ }x
\endprooftree
\caption{The descent data down the kernel $<\delta,\rho>$ of $f$ correspond to $f^\ast f_!$-algebras}
\label{Fig:BBR}
\end{center}
\end{figure}
Fig.~\ref{Fig:BBR} shows the  transpositions of the $f$-descent data using the adjunction $\rho_{!}\dashv \rho^{\ast}\colon \Pow A\to \Pow E$ and the \emph{Beck-Chevalley}\/ isomorphism\footnote{This isomorphism was attributed to unpublished work of Jon Beck by Lawvere \cite{LawvereFW:compreh} and to unpublished work of Claude Chevalley by B\'enabou and Roubaud \cite{Benabou-Roubaud}.}
\bea\label{eq:BC}
\rho_{!}\delta^{\ast}x &\cong &
f^{\ast}f_{!}
\eea
As explained in \cite{PavlovicD:Como}, the Beck-Chevalley isomorphism assures that the predicates of a hyperdoctrine $\Pow$ recognize the pullbacks in $\SSs$. In the logical notation of Fig.~\ref{Fig:Pow}, \eqref{eq:BC}  becomes
\bea
\exists e : E.\ \delta(e)=x\ \wedge\  \rho(e) = x' & \iff & f(x) = f(x')
\eea
If the subobject $E\tto{<\delta,\rho>} A\times A$ is viewed as a binary relation, and $\exists e : E.\ \delta(e)=x\ \wedge\  \rho(e) = x'$ is abbreviated to $xEx'$, then \eqref{eq:BC} becomes $xEx'\iff f(x)=f(x')$. $E$ is thus the kernel of the function $f$.  The crucial observation of \cite{Benabou-Roubaud} was that the transpositions in Fig.~\ref{Fig:BBR} map the conditions in \eqref{eq:descent},  which make $\delta^{\ast}x\tto{\ \ d\ \ } \rho^{\ast}x$ into descent data, precisely correspond to the conditions in \eqref{Fig:EM-EC}, which make $\rbadj f\lbadj f x\tto{d''} x$ into an $\lft f$-algebra for the monad $\lft f=\rbadj f\lbadj f$. The category $\Desc_{f}$ of descent data is equivalent to the category $\Emm {\Pow A} f$ of algebras. The descent condition \eqref{eq:descent} thus implies that $\Pow B$ must be equivalent to  $\Emm {\Pow A} f$. But  Prop.~\ref{Prop:three} says that $\Emm {\Pow A} f$ is equivalent to  ${\Pow B}^{\Lft f}$. In summary, a morphism $f\colon A\to B$ in $\SSs$ satisfies the descent requirement if and only if
\beq\label{eq:desb}
\Pow B\ \simeq \ \Desc_{f}\ \simeq \ \Emm {\Pow A} f \ \ \simeq\ \ {\Pow B}^{\Lft f}
\eeq 
where $\Pow A$ and $\Pow B$ are assumed to be absolutely complete.\footnote{Without this assumption, $\Pow B$ and ${\Pow B}^{\Lft f}$ must be weakly (Morita)  equivalent.} 

\para{Simple descent monads.} In words, \eqref{eq:desb} says that $\Pow A$ covers $\Pow B$ along a descent morphism $f\colon A\to B$ precisely when the objects of $\Pow B$ are retracts of the objects of $\Pow A$ along $f^{\ast}$, as specified in ${\Pow B}^{\Lft f}$. Lifting the $\rbadj f$-retractions to $\lft f$-algebras in $\Pow A$, and then transposing them to $\Pow E$ along the Beck-Chevalley and $\left(\badj \rho\right)$-adjunction yields the $f$-descent data in $\Desc_{f}$.  On the other hand, the equivalence \eqref{eq:simple-descmnd} of the retractions in ${\Pow B}^{\Lft f}$ and the $\Lft F$-algebras in ${\Pow B}^{\rgt f}$ explain the utility of descent monads and comonads  \cite{BalmerP:annalen,Caenepeel:galois,HessK:general,MesablishviliB:comonDesc}.

\subsection{Kan's adjunction}
In the final section, we consider the initial example of adjunctions. The idea of adjunction goes back at least to \'Evariste Galois, or, depending on how you think of it, much further back to Heraclitus \cite{LambekJ:Heraclitus}, and still further into the roots of logic \cite{LawvereFW:dialectica}. Yet the categorical definition of an adjunction as two categories, two functors, and two natural transformations goes back to the late 1950s, to Daniel Kan's work in homotopy theory  \cite{KanD:adj}. Kan defined the Kan extensions to capture this particular adjunction, like Eilenberg and MacLane defined categories and functors to define certain natural transformations.

\para{Simplices and the simplex category.}
Eilenberg simplified computing homology groups by decomposing topological spaces into simplices \cite{EilenbergS:singular}. An $m$-simplex is the set
\bea\label{eq:simpl}
\Delta_{[m]} & = & \left\{\vec x \in [0,1]^{m+1}\ \big|\ \sum_{i=0}^m x_i = 1\right\}
\eea
with the product topology induced by the open intervals on $[0,1]$. The relevant structure of a topological space $X$ is captured by the families of continuous maps $\Delta_{[m]} \to X$, for all $m\in \NNn$. Some of these maps do not \emph{embed}\/ some simplices into the space, but contain degeneracies or singularities; yet considering all such maps assures that any simplices that can be embedded into $X$ will be embedded by some of the maps. Since the simplicial structure is captured by each $\Delta_{[m]}$'s projections onto all $\Delta_{[\ell]}$s for $\ell \lt m$, and by $\Delta_{[m]}$'s  embeddings into all $\Delta_{[n]}$s for $n\gt m$, a simplicial structure of a space corresponds to a functor of the form \mbox{$\Delta_{[-]}: \Delta \to\Esp$}, where $\Esp$ is the category of topological spaces and continuous maps\footnote{We denote the category of topological spaces by the abbreviation $\Esp$ of the French word \emph{espace}, not just because there are other things called ${\sf Top}$ in the same contexts, but also as authors' reminder-to-self of the tacit sources of the approach \cite{GrothendieckA:SGA1,GrothendieckA:SGA4}.}, and $\Delta$ is the simplex category. Its objects are finite ordinals 
\bear
[m] & = &  \{0\lt1\lt 2\lt\cdots \lt m\}
\eear
while its morphisms are the order-preserving functions  \cite{Eilenberg-Zilber}.  All information about the simplicial structure of topological spaces is thus captured by the matrix
\bea\label{eq:simp-matrix}
\widehat \Delta \colon \Delta^o\times \Esp & \to & \Set\\
\protect \left[m\right]\times X & \mapsto & \Esp\left(\Delta_{[m]}, X\right)\notag
\eea
In the language of Sec.~\ref{Sec:FCA}, this is the \emph{``context matrix''}\/ of homotopy theory. Its \emph{``concept analysis''}\/ starts by making this matrix representable by adjoint functors along the lines of Sec.~\ref{Sec:repres}, instantiating (\ref{eq:deriv}--\ref{eq:galois-cat}) to $\Phi=\Delta$. But the matrix $\widehat \Delta$ is already representable by the functor $\Delta_{[-]}: \Delta \to\Esp$, and $\Esp$ is already complete. Leaving out of \eqref{eq:deriv} the completion of $\BBb=\Esp$, the adjunction in \eqref{eq:galois-cat} boils down to 
\beq\label{eq:galois-Kan}
\begin{tikzar}{}
\KKk\tto K\Delta \ar[mapsto]{dd} \& \sSet \ar[bend right=15]{dd}[swap]{\ladj \Delta}\ar[phantom]{dd}[description]\dashv \& 
\displaystyle \left(\Delta_{[-]}/ X\tto{\Dom} \Delta\right)
\\
\\
\displaystyle \supp\ \Big(\KKk\tto{K}\Delta \tto{\Delta_{[-]}} \Esp \Big)
\& \Esp \ar[bend right=15]{uu}[swap]{\radj \Delta}
\& X \ar[mapsto]{uu}
\end{tikzar}
\eeq
where the extension $\radj \Phi$ from \eqref{eq:galois-cat} becomes $\radj \Delta$ in \eqref{eq:galois-Kan} because
\bear
\inff\left(1\tto X \Esp \tto{_\bullet \Delta} \Dfib\diagup \Delta\right) & = & \left(\Delta_{[-]}/X\tto{\Dom} \Delta\right)
\eear
The adjunction $\left(\adj \Delta : \Esp \to \sSet\right)$ has been extensively studied. $\ladj \Delta: \sSet\to \Esp$ is the \emph{geometric realization}\/ of simplicial sets \cite{MilnorJ:geometric-realization}, whereas $\radj \Delta:\Esp \to \sSet$ is the \emph{singular decomposition}\/ of spaces. Eilenberg's singular homology was also based on it \cite{EilenbergS:singular}. Kan spelled out the concept of adjunction as the abstract structure underlying the concrete relationship between these two functors  \cite{KanD:adj,KanD:functors}. 
The idea was that the adjunction $\adj \Delta$ should provide a general method for transferring the invariants of interest between the geometric universe of spaces and an algebraic or combinatorial framework for computations. For a geometric realization $\ladj \Delta K\in \Esp$ of a complex $K\in\sSet$, the homotopy groups can indeed be computed combinatorially, from the structure of $K$ alone  \cite{KanD:combinatorial}, since the spaces in the form $\ladj \Delta K$ boil down to Whitehead's CW-complexes \cite{MilnorJ:geometric-realization,whitehead1949combinatorial}. But not all spaces are in this form.

\para{Trouble with localizations.} The upshot of Kan's adjunction $\adj \Delta\colon\Esp\to \sSet$ is that for any space $X$, we can construct a CW-complex $\rgt \Delta X = \ladj \Delta\radj\Delta X$, with a continuous map $\rgt \Delta X \tto\varepsilon X$, the counit of Kan's adjunction. This is the best approximation of $X$ by a CW-complex. When do such approximations preserve the invariants of interest? By the late 1950s, it was clear that combinatorial approximations work in many cases, certainly whenever $\varepsilon$ is invertible. In general, though,  even $\rgt \Delta \rgt \Delta X \tto \varepsilon \rgt \Delta X$ is not always invertible. The idea of approximating topological spaces by combinatorial complexes thus grew into a quest for making the units or the counits of adjunctions invertible. Which spaces have the same invariants as the geometric realizations of their singular\footnote{The word "singular" here means that the simplices, into which space may be decomposed, do not have to be embedded into it, which would make the decomposition \emph{regular}, but that the continuous maps from their geometric realizations may have \emph{singularities}.} decompositions? For particular invariants, there are direct answers \cite{Eilenberg-MacLane:relations,Eilenberg-MacLane:relationsTwo}. In general, though, the localizations at suitable spaces along suitable reflections or coreflections is the abstract form of spectral analysis and can be 
construed as an extension of concept analysis, as presented in Sec.~\ref{Sec:FCA}. Some of the most influential methods of algebraic topology can be thought of in these terms. Grossly oversimplifying, we mention three approaches.

The direct approach \cite[Vol. I, Ch.~5]{Gabriel-Zisman,BorceuxF:handbook} was to enlarge the given category by formal inverses of a family of arrows, usually called weak equivalences $\WWW$. They are made invertible in a category of fractions just like the non-zero integers are made invertible in the field of rationals. But adjoining formal inverses for a large class $\WWW$ in a large category like $\Esp$ usually leads to a category with large classes of morphisms between pairs of small objects.

Another approach \cite{dwyer2005homotopy,QuillenD:book} was to factor out the $\WWW$-arrows using two families of arrows, called fibrations $\FFF$ and cofibrations $\CCC$, such that the pairs of families $(\CCC\cap \WWW, \FFF)$ and $(\CCC, \WWW\cap\FFF)$ form weak factorization systems \cite[Part~III]{RiehlE:HT}. It follows that every morphism $f$ can be factored as in Fig.~\ref{Fig:quillen}. 
\begin{figure}[!ht]
\begin{center}
\begin{tikzar}[column sep=6em,row sep=9ex]
A\ar{d}[swap,pos=0.6]{\displaystyle f}\ar{r}[description]{\CCC\cap\WWW}  \& 
A^{\#}\ar{d}[pos=0.25,description]{\CCC\cap\FFF}[pos=0.6]{\displaystyle f^{\#}}
 \ar{dl}[pos=0.25,description]{\FFF}\\
B
\& B^{\#} \ar{l}[description]{\WWW\cap\FFF} \ar[leftarrow,crossing over]{ul}[pos=0.75,description]{\CCC}
\end{tikzar} 
\caption{A decomposition of $f$ by a model structure $(\CCC,\WWW,\FFF)$}
\label{Fig:quillen}
\end{center}
\end{figure}
Localizing every $f$ at the corresponding $\CCC\cap\FFF$-component $f^{\#}$ makes the weak equivalences $\WWW$ invertible. The triples $(\CCC,\WWW,\FFF)$ are called  \emph{model structures}\/ since they make categories into abstract models of homotopy theories. The general idea originates from the model structures in $\Esp$ and $\sSet$, making the adjunction in \eqref{eq:galois-Kan} into an equivalence. This was Dan Quillen's approach \cite{QuillenD:rational,QuillenD:book}. Here we worked out a model structure on the category of adjunctions.  

The third approach  \cite{Applegate-Tierney:models,Applegate-Tierney:iterated} tackles the task of making the arrows $\rgt \Delta X \tto\varepsilon X$ invertible by modifying the comonad $\rgt \Delta$ until it becomes idempotent and then localizing at the coalgebras of this idempotent comonad. Note that this approach does not tamper with the continuous maps in $\Esp$, be it to make some of them formally invertible, or to factor them out. The idea is that a comonad $\rgt \Delta_\infty:\Esp \to \Esp$, should localize any space $X$ at a space $\rgt \Delta_\infty X$ such that $\rgt \Delta\rgt \Delta_\infty X \stackrel \varepsilon \cong \rgt \Delta_\infty X$. It follows that $\rgt \Delta_\infty\rgt \Delta_\infty X \stackrel \varepsilon \cong \rgt \Delta_\infty X$ must hold and the comonad $\rgt \Delta_\infty$ must be idempotent. The quest for such a comonad is illustrated in Fig.~\ref{Fig:idempot}
\begin{figure}[!ht]
\[\begin{tikzar}[row sep=2.5cm,column sep=3cm]
\sSet 
\arrow[bend right = 12]{d}
\arrow[bend left = 12,leftarrow]{d}
\arrow[phantom]{d}[description]{\scriptstyle \adj\Delta}  
\arrow[bend right = 7,thin]{dr}
\arrow[bend left = 7,leftarrow,thin]{dr}
\arrow[bend right = 4,thin]{drr}
\arrow[bend left = 4,leftarrow,thin]{drr}
\arrow[phantom]{dr}[description]{\scriptstyle\Delta^0\dashv\Delta_0}
\arrow[phantom]{drr}[description]{\scriptstyle\Delta^1\dashv\Delta_1}
\arrow[bend right = 2,thin,no head]{drrr}
\arrow[bend left = 2,thin,no head]{drrr}
\arrow[phantom]{drrr}[description]{\scriptstyle\Delta^\alpha\dashv\Delta_\alpha}
\\
\Esp  \arrow[loop, out = -50, in=-130, looseness = 6]{}[swap]{\rgt \Delta} 
\arrow[bend right = 9]{r}
\arrow[bend left = 9,leftarrow]{r}
\arrow[phantom]{r}[description]{\upadj V}
\& 
\Emc \Esp \Delta 
\arrow[phantom]{r}[description,rotate=-90]{\dashv}
\arrow[loop, out = -50, in=-130, looseness = 6]{}[swap]{\rgt \Delta_0} 
\arrow[bend right = 9]{r}
\arrow[bend left = 9,leftarrow]{r}
\arrow[phantom]{r}[description]{\upadj V}
\& 
\arrow[phantom]{r}[description,rotate=-90]{\dashv}
\left(\Emc\Esp \Delta\right)^{\rgt \Delta_0}
\arrow[loop, out = -60, in=-120, looseness = 5.5]{}[swap]{\rgt \Delta_1}
\arrow[bend right = 9,dotted]{r} 
\&
\cdots \arrow[bend right = 9,dotted]{l} 
\arrow[loop, out = -55, in=-125, looseness = 12]{}[swap]{\rgt \Delta_\alpha}
\end{tikzar}\]
\caption{Iterating the comonad resolutions for $\rgt \Delta$}
\label{Fig:idempot}
\end{figure}
where $\Emc \Esp \Delta$ is the category of $\rgt\Delta$-coalgebras, $\adj V:\Esp\to \Emc \Esp \Delta$ is its final resolution, and $\Delta^0$ is the comparison functor, mapping a complex $K$ to the coalgebra $\ladj \Delta K \tto{\ladj\eta} \ladj\Delta\radj\Delta\ladj \Delta K$. Since $\sSet$ is a complete category, $\Delta^0$ has a right adjoint $\Delta_0$, and they induce the comonad $\rgt\Delta_0 = \Delta^0\Delta_0$ on $\Emc \Esp\Delta$. If $\rgt \Delta$ was idempotent, the final resolution $\adj V$ would be a coreflection, and the comonad $\rgt \Delta_0$ would be isomorphic to the identity. But $\rgt \Delta$ is not idempotent, $\rgt \Delta_0$ is not the identity, and its final resolution, generically denoted $\adj V$ again, induces a comparison functor $\Delta_1$, which induces the comonad $\rgt\Delta_1 = \Delta^1\Delta_1$, which induces the comonad $\rgt\Delta_2$, and so on. The tower of final resolutions of these monads, displayed in Fig.~\ref{Fig:idempot} $\adj V$,  continues infinitely and produces coalgebras over coalgebras over coalgebras\ldots, which span inductive towers in $\Esp$. At each step, the parts of a space that are not approximable by geometric realizations are projected away, but because of that, some previously approximable parts cease to be approximable. This keeps occurring after infinitely many steps and pushes the inductive tower into ordinal lengths, echoing Cantor's quest for accumulation points of the convergence domains of the Fourier series, which led him to discover ordinals. Since $\Esp$ is cocomplete, the process eventually stabilizes, and the resulting idempotent comonad $\rgt \Delta_\infty$ localizes at the spaces that can be analyzed combinatorially, in terms of simplicial sets alone, using Kan's adjunction $\adj \Delta$.
\begin{figure}[!ht]
\[\begin{tikzar}[row sep=2.5cm,column sep=3cm]
\sSet 
\arrow[loop, out = 135, in = 45, looseness = 4]{}[swap]{\lft \Delta} 
\arrow[bend right = 9]{d}
\arrow[bend left = 9,leftarrow]{d}
\arrow[phantom]{d}[description]{\adj \Delta}  
\arrow[bend right = 5]{r}
\arrow[bend left = 5,leftarrow]{r}
\arrow[phantom]{r}[description]{\scriptstyle\Delta^0\dashv \Delta_0}
\& 
\Emc\Esp \Delta \arrow[loop, out = 135, in = 45, looseness = 2.5]{}[swap]{\Lft \Delta}
\arrow[bend right = 9]{d}
\arrow[bend left = 9,leftarrow]{d}
\arrow[phantom]{d}[description]{\nadj \Delta}     
\arrow[leftrightarrow]{r}[description]{\mbox{\Huge$\simeq$}} 
\& \Ec \sSet \Delta 
\arrow[bend right = 9]{d}
\arrow[bend left = 9,leftarrow]{d}
\arrow[phantom]{d}[description]{\nadj \Delta}
\\
\Esp  \arrow[loop, out = -45, in=-135, looseness = 6]{}[swap]{\rgt \Delta} 
\arrow[bend right = 5]{ur}
\arrow[bend left = 5,leftarrow]{ur}
\arrow[phantom]{ur}[description]{\scriptstyle\adj V}
\arrow[bend right = 5]{r}
\arrow[bend left = 5,leftarrow]{r}
\arrow[phantom]{r}[description]{\scriptstyle H^1\dashv H_1}
\& 
\Emm \sSet \Delta 
\arrow[loop, out = -45, in=-135, looseness = 6]{}[swap]{\Rgt \Delta} 
\arrow[leftrightarrow]{r}[description]{\mbox{\Huge$\simeq$}} 
\& \Em \Esp \Delta 
\end{tikzar}\]
\caption{The nucleus of the Kan adjunction}
\label{Fig:kan-nuc}
\end{figure}

\para{The nuclear solution.} The nucleus of Kan's adjunction, displayed in Fig.~\ref{Fig:kan-nuc}, allows analyzing a much larger family of spaces combinatorially, also in terms of simplicial sets, but equipped with $\lft \Delta$-algebra structures for the monad $\lft \Delta = \radj\Delta\ladj \Delta$. This larger family of spaces arises already at the first step in Fig.~\ref{Fig:idempot}, without going down the tower. They are the spaces that can be recovered as retracts of geometric realizations of their singular decompositions. The embedding parts of the  retractions make them into $\rgt\Delta$-coalgebras in $\Emc \Esp\Delta$. Thm.~\ref{thm:Theorem} allows presenting such spaces as coalgebras over the simplicial algebras in $\Emm \sSet \Delta$. If that presentation, albeit combinatorial, seems complicated, note that Prop.~\ref{Prop:three} allows simplifying them, and provides an effective equivalence between the spaces with coalgebras in $\Emc \Esp\Delta$ and the simplicial sets with idempotents in $\sSet^{\Rgt \Delta}$.

\section{Closing}
\label{Sec:what}

\para{Categorical dimension reduction.} The task of dimension reduction often arises in mathematics: manifolds are reduced to their boundaries, curves are approximated by points, linear operators are diagonalized using spectra. Spectral methods are the tenet of concept mining and principal component analysis of data \cite{Azar}. In category theory, on the other hand, the way up into higher categories has been more attractive than the way down, for a variety of reasons. One is that categories themselves expand into the higher dimensions \cite{BaezJ-May:higher,LeinsterT:operads,LurieJ:higher,SimpsonC:higher}. Another is that Grothendieck's \emph{homotopy hypothesis}\/ \cite{GrothendieckA:pursuing,pursuing-asterisque} made higher category theory into a geometric pursuit, by some accounts subsuming homotopy theory. The nucleus results seem to suggest that there might be a way down even in category theory. On the object level, the nucleus construction ostensibly follows the familiar procedures of spectral decomposition. The nuclear adjunction induced by a categorical matrix can be construed as a diagonal matrix, spanned by the eigenvectors of the original matrix. On the meta-level, the fact that the categories of categories needed for concept analysis do not form higher categories but interact with each other within their own dimensionality seems unusual. It may signal a shortcoming of the current approach. Or it may signal the robustness and stability of categories as a mathematical tool for the application areas where the way up is the way down, like statistics, machine learning, and data analysis; the areas where the problems get simpler as we go, and not more complicated.

\para{Summary.}
We saw adjunctions in action in recommender systems in Sec.~\ref{Sec:FCA}, and in algebraic geometry and homotopy theory in Sec.~\ref{Sec:HT}. In Sections~\ref{Sec:cat}--\ref{Sec:simple}, we derived a common denominator: the nuclear adjunctions. It is a powerful tool of mathematical conceptualizations; in a certain sense universal. In the meantime, it proved its worth by untying (as inadvertently as it came about) a  fundamental problem, posed by Jim Lambek in 1966: \emph{to complete any given category under limits and colimits in such a way that any limits and colimits that may already exist are preserved} \cite{LambekJ:completions}. We call such completions \emph{tight}. Dedekind constructed the reals as a tight completion of the rationals. MacNeille generalized Dedekind's construction to posets. The task was thus to lift the Dedekind-MacNeille completions from posets to categories. It lifts readily to poset-enriched categories, with no major surprises, yet  with significant and varied applications \cite{PavlovicD:ICFCA12,PavlovicD:Samson13,PavlovicD:LICS17,WillertonS:nucleus}. Yet the ordinary, $\Set$-enriched categories offer resistance to being tightly completed under all limits and colimits, for fundamental reasons spelled out by John Isbell \cite[Thm. 3.1]{IsbellJ:no-lambek}. On the other hand, the concept analysis requires nothing less than tight completions: completing the data contexts from above and from below at once, while preserving the previous completions. And in recommender systems, avoiding blind amplifications and echoing everything requires going beyond the posetal and numeric models. So be it as it may, Isbell notwithstanding, we cannot avoid adding some limits and colimits tightly. Which ones? This is where the nucleus kicked in. The idea was sketched in \cite{PavlovicD:CALCO15}. It has been worked out in several forms, and an attempt to present it was made in \cite{PavlovicD:dede}. 

\para{Back to Escher.}
The nucleus of an adjunction, illustrated in the Introduction  by Escher's ``Drawing hands'', is a descendant of Dedekind's construction of the real numbers as \emph{cuts}\/ of the rationals. A cut is a pair $c = <\lft c, \rgt c>\in \Do\QQq\times \Up\QQq$ where $\lft c$ is the set of lower bounds of $\rgt c$, and $\rgt c$ is the set of upper bounds of $\lft c$. The ``Drawing hands'' of Dedekind's cuts are displayed in Fig.~\ref{Fig:dede-escher} on the left.
\begin{figure}[!ht]
\begin{center}
$\begin{tikzar}[row sep=1.5mm,column sep=1.5mm]
\&
\Do \QQq 
\arrow[phantom]{dddddddd}[description]{{\dashv}}   
\arrow[bend right = 13]{dddddddd}[swap]{\underline \cmn}  
\\ 
\lft c \ar[mapsto,bend right = 13]{dddddd} 
\&\& 
\lft c
\\ \hspace{3em}\\ \hspace{3em} \\ \hspace{3em} \\ \\ \\ 
\rgt c \&\& \rgt c \ar[mapsto,bend right = 13]{uuuuuu}
\\ 
\&
\Up\QQq
\arrow[bend right = 13]{uuuuuuuu}[swap]{\overline\mnd}  
\end{tikzar}$
\qquad\qquad\qquad\qquad\qquad
$\begin{tikzar}[row sep=1.5mm,column sep=1.5mm]
\&
\Ec\AAa F 
\arrow[phantom]{dddddddd}[description]{{\dashv}}   
\arrow[bend right = 13]{dddddddd}[swap]{\loadj{F}}  
\\ 
<x, \alpha_x> \ar[mapsto,bend right = 13]{dddddd} 
\&\& 
<\radj F u, \ladj F \beta^u>
\\ \hspace{3em}\\ \hspace{3em} \\ \hspace{3em} \\ \\ \\ 
<\ladj F x, \radj F \alpha_x  >\&\& <u,\beta^u> \ar[mapsto,bend right = 13]{uuuuuu}
\\ 
\&
\Em \BBb F 
\arrow[bend right = 13]{uuuuuuuu}[swap]{\roadj{F}}  
\end{tikzar}$
\caption{The nucleus of an adjunction echoes the idea of Dedekind cuts}
\label{Fig:dede-escher}
\end{center}
\end{figure}
The ``Drawing hands'' of a simple nucleus are displayed on the right. Since the components $\lft c$ and $\rgt c$ of a cut determine each other, they both carry  the same information, and either can be used alone to represent a real number. The real numbers have three isomorphic representations:
\[
\left\{\lft c \in \Do\QQq\ |\ \lft c = \overline \mnd \underline \cmn \lft c \right\} \ \ \ \ \cong\ \ \ \ 
\left\{\left<\lft c ,\rgt c\right>\in \Do\QQq\times \Up \QQq\ |\ \lft c = \overline \mnd \rgt c \ \wedge \ \underline\cmn \lft c = \rgt c \right\} \ \ \ \ \cong\ \ \ \ 
\left\{\rgt c \in \Up\QQq\ |\ \underline \cmn \overline \mnd \rgt c = \rgt c\right\}
\]
This just says that the same real number can be approached from below, from above, or tightly, from both directions. The nucleus construction shows that concepts are like real numbers: a concept can be presented as an algebra, or as a coalgebra, and they determine each other. And it shows that concepts are unlike real numbers since the algebra and the coalgebra are not each other's images, but mere retracts. As a nucleus of an adjunction, an instance of itself, the nucleus construction seems easier to use than to understand. The path forward might be to conceptualize it in its own terms. Part II of this work will describe the nucleus monad on monads and comonads, covering itself. Part III will spell out the induced model structures and localizations.

\paragraph*{Acknowledgements and apology.}
While working on the results presented in this paper, the first author was assisted by Liang-Ting Chen, Tobias Heindel, Toshiki Kataoka, Peter-Michael Seidel, and Christina Vassilakopoulou. Toshiki coauthored  \cite{PavlovicD:CALCO15}, where an early version of Prop.~\ref{prop:nuc} was proved. Peter-Michael coauthored \cite{PavlovicD:LICS17}, where parts of the proof of Thm.~\ref{thm:Theorem} were constructed. A longer version of the present paper, covering a wider family of examples and the results now left for parts II and III, has been circulated since 2019 and posted on arxiv in 2020. It also contained an account of the nucleus monads on monads and comonads, induced by the nucleus monad on adjunctions presented here. The dualities  between monads and comonads in abstract 2-categories, studied by Ross Street in the final sections of \cite{StreetR:monads}, induce closely related 2-monads directly, without a detour into adjunctions. Ross presented the result concerning the 2-monad on monads in his seminar, provided helpful comments, and considered the mentioned previous version for publication in a journal that he recommended. After a couple of years, the anonymous referees returned very succinct reviews, dispensing pedagogic style advice as a conditional acceptance recommendation. In view of the absence of substantive comments, suggesting an apparent lack of interest, we withdrew the submission, but realized that a part of the problem was our weaving and knotting of several narratives into one, each of the coauthors trusting that the other one would assure a tidy presentation. We apologize for the confusion and thank the readers and the users seeing us through it.

\bibliographystyle{alphaurl}
\bibliography{CT,HT,logic,math,PavlovicD,ref-nuc-lmcs,type}

\appendix
\section{Comprehending presheaves and matrices as discrete fibrations}\label{Appendix:compreh}

Following the step from \eqref{eq:matrp} to \eqref{eq:compreh-pos}, the comprehension correspondence \eqref{eq:compreh-pos} now lifts to 
\bea\label{eq:compreh-cat-appendix}
\Cat (\AAa^\op \times \BBb , \Set) &
\begin{tikzar}[row sep = 4em]\hspace{.1ex} \ar[bend left]{r}{\eh{(-)}} \ar[phantom]{r}[description]{\simeq} 
\& \hspace{.1ex}
\ar[bend left]{l}{\Xi}
\end{tikzar} & \Dfib\diagup \AAa\times \BBb^\op \\
\left(\AAa^\op \times \BBb\tto{\Phi} \Set \right) & \mapsto & \left(\tint\Phi\tto{\eh \Phi} \AAa \times \BBb^\op \right)
\notag \\ \notag
\left(\AAa^\op \times \BBb\tto{\Xi_E} \Set \right)
&\mapsfrom & \left(\EEe\tto E \AAa\times \BBb^\op\right)  
\eea
Transposing the arrow part of $\Phi$, which maps every pair $f\in \AAa(a,a')$ and $g\in \BBb(b',b)$ into $\Phi(a',b')\tto{\Phi_{fg}} \Phi(a,b)$, the closure property expressed by the implication in \eqref{eq:monotone} becomes the mapping
\bea
\AAa(a,a') \times \Phi(a', b') \times \BBb(b', b) & \tto{\hspace{1.5em}} & \Phi(a,b)
\eea 
The \emph{lower-upper}\/ closure property expressed by \eqref{eq:monotone} is now captured as the structure of the total category $\tint \Phi$, defined as follows:
\bea\label{eq:tint}
\left|\tint \Phi \right| & = & \coprod_{\substack{a\in \AAa\\ b\in \BBb}} \Phi(a,b)\\
\tint\Phi\left(x_{ab},x'_{a'b'}\right) & = & \left\{<f,g>\in \AAa(a,a')\times \BBb(b',b)\ |\ x = \Phi_{fg}(x') \right\}\notag
\eea
It is easy to see that the obvious projection 
\bea\label{eq:tintPhi}
\tint \Phi & \tto{\eh\Phi} & \AAa\times \BBb^\op\\
x_{ab} &\mapsto & <a,b>\notag
\eea
is a discrete fibration, i.e., an object of $\Dfib \diagup \AAa\times \BBb^\op$. In general, a functor $\FFf \tto F \CCc$ is a discrete fibration over $\CCc$ when for all $x\in \FFf$ the obvious induced functors $\FFf/x \tto{F_x} \CCc/Fx$ are isomorphisms. In other words, for every $x\in \FFf$ and every morphism $c \tto{t} Fx$ in $\CCc$, there is a unique lifting $t^! x\tto{\vartheta^t} x$ of $t$ to $\FFf$, i.e., a unique $\FFf$-morphism into $x$ such that $F(\theta^t) = t$. For a discrete fibration $\EEe\tto E \AAa\times \BBb^\op$, such liftings induce the arrow part of the corresponding presheaf
\bear
\Xi_E \colon \AAa^\op \times B & \to & \Set\\
<a,b> & \mapsto & \{x\in \EEe\ |\ Ex =<a,b>\}
\eear
because any pair of morphisms $<f,g>\in \AAa(a,a')\times \BBb^\op(b,b')$ lifts to a function $\Xi_E(f,g) = <f,g>^! :\Xi_E(a',b') \to \Xi_E(a,b)$. The equivalences in \eqref{eq:compreh-cat-appendix} thus yield an equivalent version of the category $\Mat$ of matrices:
\bea\label{eq:Mat-appendix}
|\Mat| & = & \coprod_{\AAa, \BBb\in \CAT} \Dfib\diagup \AAa \times \BBb^\op\\[2ex] 
\Mat (\Phi, \Psi) & = & \coprod_{\substack{H\in \CAT(\AAa,\CCc)\\
K\in \CAT(\BBb,\DDd)}} \Bigg(\Dfib\diagup \AAa \times \BBb^\op\Bigg)\Bigg(\Phi,(H\times K^\op)^\ast\, \Psi\Bigg)\notag
\eea
where $\Phi\in \Dfib\diagup \AAa\times \BBb^\op$, $\Psi\in \Dfib\diagup \CCc\times \DDd^\op$, and $(H\times K^\op)^\ast\, \Psi$ is a pullback of $\Psi$ along 
\bear
\left(H\times K^\op\right)\ \colon\  \AAa\times \BBb^\op &\tto{\hspace{1.5em}} & \CCc\times \DDd^\op\eear 
The notational abuse of $\Mat$ to denote both \eqref{eq:compreh-pos} and \eqref{eq:Mat-appendix} is not just technically harmless, but its tacit identification of the two sides of \eqref{eq:compreh-cat-appendix} is conceptually justified by the recurring theme of categorical comprehension, used already in the next section.

\section{Idempotents, absolute (Cauchy) completions, and weak (Morita) equivalences}\label{appendix:idempot}

\para{Idempotent splitting.} An endomorphism $\varphi:A\to A$ is \emph{idempotent}\/ if $\varphi\circ\varphi = \varphi$. A \emph{splitting} of an idempotent $\varphi$ is its epi-mono (quotient-injective) factorization
\[
\begin{tikzar}{}
A \arrow{rr}{\varphi} \arrow[twoheadrightarrow]{dr}[swap]{q} \&\& A\\ \& B \arrow[rightarrowtail]{ur}[swap]{i} 
\end{tikzar}
\]
It is easy to prove (and goes back at least as far as \cite[Sec.~IV.7.5]{GrothendieckA:SGA4}) that the following statements are equivalent:
\begin{enumerate}[(a)]
\item $\varphi\circ \varphi = \varphi$
\item $q\circ i = \id$
\item $i$ is an equalizer and $q$ is a coequalizer of $\varphi$ and the identity
\[\begin{tikzar}{}
\& B \ar[tail]{dl}[swap]{i} 
\\ A \ar[shift left = .75ex]{rr}{\varphi} \ar[shift right = .75ex]{rr}[swap]{\id} \&\& A \ar[two heads]{ul}[swap]{q} 
\end{tikzar}
\]
\end{enumerate}
The idempotent splittings are often drawn in the form $\begin{tikzar}{}
A \arrow[bend left = 15,twoheadrightarrow]{r}{q}  \&  B \arrow[bend left = 15,rightarrowtail]{l}{i} 
\end{tikzar}$, suggesting $i\circ q = \varphi$ and  $q\circ i = \id$. Since idempotent splittings are equationally characterized, and functors preserve equations, the idempotent splittings are preserved by all functors. A categorical property that is preserved by all functors is called \emph{absolute}. By (c) above, the idempotent splittings are absolute equalizers and coequalizers. It was shown in  \cite{PareR:absolute} that these are the only absolute limits or colimits. Completing a category under all absolute limits and colimits boils down to adjoining the idempotent splittings.

\para{Absolute and Cauchy completions.} 
For an arbitrary category $\CCc$, the absolute completion $\Kar\CCc$ consists of the $\CCc$-idempotents as the objects, and the idempotent-preserving homomorphisms:
\bea\label{eq:karoubi}
\rvert\Kar \CCc\rvert & = & \coprod_{A\in |\CCc|} \{A\tto\varphi A\ |\ \varphi\circ \varphi = \varphi\} \\
\Kar\CCc (A\tto\varphi A, B\tto\psi B) & = & \left\{f\in \CCc(A,B)\ \Big|\ \ \ \begin{tikzar}[row sep=1.8em,column sep=1.8em]
 A \ar{r}{f} \ar{d}{\varphi} \& B  \\ 
 A \ar{r}{f}\& B \ar{u}{\psi}
\end{tikzar}\ \right\}
\eea
Note that the homomorphism condition $\psi \circ f\circ \varphi = f$ is equivalent to the conjunction of $f\circ \varphi = f$ and $\psi \circ f = f$. In \cite[Sec.~IV.7.5]{GrothendieckA:SGA4}, the construction of the category $\Kar \CCc$ was attributed to Max Karoubi, so it came to be called the \emph{Karoubi envelope}. It also appeared as exercise 2--B in \cite[p.~61]{FreydPJ:abelian}. 

A category $\CCc$ is said to be \emph{Cauchy complete}\/ if every matrix $\CCc^\op\tto{\Gamma} \Set$ with a right adjoint $\CCc\tto{\radj \Gamma} \Set$, where $\radj \Gamma(x)$ is the set of cocones from $\Gamma$ to $x$, is representable by some $c\in \CCc$ as $\Gamma(x) = \CCc(x,c)$. The name is motivated by the observation that the corresponding property in metric spaces, viewed as enriched categories, characterizes the convergence of Cauchy sequences \cite{LawvereFW:metric}. See also \cite[Vol.~1, Sec.~7.9]{BorceuxF:handbook} and \cite{BungeM:tight}.  The following statements are equivalent for any category $\CCc$:
\begin{enumerate}[(a)]
\item $\CCc$ is Cauchy complete,
\item $\CCc$ is absolutely complete,
\item all idempotents split in $\CCc$,
\item the embedding $\CCc \inclusion \Kar \CCc$ is an equivalence of categories.
\end{enumerate}

\para{Weak (Morita) equivalences.} Categories $\CCc$ and $\DDd$ are said to be \emph{weakly equivalent}\/  when their absolute completions are strongly equivalent; i.e., $\CCc \sim \DDd$ means that $\Kar \CCc \simeq \Kar \DDd$. The weak equivalence is often named the  \emph{Morita equivalence of categories}\/ because it is also characterized by the strong equivalence of the categories $\Set^\CCc\simeq \Set^\DDd$. A  proof can be found in  \cite[Vol.~1, Sec.~7.9]{BorceuxF:handbook}. This terminology is justified by the analogy of the categories $\Set^\CCc$ and $\Set^\DDd$ with the abelian categories ${\Ab}^R$ and  ${\Ab}^S$ of $R$-modules and $S$-modules, whose equivalence was studied by Morita.\footnote{Neither Morita nor Cauchy were involved with the categorical liftings of the concepts that carry their names. History is often used for mnemonic purposes.} All of our results are valid up to weak equivalence of categories or can be construed as speaking of absolute completions. To shorten proofs, we take the latter approach, introduce absolute completions in the statements, and prove strong equivalences. Omitting this would lead to shorter statements, longer proofs, and essentially equivalent theory. 

\para{Absolute reflections.} Let $\Kar \Cat$ denote the full subcategory of $\Cat$ spanned by absolutely complete categories $\Kar \AAa$, $\Kar \BBb$, etc. Checking that the idempotent splitting construction in \eqref{eq:karoubi} induces a left adjoint to the inclusion $\Kar \Cat\inclusion \Cat$ is a standard exercise. Any functor $F:\AAa\to \BBb$ lifts to $\Kar F:\Kar \AAa \to \Kar \BBb$, where $\Kar F x = Fx$, whether $x$ is an idempotent or a homomorphism of idempotents. Any natural transformation $\tau: F\to G$ between any pair of functors $F, G:\AAa\to \BBb$ lifts to a natural family $\Kar \tau: \Kar F \to\Kar G$ comprised of
\bea\label{eq:natlift}
\Big(\Kar F \varphi \tto{\ \Kar\tau_\varphi\ } \Kar G \varphi\Big) & = & \Big(FA\tto{ F\varphi} FA \tto{\tau_A} GA\tto{G\varphi} GA
\Big)
\eea
The naturality of $\tau$ implies not only the naturality of the family $\Kar \tau$ but also the functoriality of the induced mappings $\Cat(\AAa,\BBb) \to \Cat(\Kar \AAa, \Kar \BBb)$. The idempotent splitting in \eqref{eq:karoubi} thus induces a 2-functor $\Cat\to \Kar \Cat$. This means that any monad $\left(\lft T, \eta, \mu\right)$ on $\AAa$ lifts to a monad $\left(\Kar{\lft T}, \Kar \eta, \Kar \mu\right)$ on $\Kar \AAa$. Ditto for any comonad, adjunction, matrix, etc. Hence the corresponding 2-categories $\Kar\Mnd$, $\Kar \Cmn$, $\Kar \Adj$, etc., of the various categorical structures lifted to absolute completions, in all cases fully embedded into the general categorical structures. 

\para{Notation.} Since the embedding $\Kar \Cat\inclusion\Cat$ is conservative, we write $\adj F: \Kar \BBb \to \Kar \AAa$ instead of $\adja F: \Kar \BBb\to \Kar \AAa$, since the latter provides no additional information.

\end{document}